%% file: specialbirat-10-2015.tex
\documentclass[11pt]{article} 

\pdfoutput=1

\usepackage{amsmath,amscd,amssymb}
\usepackage[applemac]{inputenc}
\usepackage[english]{babel}
\usepackage{graphicx}
\begin{document}
\title{Special birational structures on non-Kählerian complex surfaces}
\author{Georges Dloussky\\
Aix Marseille Université, CNRS, Centrale Marseille, I2M, UMR 7373,\\
39, rue F.Joliot-Curie, 13453 Marseille, France\\
  georges.dloussky@univ-amu.fr
}  
\date{ } 
\newtheorem{Lem}{Lemma \thesection.}
\newtheorem{Th}[Lem]{Theorem \thesection.}
\newtheorem{Th and Def}[Lem]{Theorem and Definition \thesection.}
\newtheorem{Cor}[Lem]{Corollary \thesection.}
\newtheorem{Def}[Lem]{Definition \thesection.}
\newtheorem{Ex}[Lem]{Example \thesection.}
\newtheorem{Exs}[Lem]{Examples \thesection.}
\newtheorem{Prob}[Lem]{Problem \thesection.}
\newtheorem{Prop}[Lem]{Proposition \thesection.}
\newtheorem{Prop and Def}[Lem]{Proposition and Definition \thesection.}
\newtheorem{Rem}[Lem]{Remark \thesection.}
\newtheorem{Lem and Def}[Lem]{Lemma and Definition \thesection.}
\newtheorem{Exer}[Lem]{Exercise \thesection.}
\newtheorem{Not}[Lem]{Notations \thesection.}
\setcounter{Lem}{0}
\setcounter{section}{-1}
\def\cal{\mathcal}
\def\frak{\mathfrak}
\def\bb{\mathbb} 
\def\a{\alpha }
\def\b{\beta }
\def\d{\delta }
\def\D{\Delta }
\def\g{\gamma }
\def\G{\Gamma }
\def\d{\delta }
\def\D{\Delta }
\def\e{\epsilon }
\def\z{\zeta }
\def\k{\kappa }
\def\l{\lambda }
\def\L{\Lambda }
\def\m{\mu }
\def\o{\omega }
\def\p{\pi }
\def\P{\Pi }
\def\s{\sigma }
\def\t{\theta }
\def\T{\Theta }
\def\f{\varphi }
\def\Card{{\rm Card\, }}
\def\ker{{\rm Ker\,}}
\def\im{{\rm Im\,}}
\def\coker{{\rm Coker\,}}
\def\codim{{\rm codim\,}}
\def\part{\partial }
\def\dps{\displaystyle }
\def\ot{\otimes }
\def\tr{{\rm tr\,}}
\def\iff{ if and only if }
\def\wh{\widehat }
\def\wt{\widetilde }
\def\la{\langle }
\def\ra{\rangle }

\input{diagram} 
\setcounter{section}{0}
\noindent
\maketitle
\abstract{We investigate the following conjecture: all compact non-Kählerian complex surfaces admit birational structures. After Inoue-Kobayashi-Ochiai, the remaining cases to study are essentially surfaces in class $VII_0^+$. We show that Kato surfaces with a cycle and only one branch of rational curves admit a special birational structure given by new normal forms of contracting germs in Cremona group $Bir(\bb P^2(\bb C))$. In particular all surfaces $S$ with GSS and second Betti number satisfying $0\le b_2(S)\le 3$ admit a birational structure. From the existence of a special birational structure we deduce a developing meromorphic mappings $\tilde S\to \bb P^2(\bb C)$ from the universal cover of $S$ to $\bb P^2(\bb C)$ which blows down an infinite number of rational curves. From this mapping we recover a GSS.}\\
\begin{center}\bf Résumé
\end{center}
On étudie la conjecture suivante: toute surface complexe compacte non kählerienne admet une structure birationnelle. D'après Inoue-Kobayashi-Ochiai, les cas restant à étudier sont essentiellement les surfaces de la classe $VII_0^+$. On démontre que les surfaces de Kato qui ont un cycle avec un seul arbre  de courbes rationnelles admettent une structure birationnelle spéciale définie par de nouvelles formes normales de germes de contractions dans le groupe de Cremona $Bir(\bb P^2(\bb C))$. En particulier toute surface $S$ contenant une coquille sphérique globale (CSG) et pour laquelle le second nombre de Betti vérifie $0\le b_2(S)\le 3$ admet une structure birationnelle. De l'existence d'une structure birationnelle on déduit une application méromorphe développante $\tilde S\to \bb P^2(\bb C)$ du revêtement universel de $S$ dans $\bb P^2(\bb C)$ qui écrase une infinité de courbes rationnelles. Cette application permet de reconstituer la CSG.\\


\section{Introduction}
What is the best $G$-structure on a compact manifold ?  The classification of Inoue-Kobayashi-Ochiai \cite{IKO,KO} shows that all compact complex non-Kählerian surfaces but some Hopf surfaces and surfaces in class VII$_0^+$ (i.e. in class VII$_0$ with $b_2>0$) admit affine structures. In view of the explicit construction of Kato surfaces (i.e. minimal surfaces $S$ containing a global spherical shell (GSS) with $b_2(S)>0$) and the particular cases of Enoki surfaces and Inoue-Hirzebruch surfaces the best $G$-structure should be obtained for a subgroup of $Bir(\bb P^2(\bb C))$. It justifies the following\\
{\bf Conjecture}: All compact complex non-Kähler surfaces admit  birational structures.\\
The conjecture is clearly satisfied for all Hopf surfaces because they are defined by an invertible contracting polynomial mapping. Remains the case of surfaces in class VII$_0^+$. Since the only known surfaces in class VII$_0^+$  are Kato surfaces and since it is conjectured that there are no others, this article focuses on the following problem: 
Do compact surfaces with GSS admit a birational structure, i.e. is there an atlas with transition mappings in Cremona group $Bir(\bb P^2(\bb C))$. As stronger requirement, is there in each conjugation class of contracting germs of the form $\Pi\s$ (or of strict germs, following Favre terminology \cite{Fav}) a birational representative~? Clearly $\Pi\s$ is birational if and only if $\s$ is birational.\\

{\bf Known results}: 
\begin{itemize}
\item If $S$ is a Enoki surface (see \cite{DK}) or a Inoue-Hirzebruch surface (see \cite{D2})  with second Betti number $b_2(S)=n$, known normal forms, namely
$$F(z_1,z_2)=(t^nz_1z_2^n+\sum_{i=0}^{n-1}a_it^{i+1}z_2^{i+1}, tz_2),\quad 0<|t|<1,$$
and
$$N(z_1,z_2)=(z_1^pz_2^q,z_1^rz_2^s),$$
respectively,  are birational. Here  $\left(\begin{array}{cc}p&q\\r&s\end{array}\right)\in Gl(2,\bb Z)$ is the composition of $n$ matrices 
$$\left(\begin{array}{cc}1&1\\0&1\end{array}\right) \quad {\rm or}\quad \left(\begin{array}{cc}0&1\\1&1\end{array}\right)$$
 with at least one of the second type.
\item If $S$ is of intermediate type (see definition in section 2), there are normal forms due to C.Favre \cite{Fav} 
$$F(z_1,z_2)=(\l z_1z_2^{\frak s}+P(z_2),z_2^k), \qquad \l\in\bb C^\star,\ \frak s\in\bb N^\star,\  k\ge 2,$$
where $P$ is a special polynomial. These normal forms are adapted to logarithmic deformations and show the existence of a foliation, however {\it are not birational}.  In \cite{OT} K.Oeljeklaus and M.Toma explain how to recover second Betti number $n$ which is now hidden and give coarse moduli spaces of surfaces with fixed intersection matrix,
\item Some special cases of intermediate surfaces are obtained from Hénon mappings $H$ or composition of Hénon mappings. More precisely,  the germ of $H$ at the fixed point at infinity is strict, hence gives a surface with a GSS \cite{HO, DO2}. These germs are birational.
\end{itemize}
\vspace{2mm}
{\bf Motivation}:\\
Let $S$ be a minimal compact complex surface with Betti numbers $b_1(S)=1$, $n=b_2(S)>0$, the class of such surfaces will be denoted VII$_0^+$. We consider the following conditions:
\begin{description}
\item{(A)} $S$ contains a global spherical shell (GSS),
\item{(B)} $S$ contains $b_2(S)$ rational curves,
\item{(C)} $S$ contains a cycle of rational curves,
\item{(D)} $S$ admits a deformation into $b_2(S)$ times blown up Hopf surfaces.
\end{description}
{\bf GSS Conjecture}: {\it  All these properties are equivalent, and any class VII$_0^+$ surface  possesses a global spherical shell (GSS) i.e. an open submanifold biholomorphic to a standard neighborhood of $S^3$ in $\bb C^2$ which does not disconnect the surface. }
\vspace{2mm}\\
We have
$$(A)\Longleftrightarrow(B)\Longrightarrow(C)\Longrightarrow(D)$$
In fact $(A)\Longrightarrow(B)$ by the construction of GSS surfaces and $(B)\Longrightarrow(A)$ by \cite{DOT3},\\
$(A)\Longrightarrow(C)$ also by construction (see \cite{D1})
The implication $(C)\Longrightarrow(D)$ has been obtained by I.\,Nakamura \cite{N1,N2}.\\
The strategy developped  in \cite{Te2,Te3} is aimed to show that any surface in VII$_0^+$ satisfies condition $(C)$, therefore the solution to the following problem would be a step toward the conjecture:\\
{\bf Problem}: Let $\cal S\to \D$ be a family of compact surfaces over the disc such that for every $u\in\D^\star$, $S_u$ contains a GSS. Does $S_0$ contain a GSS ? In other words, are the surface with GSS  closed in families ?\\
To solve this problem we have to study families of surfaces in which curves do not belong to flat families, the volume of some curves in these families may be not uniformly bounded (see \cite{DT1}) and configurations of curves change. Favre normal forms of polynomial germs associated to surfaces with GSS, cannot be used because the discriminant of the intersection form is fixed. Moreover, if using the algorithm of \cite{OT} we put $F$ under the form $\Pi\s$, $\s$ is not fixed in the logarithmic family, depends on the blown up points and degenerates when a generic blown up point approaches the intersection of two curves.\\
Therefore this article focuses on the problem of finding new normal forms of contracting germs in intermediate cases of surfaces with fixed simple birational $\s$, such that surfaces are minimal or not and intersection matrices are not fixed. Since usual holomorphic objects, curves or foliations, do not fit in global family, it turns out that these birational structures depending on a finite number of parameters (in fact the number of parameters is exactly the dimension the local moduli spaces) could be the adapted notion. Moreover, our construction gives a contracting map $G=\Pi\s$ unique up to conjugation by elements of a group $L$ of diagonal linear mappings with coefficients equal to roots of unity. The existence of our special birational structure gives rise to  developing mappings $\widetilde{Dev}_j:\tilde S\to \bb P^2(\bb C)$ which contract an infinite number of rational curves onto a point $P$. The inverse image  of a small sphere by $\widetilde{Dev}_j$ gives a spherical shell in $\tilde S$, hence a GSS in $S$. This observation will be useful to prove the GSS conjecture.\\
\vspace{2mm}\\
This article is organized in the following way:\\
In {\bf section 2}, we introduce general notions on $G$-structures when $G=Bir(\bb P^2(\bb C))$, developing meromorphic mappings and recall known results on affine and projective structures, these being particular cases of birational structures.\\
In {\bf section 3}, we recall, in order to be self-contained, basic facts on surfaces with global spherical shells (GSS) which will be used all along this article, large families of marked surfaces with GSS which have been introduced in \cite{D6}. The proof of the main results hinge upon the fact that in all conjugation class of contracting germs of the form $\Pi\s$ there is a Favre germ. In order to be complete and to have clear notations we recall results on Oeljeklaus-Toma logarithmic moduli spaces of Kato surfaces with fixed intersection matrix of rational curves \cite{OT} with a slight modification (see Remark 3.\ref{erreur}). Favre germs $F$ which correspond to Kato surfaces which have a cycle with $\rho$ branches split into $\rho$ polynomial germs of simple type $F= F_1\circ\cdots\circ F_\rho$. Let $\cal F(\frak s,k,j)$ be the family of Favre germs of simple type (see Def. 3.\ref{Def 4.7}),
$$F(z_1,z_2)=(\l z_1z_2^\frak s  +\sum_{i=j}^{\frak s} b_iz_2^i+cz_2^{\frac{\frak s  k}{k-1}},\ z_2^k)$$
The associated surfaces have a cycle with exactly one branch.\\
In {\bf section 4} new germs are defined,  obtained by composition of $n$ blowing-ups ($2n$ parameters) and, if global twisted vector fields are expected to exist, of an extra invertible polynomial mapping tangent to the identity (the extra parameter $a_{l+K}$). If the surface contains a cycle of rational curves with only one branch this class of birational contracting germs is denoted by 
$\cal G=\cal G(p,q,r,s,l)$ and have the following form
\begin{itemize}
\item if the first blowing-up is not generic
$$G(z)=\left(z_1^{p+rl}z_2^{q+sl}+\sum_{i=0}^{l-1}a_i\bigl(z_1^rz_2^s\bigr)^{i+1}+a_{l+K}\bigl(z_1^rz_2^s\bigr)^{l+K+1},z_1^rz_2^s\right),$$
where  $K=\max\left\{0,\left[\frac{l-d}{r+s-1}\right]\right\}$, $a_0\in\bb C^\star$, $a_i\in\bb C$, $i=1,\ldots,l-1,l+K$,  and
\item if the first blowing-up is generic
$$G(z_1,z_2)=\left(\Bigl(z_1z_2^l+\sum_{i=0}^{l-1}a_iz_2^{i+1}+a_{l+K}z_2^{l+K+1}\Bigr)^pz_2^q, \ \Bigl(z_1z_2^l+\sum_{i=0}^{l-1}a_iz_2^{i+1}+a_{l+K}z_2^{l+K+1}\Bigr)^rz_2^s\right)$$
\end{itemize}
Among the blowing-ups there are $l$ generic blowing-ups, and $n-l$ non generic, determined by the matrix 
$$\left(\begin{array}{cc}p&q\\r&s\end{array}\right)\in Gl(2,\bb Z).$$
\\
First, we establish a correspondance between both families $\cal F(\frak s,k,j)$ and $\cal G=\cal G(p,q,r,s,l)$ when they correspond to the same sequence of blowing-ups giving the same intersection matrix of the rational curves. We provide precise relations between the integers involved in the construction and we explicit conditions which insure the existence of global vector fields.\\
We denote by $\Phi=\Phi(p,q,r,s,l)$ the group  of the  germs of biholomorphisms $\f:(\bb C^2,0)\to (\bb C^2,0)$ for which there exists $G,G'\in \cal G$ such that $G'=\f^{-1}G\f\in\cal G$.   Let
$L:=L(p,q,r,s,l)$ be the group of diagonal linear mappings $\f_{A,B}(z_1,z_2)=(Az_1,Bz_2)$ where $A,B$ satisfy the condition
$$B=A^rB^s,\quad A=A^{p+rl}B^{q+sl}$$
Then the following holds (see Prop. 4.\ref{presque-unicite} for a more detailed statement):
\begin{Prop}[unicity] There is an exact sequence 
$$0\to (\bb C,+)\to \Phi\to L\to {Id}$$
Moreover if $a_{l+K}=0$, then $\Phi=L$, i.e. the birational germ $G$ is unique up to a conjugation by a diagonal linear mapping whose coefficients are roots of unity.
\end{Prop}

Moving the parameters we have large families of surfaces with base $B_J$.  Is the canonical image of a stratum $B_{J,M}$ of surfaces with fixed intersection matrix $M$ in the Oeljeklaus-Toma coarse moduli space open ?  Do we obtain all possible surfaces ?\\
We know by \cite{D6} that outside the hypersurface $T_{J,\s}$ the family is versal. Here we show that $T_{J,\s}$ is a ramification hypersurface, in particular the canonical mapping from a stratum $B_{J,M}$ to the Oeljeklaus-Toma coarse moduli space is a ramified covering,  the mapping is surjective and vanishing of cohomology classes is due to ramification phenomena at $T_{J,\s}\cap B_{J,M}$. More precisely  (see section 3.4) we have the existence theorem
\begin{Th}[Main theorem]  Denote $\frak s:=p+q+l-1$ and $d:=(r+s)-(p+q)$. We choose $a_0\in\bb C^\star$ and $\e$ such that $\e^{r+s-1}=1$. Then\\
A) If $r+s-1$ does not divide $l-d$ or $\l\neq 1$ there is a bijective polynomial mapping 
$$\begin{array}{cccc}
f_{a_0,\e}:& \bb C^{l-1}&\longrightarrow&  \bb C^{l-1}\\
&a=(a_1,\ldots,a_{l-1})&\longmapsto&\Bigl(b_{p+q+1}(a),\ldots,b_{p+q+l-1}(a)\Bigr)
\end{array}$$
  such that 
$$G(z_1,z_2)=\left(\Bigl(z_1z_2^l+\sum_{i=0}^{l-1}a_iz_2^{i+1}\Bigr)^pz_2^q, \ \Bigl(z_1z_2^l+\sum_{i=0}^{l-1}a_iz_2^{i+1}\Bigr)^rz_2^s\right)$$
 is conjugated to the polynomial germ
$$F(z_1,z_2)=\Bigl(\l z_1z_2^\frak s+\sum_{i=p+q}^\frak s b_iz_2^i,\  z_2^{r+s}\Bigr),$$
where $\l$ depends only on $a_0$ by 4.\ref{dependancedekappa}.\\
B) If $l-d=K(r+s-1)$ and $\l=1$,  there is a bijective polynomial mapping 
$$\begin{array}{cccc}
f_{a_0,\e}:& \bb C^{l-1}\times \bb C&\longrightarrow&  \bb C^{l-1}\times \bb C\\
&a=(a_1,\ldots,a_{l-1},a_{l+K})&\longmapsto&\Bigl(b_{p+q+1}(a),\ldots,b_{p+q+l-1}(a),c(a)\Bigr)
\end{array}$$
such that
$$G(z_1,z_2)=\left(\Bigl(z_1z_2^l+\sum_{i=0}^{l-1}a_iz_2^{i+1}+a_{l+K}z_2^{2l+K}\Bigr)^pz_2^q, \ \Bigl(z_1z_2^l+\sum_{i=0}^{l-1}a_iz_2^{i+1}+a_{l+K}z_2^{2l+K}\Bigr)^rz_2^s\right)$$
 is conjugated to the polynomial germ
 $$F(z_1,z_2)=\Bigl(\l z_1z_2^\frak s+\sum_{k=p+q}^\frak s b_kz_2^k+ cz_2^{\frac{\frak s k(S)}{k(S)-1}}, z_2^{r+s}\Bigr).$$
\end{Th}
\begin{Cor} Let $S$ in class VII$_0^+$ containing a GSS. Suppose that the dual graph of the rational curves contains a cycle with only one branch, then $S$ admits a birational structure. In particular  Kato surfaces admit birational structures provided that $b_2(S)\le 3$.
\end{Cor}
In there are $\rho>1$ branches, there is for each intersection matrix $M$ an open set in the moduli space of Kato surfaces with intersection matrix $M$ obtained by birational germs $G_1\circ\cdots\circ G_\rho$ (see Cor. 1.3 in \cite{D6}).\\
In last section we show how to recover GSS from the existence of developing mappings $\widetilde{Dev}_j:\tilde S\to \bb P^2(\bb C)$.

\section{Birational structures on complex manifolds}
\subsection{preliminaries}
Here are classical definitions as in \cite{IKO,Kl}, in the context of complex manifolds. Notice that $\dim X=\dim Y$.
\begin{Def} Let $Y$ be a complex manifold, $G$ a Lie group acting holomorphically on $Y$ on the left and $X$ a complex manifold of dimension $n$. A $(G,Y)$-structure on $X$ is a maximal atlas of $X$, $\phi_i:U_i\to Y$ such that transition maps 
$$\phi_{ij}:=\phi_i\circ \phi_j^{-1}:\phi_j(U_i\cap U_j)\to \phi_i(U_i\cap U_j)$$
are locally elements of $G$.\\
Given two $(G,Y)$-manifolds, a $(G,Y)$-morphism  $f:X_1\to X_2$ is a holomorphic mapping such that for any charts $\phi_i:U_i\to Y$, $\psi_j:V_j\to Y$ of $X_1$ and $X_2$ respectively  and every connected component $C$ of $U_i\cap f^{-1}(V_j)$, there exists $g\in G$ such that 
$$f_{\mid C}=\psi_j^{-1}\circ g\circ \phi_i.$$
An affine structure (resp. a projective structure) on $X$ is a $(G,Y)$-structure where $Y=\bb C^n$ and $G$ is the affine group $A(n,\bb C)=Gl(n,\bb C)\rtimes \bb C^n$ (resp. $Y=\bb P^n(\bb C)$ and $G=\bb P Gl(n+1,\bb C)$).
\end{Def}
 If $f:X_1\to X_2$ is a local diffeomorphism and $X_2$ is a $(G,Y)$-manifold, there exists a unique $(G,Y)$-structure on $X_1$ such that $f$ is a morphism of $(G,Y)$-manifolds. In particular if $f$ is a non ramified covering, $X_1$ has a canonical $(G,Y)$-structure.\\
Taking now $Y=\bb P^n(\bb C)$ and $G=Bir(\bb P^n(\bb C))$, $G$ is neither an algebraic group nor a finite dimensional Lie group \cite{BF}. Therefore we extend the previous definition:
\begin{Def} Let $X$ be a  complex manifold of dimension $n$. We say that $X$ admits a  birational structure if there is an atlas $(U_i,\f_i)_{i\in I}$ such that  holomorphic transition maps $b_{ij}:=\f_i\circ\f_j^{-1}:\f_j(U_i\cap U_j)\to \f_i(U_i\cap U_j)$ are the restriction of birational mappings of $\bb P^n(\bb C)$. 
\end{Def}
Affine or projective structures are birational structures. If $X$ admits a birational structure  and $\Pi:X'\to X$ is a blowup, then $X'$ admits a unique birational structure such that $\Pi$ is a $(Bir(\bb P^n(\bb C),\bb P^n(\bb C))$-morphism. 

\begin{Ex} Let $X$ be compact riemann surface, then $X$ admits a birational structure. In fact, if $g(X)\ge 2$, $X$ is the quotient of the upper half-plane $\bb H=\{z\in\bb C\mid \Im z>0\}$ by a Fuchsian group, hence a subgroup of $PSL(2,\bb R)$ which acts on $\bb H$ as $\left(\begin{matrix}a&b\\c&d\end{matrix}\right).\,z=\frac{az+b}{cz+d}$.
\end{Ex}

As when $G$ is a Lie group we have a developing mapping $Dev:\tilde X\to \bb P^n(\bb C)$, however $Dev$ is now meromorphic.
\begin{Lem} \label{Extbir} Let $X$ be a complex manifold endowed with a birational structure $(U_i,\f_i)_{i\in I}$ and $p:\tilde X\to X$ its universal covering space. Let $x\in X$, $\G=\pi_1(X,x)$ be the fundamental group with base point $x$. Then for each $x_0\in p^{-1}(x)$ there is a $\G$-equivariant meromorphic developing mapping $Dev_{x_0}:\tilde X\to \bb P^n(\bb C)$, in other words there is a group morphism $h:\G\to Bir(\bb P^n(\bb C))$ such that
$$\forall \g\in\G, \quad Dev_{x_0}\circ \g=h(\g)\circ Dev_{x_0}.$$
Moreover $Dev_{x_0}$ is holomorphic in a neighbourhood of $x_0$.
\end{Lem}
Proof: It is sufficient to prove the extension along any path with base point $x_0\in p^{-1}(x)$.  Let $x_1\in \tilde X$ and $\g:[0,1]\to \tilde X$ a path joining $x_0=\g(0)$ to $x_1=\g(1)$. We cover $\g([0,1])$ by open domains  of charts $(U_0,\f_0),\ldots,(U_p,\f_p)$, such that $U_i\cap U_j\neq\emptyset$ if and only if $0\le i\le p-1$ and $j=i+1$. We prove by induction on $1\le j\le p$ that $\f_0:U_0\to \bb P^n(\bb C)$ admits a meromorphic extension $Dev_{x_0}$ on $U_0\cup \cdots \cup U_j$. For  $j=1$, setting ${Dev_{x_0}}_{\mid U_0}=\f_0$ and ${Dev_{x_0}}_{\mid U_1}=b_{01}\circ \f_1$. Let $b_{i, i+1}:=\f_{i}\circ\f_{i+1}^{-1}:\f_{i+1}(U_{i , i+1})\to \f_{i}(U_{i, i+1})$. By assumption $b_{i, i+1}$ extends birationally to $\bb P^n(\bb C)$. We suppose that $Dev_{x_0}$ has been extended along $U_0\cup\cdots\cup U_{j-1}$ for $j\ge 1$ setting
 $${Dev_{x_0}}_{\mid U_{j-1}}:U_{j-1} \to \bb P^n(\bb C), \quad x\mapsto Dev_{x_0}(x)= b_{01}\circ\cdots\circ b_{j-2,j-1}\circ\f_{j-1}(x)$$
and we define
$${Dev_{x_0}}_{\mid U_j}:U_j \to \bb P^n(\bb C), \quad x\mapsto Dev_{x_0}(x)=b_{01}\circ\cdots\circ b_{j-1,j}\circ\f_j(x)$$
For $x\in U_{j-1}\cap U_j$,  we have ${Dev_{x_0}}_{\mid U_{j-1}}(x)={Dev_{x_0}}_{\mid U_j}(x)$. \hfill$\Box$\\

From the lemma we obtain immediately
\begin{Prop} If $X$ is compact simply connected of dimension $n$ and admits a birational structure then the algebraic dimension of $X$ is equal to $n$.
\end{Prop}
\begin{Cor} A non projective K3 surface has no birational structure.
\end{Cor}

\subsection{Birational structures on non-Kähler complex surfaces}
\begin{Th}[\cite{BHPV84,KO60,KO64}] Any compact complex  non-Kählerian surface has a unique minimal model $X$ in the following classes
\begin{itemize}
\item Class VI$_0$, $b_1(X)$ is odd and geometric genus satisfies $p_g>0$. $X$ is an elliptic surface and admits, by \cite{IKO}, an holomorphic affine structure.
\item Class VII$_0$, $b_1(X)=1$, $p_g=0$, and Kodaira dimension $\kappa(X)\le 0$. 
\begin{description}
\item{(i)} If $\kappa(X)=-\infty$, $b_2(X)=0$ and $X$ contains a curve, then $X$ is a Hopf surface \cite{KO60} and admits a finite covering by a primary Hopf surface. Any primary Hopf surface is isomorphic to $\bb C^2\setminus\{0\}/H$ where $H$ is an infinite cyclic group generated by a contraction
$$g:(z_1,z_2)\mapsto (\a z_1+\l z_2^m,\b z_2), \quad 0<|\a|\le |\b|<1,\quad  (\b^m-\a)\l=0,\quad  m\ge 1.$$
If $(m-1)\l=0$, $X$ admits a holomorphic affine structure \cite{KO}, p93. In all cases the contraction is an invertible polynomial mapping hence birational, therefore $X$ admits a birational structure.
\item{(ii)} If $\kappa(X)=-\infty$, $b_2(X)=0$ and $X$ contains no curve, then $X$ is a Inoue surface by \cite{Te94} and admits a holomorphic affine structure \cite{IKO}.
\item{(iii)} If $\kappa(X)=0$, then $X$ is a secondary Kodaira surface \cite{BHPV84}, $b_2(X)=0$. By \cite{IKO}, $X$ admits an affine holomorphic structure.
\item{(iv)} Surfaces with $b_2(X)>0$. The only known surfaces are Kato surfaces.
\end{description}
\end{itemize}
\end{Th}
All compact non-Kähler surfaces admit affine structures but some Hopf surfaces and surfaces in class VII$_0^+$. \\
{\bf Conjecture}: Any compact non-Kähler complex surface admits a birational structure.

\section{Surfaces with Global Spherical Shells}

\subsection{Basic constructions}\label{basicconstructions}
\begin{Def} Let $S$ be a compact complex surface. We say that $S$ contains a global spherical shell, if there is a biholomorphic map $\f:U\to S$ from a neighbourhood $U\subset \bb C^2\setminus\{0\}$ of the sphere $S^3$ into $S$ such that $S\setminus \f(S^3)$ is connected.
\end{Def}
Hopf surfaces are the simplest examples of surfaces with GSS.\\

Let $S$ be a surface containing a GSS with $n=b_2(S)$. It is known that $S$ contains $n$ rational curves and to each curve it is possible to associate a contracting germ of mapping $F=\Pi\s=\Pi_0\cdots\Pi_{n-1}\s:(\bb C^2,0) \to (\bb C^2,0)$ where $\Pi=\Pi_0\cdots\Pi_{n-1}:B^\Pi\to B$ is a sequence of $n$ blowing-ups  and $\s$ is a germ of isomorphism (see \cite{D1}). \begin{Def}\label{Un-arbre-RecEnoki} Let $S$ be a  surface containing a GSS, with $n=b_2(S)$. A {\bf Enoki covering} of $S$ is an open covering $\cal U=(U_i)_{0\le i\le n-1}$ obtained in the following way:
\begin{itemize}
\item $W_0$ is the  ball of radius $1+\e$  blown up at the origin, $C_0=\Pi_0^{-1}(0)$, $B'_0\subset \subset B_0$ are small balls centered at $O_0=(a_0,0)\in W_0$, $U_0=W_0\setminus B'_0$,
\item For $1\le i\le n-1$, $W_i$ is the  ball $B_{i-1}$  blown up at $O_{i-1}$, $C_i=\Pi_i^{-1}(O_{i-1})$, $B'_{i}\subset \subset B_{i}$  are small balls centered at $O_{i}\in W_i$, $U_i=W_i\setminus B'_i$.
\end{itemize}
The  pseudoconcave boundary of $U_i$ is patched with the pseudoconvex boundary of $U_{i+1}$ by $\Pi_i$, for $i=0,\ldots,n-2$ and the pseudoconcave boundary of $U_{n-1}$ is patched with the pseudoconvex boundary of $U_0$ by $\s\Pi_0$, where
$$\begin{array}{cccc}
\s:&B(1+\e)&\to&W_{n-1}\\
&z=(z_1,z_2)&\mapsto&\s(z)
\end{array}$$
is biholomorphic on its image, satisfying $\s(0)=O_{n-1}$. 
\end{Def}

If we want to obtain a minimal surface, the sequence of blowing-ups has to be made in the following way:
\begin{itemize}
\item $\Pi_0 $ blows up the origin of the two dimensional unit ball $B$, 
\item $\Pi_1$ blows up a point $O_0\in C_0=\Pi_0^{-1}(0)$,\ldots 
\item $\Pi_{i+1}$ blows up a point $O_{i}\in C_{i}=\Pi_{i}^{-1}(O_{i-1})$, for $i=0,\ldots,n-2$, and 
\item $\s:\bar B\to B^\Pi$ sends isomorphically a neighbourhood of $\bar B$ onto a small ball in $B^\Pi$ in such a way that $\s(0)\in C_{n-1}$. 
\end{itemize}
Each $W_i$ is covered by two charts with coordinates $(u_i,v_i)$ and $(u'_i,v'_i)$ in which $\Pi_i$ writes $\Pi_i(u_i,v_i)=(u_iv_i+a_{i-1},v_i)$ and $\Pi_i(u'_i,v'_i)=(v'_i+a_{i-1},u'_iv'_i)$. In these charts the exceptional curves has always the equations $v_i=0$ and $v'_i=0$.\\
A  blown up point $O_i\in C_i$ will be called {\bf generic} if it is not at the intersection of two curves.
The data $(S,C)$ of a surface $S$ and of a rational curve in $S$ will be called a {\bf marked surface}.

We assume that $S$ is minimal and that we are in the intermediate case, therefore there is at least one blowing-up at a generic point, and one at the intersection of two curves (hence $n\ge 2$). If there is only one branch i.e. one regular sequence and if we choose $C_0$ as being the curve which induces the root of the branch, we suppose that
\begin{itemize}
\item
 $\Pi_1$ is a generic blowing-up, 
 \item $\Pi_{n-1}$ blows-up the intersection of $C_{n-2}$ with another rational curve and
 \item $\s(0)$ is one of the two intersection points of $C_{n-1}$ with the previous curves.
 \end{itemize}
The  Enoki covering is obtained as in the following picture:
\begin{center}\label{Enoki}
\includegraphics[width=12cm]{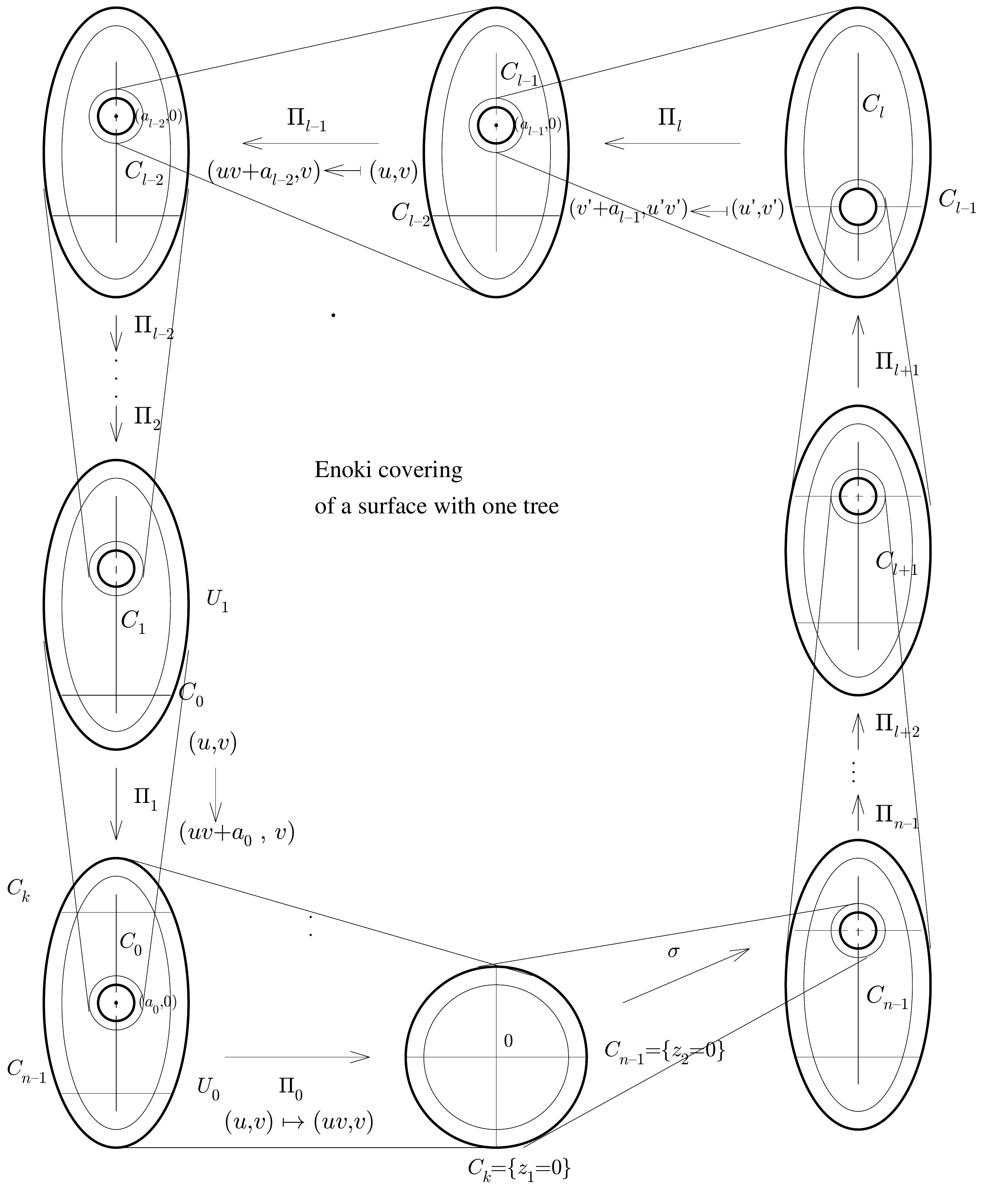}
\end{center}
where 
\begin{itemize}
\item $1\le l\le n-1$ and $n\ge 2$. If all, but one, blowing-ups are generic, then $l=n-1$ 
\item For $i=1,\ldots,l-1$, $\Pi_i(u_i,v_i)=(u_iv_i+a_{i-1},v_i)$ are generic blowing-ups,
\item $\Pi_{l}(u'_{l},v'_{l})=(v'_{l}+a_{l-1},u'_{l}v'_{l})$ is also generic, but $O_l$ is the origin of the chart $(u'_l,v'_l)$,
\item For $i=l+1,\ldots,n-1$, $\Pi_i(u_i,v_i)=(u_iv_i,v_i)$ or $\Pi_i(u'_i,v'_i)=(v'_i,u'_iv'_i)$ are blowing-ups at the intersection of two curves.
\end{itemize}

\subsection{Large families  of marked surfaces\label{logfamilies}}
With the previous notations, we consider global families of minimal compact surfaces with the same charts, parameterized by the coordinates of the blown up points on the successive exceptional curves obtained in the construction of the surfaces and such that any marked surface with GSS $(S,C_0)$ belongs to at least one of these families.  More precisely,  let $F(z)=\Pi_0\cdots\Pi_{n-1}\s(z)$ be a germ associated to any marked surface $(S,C_0)$ with $tr(S)=0$. In order to fix the notations we suppose that $C_0=\Pi_0^{-1}(0)$ meets two other curves (see the picture after definition \thesection.\ref{Un-arbre-RecEnoki}), hence $\s(0)$ is the intersection of $C_{n-1}$ with another curve. We suppose that
 $$\part_1\s_2(0)=0.$$
  We denote by $I_\infty(C_0)\subset \{0,\ldots,n-1\}$ the subset of indices which correspond to blown up points at infinity, that is to say,
$$I_\infty(C_0):=\bigl\{i\mid O_i\ {\rm is\ the\ origin\ of\ the\ chart}\ (u'_i,v'_i)\bigr\}.$$
Each generic blow-up
$$\Pi_i(u_i,v_i)=(u_iv_i+a_{i-1},v_i)\quad {\rm or}\quad \Pi_i(u'_i,v'_i)=(v'_i+a_{i-1},u'_iv'_i)$$
 may be deformed moving the blown up point $(a_{i-1},0)$. If we do not want to change the configuration  we take
 $$\begin{array}{l}
 {\rm for\  all}\  \k=0,\ldots,\rho-1\quad {\rm (with }\ n_0=0{\rm)},\\
 \\
 \hspace{20mm}\left\{
 \begin{array}{l}
 a_{n_1+\cdots+n_\kappa}\in \bb C^\star,\\
 \\
 \forall i,\  1\le i\le l_{\k} -1,\qquad  a_{n_1+\cdots+n_\kappa +i}\in \bb C, \\
 \\
  \forall j,\  0\le j\le n_{\k+1}-l_{\k}-1,\qquad a_{n_1+\cdots+n_\kappa+l_{\k}+j}=0.
  \end{array}\right.
  \end{array}$$
  The mapping $\s$ is supposed to be fixed.  We obtain a large family of compact surfaces which contains $S$ such that all the surfaces $S_a$ have the same intersection matrix $$M=M(S_a)=M(S),$$
   therefore are logarithmic deformations. For $J=I_\infty(C_0)$ we denote 
  this family
$$\Phi_{J,M,\s}:\cal S_{J,M,\s}\to B_{J,M}$$
where
$$\begin{array}{l}
B_{J,M}\\
\\
:=\bb C^\star\times \bb C^{l_0-1}\times\{0\}^{n_1-l_0}\times\cdots\times \bb C^\star\times \bb C^{l_\kappa -1}\times\{0\}^{n_{\k+1}-l_\k}\times\cdots\times \bb C^\star\times\bb C^{l_{\rho-1}-1}\times\{0\}^{n_\rho-l_{\rho-1}}\\
\\
\simeq\bb C^\star\times \bb C^{l_0-1}\times\cdots\times \bb C^\star\times \bb C^{l_\kappa -1}\times\cdots\times \bb C^\star\times\bb C^{l_{\rho-1}-1}\end{array}$$
and $n_1+\cdots+n_\rho=n$.\\

In $\cal S_{J,M,\s}$ there is a flat family of divisors $\cal D\subset \cal S$ with irreducible components 
$$\cal D_i, \quad i=0,\ldots,n-1,$$
such that for every $a\in B_{J,M}$, $M=(D_{i,a}.D_{j,a})_{0\le i,j\le n-1}$.
  We may extend this family towards smaller or larger strata which produce minimal surfaces:
  \begin{itemize}
  \item On one hand, {\bf towards a unique Inoue-Hirzebruch surface}:
  Over 
  $$\bb C^{l_0}\times\{0\}^{n_1-l_0}\times\cdots\times \bb C^{l_\kappa}\times\{0\}^{n_{\k+1}-l_\k}\times\cdots\times \bb C^{l_{\rho-1}}\times\{0\}^{n_\rho-l_{\rho-1}}\simeq \bb C^{l_0}\times\cdots\times \bb C^{l_\kappa}\times\cdots\times \bb C^{l_{\rho-1}},$$ 
$$\Phi_{J,\s}:\cal S_{J,\s}\to \bb C^{l_0}\times\cdots\times \bb C^{l_\kappa}\times\cdots\times \bb C^{l_{\rho-1}}.$$
If for an index  $\kappa$, $a_{n_1+\cdots+n_\kappa}=0$, there is a jump in the configuration of the curves. For instance, if for all $\kappa$, $\kappa=0,\ldots,\rho-1$
 $$a_{n_1+\cdots+n_\kappa}=\cdots=a_{n_1+\cdots+n_\kappa+l_\kappa-1}=0$$
 we obtain a Inoue-Hirzebruch surface.  To be more precise the base 
 $$\bb C^{l_0}\times\cdots\times \bb C^{l_\kappa}\times\cdots\times \bb C^{l_{\rho-1}}$$
 splits into locally closed submanifolds called {\bf strata}
\begin{itemize}
\item the Zariski open set $\bb C^\star\times\bb C^{l_0-1}\times\cdots\times \bb C^\star\times\bb C^{l_\k-1}\times\cdots\bb C^\star\times\bb C^{l_{\rho-1}-1}$,
\item $\rho=C^1_\rho$ codimension one strata 
$$\bb C^\star\times\bb C^{l_0-1}\times\cdots\times \{0\}\times \bb C^\star\times\bb C^{l_\k-2}\times\cdots\times \bb C^\star\times\bb C^{l_{\rho-1}-1}, \quad 0\le \k\le \rho-1,$$
\item $C^p_{\rho+p-1}$ codimension $p$ strata, $1\le p:=p_0+\cdots +p_{\rho-1}\le l_0+\cdots+l_{\rho-1}$,
$$\{0\}^{p_0}\times \bb  C^\star\times \bb C^{l_0-p_0-1}\times\cdots\times\{0\}^{p_\k}\times\bb C^\star \times \bb C^{l_\k-p_\k-1}\times\cdots\times\{0\}^{p_{\rho-1}}\times \bb C^\star\times \bb C^{l_{\rho-1}-p_{\rho-1}-1}$$
\end{itemize}
 \item On second hand, {\bf towards Enoki surfaces}. If for all indices such that $O_i$ is at the intersection of two rational curves, in particular for $i\in J$, the blown up point $O_i$  is moved to $O_i=(a_i,0)$ with $a_i\neq 0$, all the blown up points become generic, the trace of the contracting germ is different from $0$. We obtain also all the intermediate configurations.\\
\begin{Prop}[\cite{D6} Prop.2.6] \label{towardsEnoki} There is a monomial holomorphic function $t:\bb C^{\Card J}\to \bb C$ depending on the variables $a_j$, $j\in J$ such that over $B_J:=\{|t(a)|<1\}\subset \bb C^n$, the family $\Phi_{J,\s}:\cal S_{J,\s}\to B_{J}$ may be extended and for every $a\in B_J$, $t(a)=\tr(S_a)$. 
\end{Prop}
\end{itemize}
 Remain non minimal surfaces: we still extend the previous family on a small neighbourhood $\wh{B_J}$ of $B_J$, moving the blown up point transversally to the exceptional curves $C_i=\{v_i=0\}\cup\{v'_i=0\}$, introducing $n$ new parameters
 $$\Pi_i(u_i,v_i)=(u_iv_i+a_{i-1},v_i+b_{i-1}),  \quad {\rm or}\quad \Pi_i(u'_i,v'_i)=(v'_i+a_{i-1},u'_iv'_i+b_{i-1}), \quad |b_{i-1}|<<1,$$
 we obtain 
 $$\wh{\Phi}_{J,\s}:\wh{\cal S}_{J,\s}\to \wh{B}_J,$$
 with $\dim \wh{B_J}=2n=2b_2$.
 Since for any $(a,b)\in \wh{B_J}$, $h^1(S_{a,b},\T_{a,b})=2b_2(S_{a,b})+h^0(S_{a,b},\T_{a,b})$, there are some questions:
 \begin{itemize}
 \item Are the parameters $a_i,b_i$, $i=0,\ldots,n-1$, effective ? By \cite{D6}, they are generically effective.
 \item Which parameter to add when $h^1(S_{a,b},\T_{a,b})=2b_2(S_{a,b})+1$ in order to obtain a complete family ?
\item If we choose $\s=Id$ or more generally an invertible polynomial mapping, we obtain a birational polynomial germs. Does this families contain all the isomorphy classes of surfaces with fixed intersection matrix $M$ ?  
\end{itemize}

\begin{Rem} It is difficult to determine the maximal domain $\wh{B}_J$ over which $\wh{\Phi}_{J,\s}$ may be defined. When the surface is minimal, i.e. when $b=(b_0,\ldots,b_{n-1})=0$, $F_{a,b}(0)=0$. However, when $b\neq 0$, the fixed point $\z=(\zeta_1,\zeta_2)$  moves and the existence condition for the corresponding surface is that the eigenvalues $\l_1$ and $\l_2$ of $DF_{a,b}(\z)$ satisfy $|\l_i|<1$, $i=1,2$.
\end{Rem}

\subsection{Oeljeklaus-Toma logarithmically versal family}\label{sectionOT}
The goal is to compare the Oeljeklaus-Toma logarithmic families of surfaces with the strata in large families of surfaces $\Phi_{J,M\s}:\cal S_{J,M,\s}\to B_{J,M}$ which have the same intersection matrix $M$. In the case of surfaces with only one branch it turns out that we obtain all the surfaces.\\
We recall the results of \cite{OT} used in the sequel with a small correction described in the remark 3.\ref{erreur}. \\
All surfaces of intermediate type may be obtained from a polynomial germ in the following normal form obtained by \cite{Fav} and improved by \cite{OT}. 
$$F(z_1,z_2)=(\l z_1z_2^\frak s +P(z_2)+cz_2^{\frac{\frak s  k}{k-1}},z_2^k)\leqno{(CG)}$$
where $k,\frak s  \in\bb Z$, $k>1$, $\frak s  >0$, $\l\in\bb C^\star$, 
$$P(z_2)=c_jz_2^j+c_{j+1}z_2^{j+1}+\cdots+c_\frak s  z_2^\frak s  $$
is a complex polynomial satisfying the conditions
$$0<j<k, \quad j\le \frak s  ,\quad c_j=1,\quad c\in \bb C, \quad \gcd\{k,m\mid c_m\neq 0\}=1$$
with $c=0$ whenever $\frac{\frak s  k}{k-1}\not\in\bb Z$ or $\l\neq 1$.
\begin{Lem}[\cite{OT},\S 4] \label{condconj} Two polynomial germs $F$ and 
$$\tilde F(z_1,z_2)=\Bigl(\tilde\l z_1z_2^{\widetilde{ \frak s}  }+\tilde P(z_2)+\tilde cz_2^{\frac{\widetilde{ \frak s}  \tilde k}{\tilde k-1}},\ z_2^{\tilde k}\Bigr),$$
in normal form $(CG)$ are conjugated if and only if there exists $\e\in\bb C$, $\e^{k-1}=1$ such that
$$\tilde k=k,\quad \widetilde {\frak s}  =\frak s  ,\quad \tilde\l=\e^\frak s  \l,\quad \tilde P(z_2)=\e^{-j}P(\e z_2), \quad \tilde c=\e^{\frac{\frak s  k}{k-1}}c.$$
\end{Lem}

Intermediate surfaces admitting a global non-trivial twisted vector field or a non-trivial section of the anticanonical line bundle are exactly those for which  $(k-1) \mid \frak s  $. When moreover $\l=1$ we have a non-trivial global vector field.
\begin{Def} Let $S$ be a surface containing a GSS. The least integer $\mu\ge 1$ such that there exists $\kappa\in\bb C^\star$ for which
$$H^0(S,K_S^{-\mu}\ot L^\kappa)\neq 0$$
is called the {\bf index} of $S$. 
\end{Def}
If $S$ is defined by the polynomial germ $$F(z_1,z_2)=(\l z_1z_2^\frak s  +P(z_2)+cz_2^{\frac{\frak s  k}{k-1}},z_2^k)\leqno{(CG)}$$
then by \cite{OT} Remark 4.5,
$$index(S):=\mu=\frac{k-1}{gcd(k-1,\frak s  )}.$$
Notice that these germs show the existence of a foliation whose leaves are defined by $\{z_2=constant\}$, however {\it they are not birational}.\\
The set of polynomial germs
$$F(z_1,z_2)=(\l z_1z_2^\frak s  +P(z_2),\ z_2^k)$$
with $c=0$ are called in pure normal form.
\begin{Def}[\cite{OT} Def 4.7] \label{Def 4.7} For fixed $k$ and $\frak s  $ and for a polynomial germ
$$F(z_1,z_2)=(\l z_1z_2^\frak s  +P(z_2)+cz_2^{\frac{\frak s  k}{k-1}},z_2^k)\leqno{(CG)}$$
we define inductively the following finite sequences of integers
$$j=:m_1<\cdots <m_\rho\le \frak s  ,\quad {\rm and} \quad k>i_1>i_2>\cdots >i_\rho=1,$$
by:
\begin{description}
\item{(i)} $m_1:=j$, $i_1:=\gcd(k,m_1)$,
\item{(ii)} $m_\a:=\min\bigl\{m>m_{\a-1} \mid c_m\neq 0, \gcd(i_{\a-1},m)<i_{\a-1}\bigr\}$, $i_\a=\gcd(k,m_1,\ldots,m_\a)=\gcd(i_{\a-1},m_\a)$,
\item{(iii)} $1=i_\rho:=\gcd(k,m_1,\ldots,m_{\rho-1},m_\rho)<\gcd(k,m_1,\ldots,m_{\rho-1})$.
\end{description}
We call $(m_1,\ldots,m_\rho)$ the {\bf type} of $F$ and $\rho$ the {\bf length of the type}. If $\rho=1$, we say that $F$ is of {\bf simple type}.
\end{Def}
By \cite{OT}, \S 6, the length of the  type is exactly the number $\rho$ of branches  previously introduced.
\begin{Rem}\label{erreur} 1) If the length is $\rho=1$, then $\gcd(k,j)=1$ and there is no extra condition on the coefficients $c_{j+1},\ldots,c_\frak s  $, therefore the parameter space of polynomial germs in pure form with integers $k,\frak s  $ and type $j$ is 
$$U_{k,\frak s  ,j}=\bb C^\star\times \bb C^{\frak s  -j}.$$
If the length of the type is $\rho\ge 2$, notice that by definition, we have $c_{m_\a}\in \bb C^\star$, $\a=1,\ldots,\rho$, and $c_{m_1}=c_j=1$, however  between $c_{m_\a}$ and $c_{m_{\a+1}}$, the coefficients
$$c_{m_\a+i_\a},c_{m_\a+2i_\a},\ldots,c_{m_\a+\left[\frac{m_{\a+1}-m_\a}{i_\a}\right]i_\a}\in \bb C$$
may take any value, 
but all the other coefficients from $c_{m_\a+1}$ to $c_{m_{\a+1}-1}$ should vanish.
Let
$$\e(k,\frak s  ,m_1,\ldots,m_\rho):=\sum_{\a=1}^{\rho-1}\left[\frac{m_{\a+1}-m_\a}{i_\a}\right] + t-m_\rho$$
then the parameter space of all the germs $F$ with the same integers $\frak s  ,k$ and of the same type $(m_1,\ldots,m_\rho)$ in pure form are parameterized by 
$$(\bb C^{\star})^\rho\times \bb C^{\e(k,\frak s  ,m_1,\ldots,m_\rho)}.$$
\end{Rem}
There exists a family of surfaces
 $$\cal S_{k,\frak s  ,m_1,\ldots,m_\rho}\to (\bb C^{\star})^\rho\times \bb C^{\e(k,\frak s  ,m_1,\ldots,m_\rho)}$$
 such that for every $u\in (\bb C^{\star})^\rho\times \bb C^{\e(k,\frak s  ,m_1,\ldots,m_\rho)}$, $S_u$ is associated to the germ $F_u$. We have

\begin{Th and Def}[\cite{OT}, thm 7.13] With the above notations we have:
\begin{description}
\item{$\bullet$} If $k-1$ does not divide $\frak s  $, the family
$$\cal S_{k,\frak s  ,m_1,\ldots,m_\frak s  }\to (\bb C^{\star})^\rho\times \bb C^{\e(k,\frak s  ,m_1,\ldots,m_\rho)}=:U_{k,\frak s  ,m_1,\ldots,m_\rho}$$
is logarithmically versal at every point and  contains all surfaces with parameters $\frak s  ,k$ and type $(m_1,\ldots,m_\rho)$.
\item{$\bullet$} If $k-1$ divides $\frak s  $, the family 
$$\cal S_{k,\frak s  ,m_1,\ldots,m_\rho}\to (\bb C^{\star})^\rho\times \bb C^{\e(k,\frak s  ,m_1,\ldots,m_\rho)}\times \bb C=:U_{k,\frak s  ,m_1,\ldots,m_\rho}$$
\begin{itemize}
\item is  logarithmically complete at every point,
\item is logarithmically versal at every point of 
$$U^{\l=1}_{k,\frak s  ,m_1,\ldots,m_\rho}:=(\bb C^\star)^{\rho-1}\times \bb C^{\e(k,\frak s  ,m_1,\ldots,m_\rho)}\times \bb C$$
 and its restriction 
 $$\cal S_{k,\frak s  ,m_1,\ldots,m_\rho}\to U^{\l=1}_{k,\frak s  ,m_1,\ldots,m_\rho}$$
  contains all surfaces with parameters $\frak s  ,k$ and type $(m_1,\ldots,m_\rho)$ admitting a non-trivial global vector field,
\end{itemize}
Moreover 
 its restrition 
$$\cal S_{k,\frak s  ,m_1,\ldots,m_\rho}\to \bb C\setminus\{0,1\}\times(\bb C^{\star})^{\rho-1}\times \bb C^{\e(k,\frak s  ,m_1,\ldots,m_\rho)}:=U^{\l\neq 1,c=0}_{k,\frak s  ,m_1,\ldots,m_\rho}$$
is logarithmically versal at every point and contains all surfaces with parameters $k$, $\frak s  $ and type $(m_1,\ldots,m_\rho)$ without non-trivial global vector fields.
\end{description}
We shall call this family the {\bf Oeljeklaus-Toma logarithmic family of parameters $k$, $\frak s  $ and type $(m_1,\ldots,m_\rho)$}.\\
\end{Th and Def}
By lemma 3\ref{condconj}, for fixed $k,\frak s  $ and type $(m_1,\ldots,m_\rho)$, $\bb Z/(k-1)$ acts on the germs in pure normal form.  By \cite{OT} (7.14),\\
\begin{description}
\item{$\bullet$} $\cal M_{k,\frak s  ,m_1,\ldots,m_\rho}:=U_{k,\frak s  ,m_1,\ldots,m_\rho}/ \bigl(\bb Z/(k-1)\bigr)$\qquad  if $k-1$ does not divide $\frak s  $,
\item
\item{$\bullet$} $\left\{\begin{array}{l}
\cal M^{\l\neq 1,c=0}_{k,\frak s  ,m_1,\ldots,m_\rho}:=U^{\l\neq 1,c=0}_{k,\frak s  ,m_1,\dots,m_\rho}/ \bigl(\bb Z/(k-1)\bigr)\\
\\
\\
 \cal M^{\l=1}_{k,\frak s  ,m_1,\ldots,m_\rho}:=U^{\l=1}_{k,\frak s  ,m_1,\ldots,m_\rho}/ \bigl(\bb Z/(k-1)\bigr),\end{array}\right.$ if $k-1$ divides $\frak s  $,
\end{description}
\vspace{2mm}
are coarse moduli spaces, moreover the canonical mappings are ramified covering spaces. By lemma 3.\ref{condconj}, the ramification set is the union $T_{k,\frak s,m_1,\ldots,m_\rho}$ (resp. $T_{k,\frak s,m_1,\ldots,m_\rho}^{\l\neq 1,c=0}$, $T_{k,\frak s,m_1,\ldots,m_\rho}^{\l=1}$) of  hypersurfaces $\{c_i=0\}$, with $j+1\le i\le \frak s$ such that $c_i\in\bb C$, in particular
$$U_{k,\frak s,m_1,\ldots,m_\rho}\setminus T_{k,\frak s,m_1,\ldots,m_\rho}\to \cal M_{k,\frak s,m_1,\ldots,m_\rho}$$
$$U_{k,\frak s,m_1,\ldots,m_\rho}^{\l\neq 1,c=0}\setminus T_{k,\frak s,m_1,\ldots,m_\rho}^{\l\neq 1,c=0}\to \cal M_{k,\frak s,m_1,\ldots,m_\rho}^{\l\neq 1,c=0}$$
$$U_{k,\frak s,m_1,\ldots,m_\rho}^{\l=1}\setminus T_{k,\frak s,m_1,\ldots,m_\rho}^{\l=1}\to \cal M_{k,\frak s,m_1,\ldots,m_\rho}^{\l=1}$$
are non ramified covering spaces having $k-1$ sheets.

\begin{Rem} When $k-1$ divides $\frak s  $,  all the surfaces over the fiber $(\l,a,b)\times \bb C$ with $(\l,a,b)\in \bb C\setminus\{0,1\}\times (\bb C^\star)^{\rho-1}\times \bb C^{\e(k,\frak s  ,m_1,\ldots,m_\rho)}$ are isomorphic. Moreover
$$U_{k,\frak s  ,m_1,\ldots,m_\rho}^{\l=1}/\bigl(\bb Z/(k-1)\bigr)\cup U^{\l\neq 1,c=0}_{k,\frak s  ,m_1,\ldots,m_\rho}/\bigl(\bb Z/(k-1)\bigr)$$
is not separated. In fact, denote by
$$F_{\l,c}(z_1,z_2)=(\l z_1z_2^\frak s  +P(z_2)+cz_2^{\frac{\frak s  k}{k-1}},z_2^k).$$
Then any neighbourhood of $F_{1,c}$ with $c\neq 0$ meets any neighbourhood of $F_{1,0}$ because if $\l\neq 1$,
$$F_{\l,c}\sim F_{\l,0}.$$
\end{Rem}

\begin{Prop and Def} If $k-1$ divides $\frak s  $, the restriction 
$$\cal S^0_{k,\frak s ,m_1,\ldots,m_\rho}\to (\bb C^{\star})^{\rho}\times \bb C^{\e(k,\frak s  ,m_1,\ldots,m_\rho)}:=U^{c=0}_{k,\frak s  ,m_1,\ldots,m_\rho}$$
of the family
$$\cal S_{k,\frak s  ,m_1,\ldots,m_\rho}\to (\bb C^{\star})^{\rho}\times \bb C^{\e(k,\frak s  ,m_1,\ldots,m_\rho)}\times \bb C:=U_{k,\frak s  ,m_1,\ldots,m_\rho}$$
will be called the {\bf Oeljeklaus-Toma family of pure surfaces}.
 It
is versal at every point of 
$$\bb C\setminus\{0,1\}\times(\bb C^{\star})^{\rho-1}\times \bb C^{\e(k,\frak s  ,m_1,\ldots,m_\rho)}$$
 and effective at every point of 
 $$\{1\}\times(\bb C^{\star})^{\rho-1}\times \bb C^{\e(k,\frak s  ,m_1,\ldots,m_\rho)}.$$
\end{Prop and Def}
Since the hypersurface $(\bb C^{\star})^{\rho}\times \bb C^{\e(k,\frak s  ,m_1,\ldots,m_\rho)}\times\{0\}$ is invariant under the action of $\bb Z/(k-1)$ by (\ref{condconj}), the projection 
$$pr:(\bb C^{\star})^{\rho}\times \bb C^{\e(k,\frak s  ,m_1,\ldots,m_\rho)}\times\bb C\to (\bb C^{\star})^{\rho}\times \bb C^{\e(k,\frak s  ,m_1,\ldots,m_\rho)}\times\{0\}$$
 induces a holomorphic mapping
$$p:(\bb C^{\star})^{\rho}\times \bb C^{\e(k,\frak s  ,m_1,\ldots,m_\rho)}\times\bb C/\bigl(\bb Z/(k-1)\bigr)\to (\bb C^{\star})^{\rho}\times \bb C^{\e(k,\frak s  ,m_1,\ldots,m_\rho)}\times\{0\}/\bigl(\bb Z/(k-1)\bigr).$$

\section{Special birational structures on compact surfaces and birational germs}

\subsection{Birational germs associated to marked surfaces with one branch}
\subsubsection{Invariants and geometric properties}
In this section we define new normal forms of contracting germs, then we determine geometric properties and  conditions for the existence of global vector fields.\\

Let $(S,C_0)$ be a marked surface with GSS and let $M$ be the intersection matrix of the rational curves. We suppose that $C_0$ is the root of the unique branch (see picture in section 3.1). Then
we have
$$\Pi_l\cdots\Pi_{n-1}(u'',v'')=\left(u''^pv''^q+a_{l-1},u''^rv''^s\right)$$
where $(u'',v'')=(u,v)$ or $(u'',v'')=(u',v')$, $\left(\begin{array}{cc}p&q\\r&s\end{array}\right)$
is the composition of matrices  $A=\left(\begin{array}{cc}1&1\\0&1\end{array}\right)$ or $A'=\left(\begin{array}{cc}0&1\\1&1\end{array}\right)$, the last one being equal to $A'$. We set 
$$\d:=ps-qr=\pm 1,$$
$$1\le d:=(r+s)-(p+q)< r+s.$$
Moreover
$$\Pi_0\cdots\Pi_{l-1}(u,v)=\left(uv^l+\sum_{i=0}^{l-2}a_iv^{i+1},v\right)$$
Hence 
$$G(z)=\Pi\s(z)=\left(\s_1(z)^{p+rl}\s_2(z)^{q+sl}+\sum_{i=0}^{l-1}a_i\Bigl(\s_1(z)^r\s_2(z)^s\Bigr)^{i+1},\s_1(z)^r\s_2(z)^s\right),$$
where $\s$ is a germ of biholomorphism. \\

 If there is no global vector fields the number of parameters given by the blown up points is $2n$ as the expected number of parameters of the  versal deformation, therefore the question arises to know if with $\s=Id$ we obtain locally versal families. If there are non trivial global vector fields  we need (at least) an extra parameter.  We add this parameter  by the composition $\bar\s\Pi_0\cdots\Pi_{l-1}\Pi_l\cdots\Pi_{n-1} Id$ where
$$\bar\s(u,v)=(u+a_{l+K}v^{l+K+1},v),\quad K\ge 0,$$
where $K$ will be chosen in proposition \ref{CNSchamp}.
We obtain a new mapping (denoted in the same way)
$$\begin{array}{lcl}
G(z)&=&\bar\s\Pi_0\cdots\Pi_{l-1}\Pi_l\cdots\Pi_{n-1} Id(z)\\
&&\\
&=&\dps \left(z_1^{p+rl}z_2^{q+sl}+\sum_{i=0}^{l-1}a_i\bigl(z_1^rz_2^s\bigr)^{i+1} + a_{l+K}(z_1^rz_2^s)^{l+K+1},z_1^rz_2^s\right).
\end{array}$$
 We obtain large families $\cal S_{J,\s_{a_{l+K}}\to B_{J}}$ and we shall prove that the stratum $B_{J,M}$ is a ramified covering over the OT moduli space of  marked surfaces with GSS and intersection matrice $M$.

\begin{Lem}\label{existchamp} Let 
$$G(z)=\Pi\s(z)=\left(z_1^{p+rl}z_2^{q+sl}+\sum_{i=0}^{l-1}a_i\bigl(z_1^rz_2^s\bigr)^{i+1}+a_{l+K}\bigl(z_1^rz_2^s\bigr)^{l+K+1},z_1^rz_2^s\right),$$
 then the associated surface $S=S(G)$ admits a non trivial global twisted vector field if and only if
 $$u=\frac{p+s+rl-1-\d }{r+s-1},\quad v=\frac{r+q+sl-1+\d }{r+s-1},\quad{\rm where}\quad \d :=ps-qr$$
 are positive integers. Moreover this twisted vector field is a global vector field if and only if 
 $$\d a_0^uk(S)=1.$$
\end{Lem}
Proof: We have by a straightforwad computation
$$\det DG(z)=(ps-qr)z_1^{p+r(l+1)-1}z_2^{q+s(l+1)-1}.$$
By \cite{DO1}, there exists a non trivial global twisted vector field $\t\in H^0(S,\T\ot L^\l)$ on $S$ if and only if there is a global twisted  section of the anticanonical bundle $\o\in H^0(S,K^{-1}\ot L^\k)$. Moreover the twisting factors satisfy the relation  $\l=k(S)\k$. The section $\t$ is a global vector field if $\l=1$ i.e.
\begin{equation}\label{kS}
\k=\frac{1}{k(S)}
\end{equation} Such a section exists if and only if there is a germ of $2$-vector field
(denoted in the same way)
$$\o(z)=z_1^uz_2^vA(z)\frac{\part}{\part z_1}\wedge \frac{\part}{\part z_2}$$
where $A(0)\neq 0$ such that $\o(G(z))=\k \det DG(z)\o(z)$, or equivalently,
$$(a_0z_1^rz_2^s+\cdots)^u(z_1^rz_2^s)^vA(G(z))=\k(ps-qr)z_1^{p+r(l+1)-1+u}z_2^{q+s(l+1)-1+v}A(z).$$
Comparing terms of lower degree, we obtain the necessary condition
$$a_0^u(z_1^rz_2^s)^{u+v}=\k(ps-qr)z_1^{p+r(l+1)-1+u}z_2^{q+s(l+1)-1+v}$$
therefore $u$ and $v$ satisfy the linear system
$$\left\{\begin{array}{lcl} r(u+v)&=&p+r(l+1)-1+u\\&&\\s(u+v)&=&q+s(l+1)-1+v\end{array}\right.$$
The determinant of the system is $\D=-r-s+1<0$ and the solution is
$$u=\frac{p+s+rl-1-\d }{r+s-1},\quad v=\frac{r+q+sl-1+\d }{r+s-1},\quad{\rm where}\quad \d :=ps-qr=\pm1.$$
Since $u$ and $v$ are the vanishing orders  of $\o$ along the curves, a necessary condition for the existence of $\o$ is that  $u$ and $v$ are positive integers. Cancelling the common factors we obtain
$$a_0^u=\k\d $$
and with relation (\ref{kS})
$$\k=\d a_0^u=\frac{1}{k(S)}.$$ 

 If $u$ and $v$ are integers,
$$(a_0+\cdots)^uA(G(z))=\k\d A(z),$$
with $a_0\neq 0$. Setting 
$$1+f(z)=\frac{\k\d }{(a_0+\cdots)^u},$$
we have 
$$A(G(z))=(1+f(z))A(z)$$
Therefore
$$ A(z)=\frac{A(0)}{\dps\prod_{j=0}^\infty \bigl(1+f(G^j(z))\bigr)},$$
the infinite product converges because $G$ is contractant.  This proves the existence of $\o$.\hfill$\Box$

\begin{Prop}\label{CNSchamp} Let 
$$G(z)=\Pi\s(z)=\left(z_1^{p+rl}z_2^{q+sl}+\sum_{i=0}^{l-1}a_i\bigl(z_1^rz_2^s\bigr)^{i+1}+a_{l+K}\bigl(z_1^rz_2^s\bigr)^{l+K+1},z_1^rz_2^s\right),$$
  and let $S=S(G)$ be  the associated surface. Then
the surface $S(G)$  admits a non trivial global twisted vector field if and only if there exists an integer $k\ge 0$ such that
 $$l=d+k(r+s-1), $$
 If this condition is fulfilled, we choose 
 $$K=k$$
and $S(G)$ admits a non trivial vector field if and only if  for $u=\frac{p+s+rl-1-\d }{r+s-1}\in\bb N^\star$,
 $$\d a_0^uk(S)=1.$$
\end{Prop}
Proof: With notations of lemma \thesection.\ref{existchamp}, we have to show that $u$ and $v$ are integers if and only if $l=d+k(r+s-1)$.\\
If $u$ and $v$ are integers, 
$$u+v= l+1+\frac{p+q+l-1}{r+s-1}\in\bb N,$$
where $p+q<r+s$. 
Therefore, $l=d+k(r+s-1)$. Conversely, if $l=d+k(r+s-1)$, it is easy to check that $u$ and $v$ are positive integers and the proof is left to the reader.\hfill$\Box$

\begin{Prop} \label{kSrs} Let 
$$G(z)=\Pi\s(z)=\left(z_1^{p+rl}z_2^{q+sl}+\sum_{i=0}^{l-1}a_i\bigl(z_1^rz_2^s\bigr)^{i+1}+a_{l+K}\bigl(z_1^rz_2^s\bigr)^{l+K+1},z_1^rz_2^s\right),$$
and $S=S(G)$ the associated surface. Then
$$k(S)=r+s.$$
\end{Prop}
Proof: The dual graph of the curves is composed of a cycle with (here) only one chain of rational curves called the tree or the branch. The proof is achieved by induction on the number $N\ge 1$ of singular sequences. We denote as in \cite{D1}  $$a(S)=(s_{k_1}\cdots s_{k_N}r_l),$$
 where for any $k\ge 1$, $s_k$ is the singular $k$-sequence $s_k=(k+2,2,\ldots,2)$ and $r_l$ is the regular $l$-sequence $r_l=(2,\ldots,2)$.
 We have
 $$\left(\begin{array}{cc}p&q\\r&s\end{array}\right)=  \left(\begin{array}{cc}0&1\\1&k_1\end{array}\right)\cdots \left(\begin{array}{cc}0&1\\1&k_N\end{array}\right)$$
 and for any $1\le i\le N$ we set
 $$\left(\begin{array}{cc}p_i&q_i\\r_i&s_i\end{array}\right)=\left(\begin{array}{cc}p_i(k_1,\ldots,k_i)&q_i(k_1,\ldots,k_i)\\r_i(k_1,\ldots,k_i)&s_i(k_1,\ldots,k_i)\end{array}\right)=  \left(\begin{array}{cc}0&1\\1&k_1\end{array}\right)\cdots \left(\begin{array}{cc}0&1\\1&k_i\end{array}\right),$$
 therefore
\begin{equation}\label{equation2}
\left(\begin{array}{cc}p_i&q_i\\r_i&s_i\end{array}\right)=\left(\begin{array}{cc}q_{i-1}&p_{i-1}+k_iq_{i-1}\\s_{i-1}&r_{i-1}+k_is_{i-1}\end{array}\right)
\end{equation}

  If $N=1$, dual graph of the curves is
\begin{center}
\includegraphics[width=6cm]{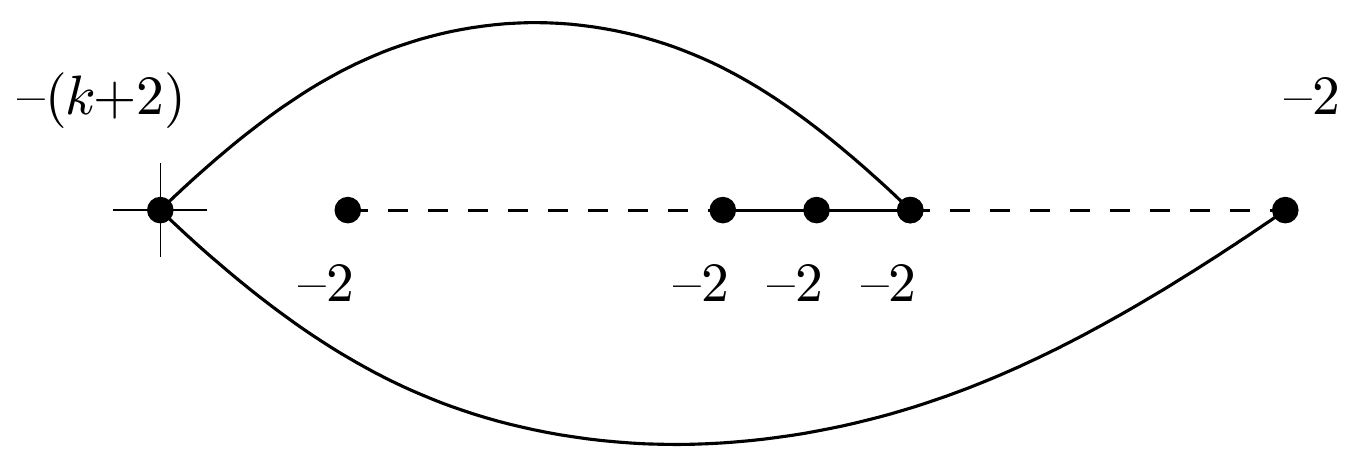}
\end{center}
the (opposite) intersection matrix of the (unique) branch is the matrix  of a chain of length $k$
$$\d_k=\left|\begin{array}{rrrrrr}
2&-1&0&\ldots&\ldots&0\\
-1&2&-1&\ddots&&\vdots\\
0&-1&2&\ddots&\ddots&\vdots\\
\vdots&&\ddots&\ddots&\ddots&0\\
\vdots&&&\ddots&\ddots&-1\\
0&\ldots&\ldots&0&-1&2\\
\end{array}\right|$$
We have $\d_k=k+1$ and by \cite{D5} thm 3.20, $k(S)$ is equal to $\d_k$. Now here
$$\left(\begin{array}{cc}p&q\\r&s\end{array}\right)=  \left(\begin{array}{cc}0&1\\1&k\end{array}\right)$$
therefore the result is checked for $N=1$.\\
If $N=2$, the sequence of opposite self-intersections of the curves in the branch is
$$\underbrace{2\cdots 2}_{k_1-1}\ (k_2+2)$$
On one hand
$$\left(\begin{array}{cc}p&q\\r&s\end{array}\right)=  \left(\begin{array}{cc}0&1\\1&k_1\end{array}\right)\left(\begin{array}{cc}0&1\\1&k_2\end{array}\right)= \left(\begin{array}{cc}1&k_2\\k_1&1+k_1k_2\end{array}\right)$$
On second hand, the order of the (opposite) intersection matrix of the branch is $k_1$. By \cite{D5} thm 3.20, 
$$k(S)=\left|\begin{array}{rrrc}
2&-1\\
-1&\ddots&\ddots\\
&\ddots&2&-1\\
&&-1&k_2+2
\end{array}\right|=k_1k_2+k_1+1=r+s.$$

\begin{itemize}
\item If $N=2\nu$, the sequence of opposite self-intersections of the curves in the branch is
$$\underbrace{2\cdots 2}_{k_1-1}\ (k_2+2)\ \underbrace{2\cdots 2}_{k_3-1}\ \cdots\cdots\cdots \underbrace{2\cdots 2}_{k_{2\nu -1}-1}\ (k_{2\nu}+2)$$
\item If $N=2\nu +1$,  the sequence of opposite self-intersections is
$$\underbrace{2\cdots 2}_{k_1-1}\ (k_2+2)\ \underbrace{2\cdots 2}_{k_3-1}\ \cdots\cdots\cdots \underbrace{2\cdots 2}_{k_{2\nu -1}-1}\ (k_{2\nu}+2)\ \underbrace{2\cdots 22}_{k_{2\nu +1}}$$
\end{itemize}

\begin{itemize}
\item If $N=2\nu$, we have
$$\left(\begin{array}{cc}p&q\\r&s\end{array}\right)=  \left(\begin{array}{cc}0&1\\1&k_1\end{array}\right)\cdots \left(\begin{array}{cc}0&1\\1&k_{2\nu}\end{array}\right)$$

the determinant of the opposite self-intersection matrix of the branch is
$$\d(k_1,\ldots,k_{2\nu})=\left| \begin{array}{cc}\vline \begin{array}{ccc}\hline&&\\&D&\\
&&\\ \hline\end{array}\vline&\begin{array}{l}\\ \\-1\end{array}\\
\begin{array}{rrr}&&\\&&-1\end{array}&\\
&k_{2\nu}+2
\end{array}\right|$$

where $D=D(k_1,\ldots,k_{2\nu-1})$ is the block corresponding to 
$$\underbrace{2\cdots 2}_{k_1-1}\ (k_2+2)\ \underbrace{2\cdots 2}_{k_3-1}\ \cdots\cdots\cdots \underbrace{2\cdots 2}_{k_{2\nu -1}-1}$$

We have by \cite{D5},   the induction hypothesis and relations (\ref{equation2}),
$$\begin{array}{lcl}
k(S)&=&\d(k_1,\ldots,k_{2\nu})=k_{2\nu}\det D(k_1,\ldots,k_{2\nu-1}) + \det D(k_1,\ldots,k_{2\nu-1}+1)\\
&&\\
& = &k_{2\nu}\Bigl(r(k_1,\ldots,k_{2\nu-2},k_{2\nu-1}-1)+s(k_1,\ldots,k_{2\nu-2},k_{2\nu-1}-1)\Bigr)\\
&&\\
&& + r(k_1,\ldots,k_{2\nu-1})+s(k_1,\ldots,k_{2\nu-1})\\
&&\\
&=&k_{2\nu}\Bigl(r(k_1,\ldots,k_{2\nu-2},k_{2\nu-1})+s(k_1,\ldots,k_{2\nu-2},k_{2\nu-1})-s(k_1,\ldots,k_{2\nu-2})\Bigr)\\
&&\\
&& + r(k_1,\ldots,k_{2\nu-1})+s(k_1,\ldots,k_{2\nu-1})\\
&&\\
&=&k_{2\nu}s(k_1,\ldots,k_{2\nu-2},k_{2\nu-1}) + r(k_1,\ldots,k_{2\nu-1})+s(k_1,\ldots,k_{2\nu-1})\\
&&\\
&=&r(k_1,\ldots,k_{2\nu})+s(k_1,\ldots,k_{2\nu})=r+s.
\end{array}$$
\item If $N=2\nu+1$, we follow similar arguments:\\
Let $D$ be the matrix of the chain
$$\underbrace{2\cdots 2}_{k_1-1}\ (k_2+2)\ \underbrace{2\cdots 2}_{k_3-1}\ \cdots\cdots\cdots \underbrace{2\cdots 2}_{k_{2\nu -3}-1}\ (k_{2\nu-2}+2)$$
then by \cite{D5}, $k(S)=\d(k_1,\ldots,k_{2\nu+1})$ and 

$$\d(k_1,\ldots,k_{2\nu+1})=\left| \begin{array}{cl}\vline \begin{array}{ccc}\hline&&\\&D&\\
&&\\ \hline\end{array}\vline&\begin{array}{l}\\ \\-1\end{array}\\
\begin{array}{rrrr}&&&-1\\&&\\&&\\&&\\&&\\&&\\&&\\&&\\&&\\&&\\&&\end{array}&
\vline\begin{array}{ccccccc}\hline &\\  2&-1\\-1&\ddots&\ddots\\&\ddots&2&-1\\&&\ddots&k_{2\nu}+2&\ddots\\&&&-1&2&\ddots\\&&&&\ddots&\ddots&-1\\&&&&&-1&2\\ \hline
\end{array}\vline
\end{array}\right| \begin{array}{l}\scriptstyle 1\\  \\\scriptstyle\sum_{i=1}^{\nu-1}k_{2i-1}\\ \\ \\ \\ \\ \\ \\\scriptstyle\sum_{i=1}^{\nu}k_{2i-1}\\ \\ \\ \\ \scriptstyle\sum_{i=1}^{\nu+1}k_{2i-1}
\end{array}$$

$$\begin{array}{lcl}
&=&k_{2\nu}(k_{2\nu+1}+1)\left| \begin{array}{cl}\vline \begin{array}{ccc}\hline&&\\&D&\\
&&\\ \hline\end{array}\vline&\begin{array}{l}\\ \\-1\end{array}\\
\begin{array}{rrrr}&&&-1\\&&\\&&\\&&\end{array}&
\vline\begin{array}{cccc}\hline &\\  2&-1\\-1&\ddots&\ddots\\&\ddots&\ddots&-1\\&&-1&2\\ \hline
\end{array}\vline
\end{array}\right|  \begin{array}{l}\scriptstyle 1\\   \\ \scriptstyle\sum_{i=1}^{\nu-1}k_{2i-1}\\\scriptstyle\sum_{i=1}^{\nu-1}k_{2i-1}+1 \\ \\ \\ \\   \\ \scriptstyle\sum_{i=1}^{\nu}k_{2i-1}-1\end{array}\\
&&\\
&&+\left| \begin{array}{cl}\vline \begin{array}{ccc}\hline&&\\&D&\\
&&\\ \hline\end{array}\vline&\begin{array}{l}\\ \\-1\end{array}\\
\begin{array}{rrrr}&&&-1\\&&\\&&\\&&\end{array}&
\vline\begin{array}{cccc}\hline &\\  2&-1\\-1&\ddots&\ddots\\&\ddots&\ddots&-1\\&&-1&2\\ \hline
\end{array}\vline
\end{array}\right|  \begin{array}{l}\scriptstyle 1\\   \\ \scriptstyle\sum_{i=1}^{\nu-1}k_{2i-1}\\\scriptstyle\sum_{i=1}^{\nu-1}k_{2i-1}+1 \\ \\ \\ \\   \\ \scriptstyle\sum_{i=1}^{\nu}k_{2i+1}\end{array}\\
\end{array}$$

$$\begin{array}{lcl}
&=&k_{2\nu}(k_{2\nu+1}+1)\d(k_1,\ldots,k_{2\nu-2},k_{2\nu-1}-1)+\d(k_1,\ldots,k_{2\nu-2},k_{2\nu-1}+k_{2\nu+1})\\
&&\\
&=&k_{2\nu}(k_{2\nu+1}+1)\Bigl( r_{2\nu-1}+s_{2\nu-1}-s_{2\nu-2}\Bigr)+s_{2\nu-2}+r_{2\nu-2}\\
&&\\
&&+(k_{2\nu-1}+k_{2\nu+1})s_{2\nu-2}.
\end{array}$$
\end{itemize}
A straightforward computation show that this last expression is equal to $r_{2\nu+1}+s_{2\nu+1}$.

\hfill$\Box$
\begin{Cor} The index of the surface $S(G)$ is
$$Index(S)=\frac{r+s-1}{gcd\{r+s-1,p+q+l-1\} }.$$
\end{Cor}
\begin{Cor} \label{existchamp+} Suppose that $l=d+k(r+s-1)$, then $S$ admits a non trivial global vector field if and only if
$$1-\d (r+s)a_0^{(k+1)r-p+1}=0.$$
\end{Cor}
Proof: If $l=d+k(r+s-1)$, it is easy to check that 
$$u=\frac{p+s+rl-1-\d }{r+s-1}=(k+1)r-p+1.$$
 By propositions \thesection.\ref{CNSchamp} and \thesection.\ref{kSrs},  we have the result.\hfill$\Box$\\

\begin{Not} We denote by $\cal G=\cal G(p,q,r,s,l)$ the family of contracting birational mappings 
$$G(z)=\left(z_1^{p+rl}z_2^{q+sl}+\sum_{i=0}^{l-1}a_i\bigl(z_1^rz_2^s\bigr)^{i+1}+a_{l+K}\bigl(z_1^rz_2^s\bigr)^{l+K+1},z_1^rz_2^s\right),$$
where  $K=\max\left\{0,\left[\frac{l-d}{r+s-1}\right]\right\}$, $a_0\in\bb C^\star$, $a_i\in\bb C$, $i=1,\ldots,l-1,l+K$,  and by $\Phi=\Phi(p,q,r,s,l)$ the group  of  germs of biholomorphisms $\f:(\bb C^2,0)\to (\bb C^2,0)$ for which there exists $G,G'\in \cal G$ such that $G'=\f^{-1}G\f\in\cal G$.   Let
$L:=L(p,q,r,s,l)$ be the group of diagonal linear mappings $\f_{A,B}(z_1,z_2)=(Az_1,Bz_2)$ where $A,B$ satisfy the condition
$$B=A^rB^s,\quad A=A^{p+rl}B^{q+sl}$$
\end{Not}
\begin{Lem} \label{automorphismesdiagonaux} 1) The group $L$ is a subgroup of $\bb U_{p+s+rl-\d -1}\times \bb U_{p+s+rl-\d -1}$, where for any $m\in\bb N^\star$, $\bb U_m$ is the group of $m$-roots of unity.\\
2) The group $L$ operates on $\cal G$; more precisely if $\f_{A,B}\in L$ and
$$G(z)=\left(z_1^{p+rl}z_2^{q+sl}+\sum_{i=0}^{l-1}a_i\bigl(z_1^rz_2^s\bigr)^{i+1}+a_{l+K}\bigl(z_1^rz_2^s\bigr)^{l+K+1},z_1^rz_2^s\right),$$
then
$$G'(z)=\f_{A,B}^{-1}G\f_{A,B}(z)=\left(z_1^{p+rl}z_2^{q+sl}+\sum_{i=0}^{l-1}a'_i\bigl(z_1^rz_2^s\bigr)^{i+1}+a'_{l+K}\bigl(z_1^rz_2^s\bigr)^{l+K+1},z_1^rz_2^s\right),$$
where 
$$Aa'_i=B^{i+1}a_i,\quad{\rm for}\quad i=0,\ldots,l-1,l+K.$$
 In particular $L$ is an abelian subgroup of $\Phi$.
\end{Lem}
The proof is easy and left to the reader.\hfill$\Box$\\

\subsubsection{Moduli spaces of birational mappings}
We want to determine the equivalence classes of the birational mappings $G$, or, that is equivalent, the fibers of the canonical morphism to the OT moduli space. Let 
$$G(z)=\Pi\s(z)=\left(z_1^{p+rl}z_2^{q+sl}+\sum_{i=0}^{l-1}a_i\bigl(z_1^rz_2^s\bigr)^{i+1}+a_{l+K}\bigl(z_1^rz_2^s\bigr)^{l+K+1},z_1^rz_2^s\right),$$
$$ G'(z)=\Pi'\s'(z)=\left(z_1^{p+rl}z_2^{q+sl}+\sum_{i=0}^{l-1}a'_i\bigl(z_1^rz_2^s\bigr)^{i+1}+a'_{l+K}\bigl(z_1^rz_2^s\bigr)^{l+K+1},z_1^rz_2^s\right)$$
be two such birational germs and suppose that there exists  a germ of biholomorphism $\f$ such that
$G'\circ\f=\f\circ G$. Since the degeneration set $\{z_1z_2=0\}$ is invariant and $\f$ cannot swap the rational curves, $\f$ has the form
$$\f(z_1,z_2)=\bigl(Az_1(1+\t(z)),\  Bz_2(1+\mu(z))\bigr).$$
\\
We have
$$\f(G(z))=\left(A\Bigl(z_1^{p+rl}z_2^{q+sl}+\dps\sum_{i\in\{0,\ldots,l-1,l+K\}}a_i\bigl(z_1^rz_2^s\bigr)^{i+1}\Bigr)\bigl(1+\t(G(z))\bigr), Bz_1^rz_2^s\bigl(1+\mu(G(z))\bigr)\right)
$$
$$\begin{array}{lcl}
G'(\f(z))&=&\Bigl(\Bigl[Az_1(1+\t(z))\Bigr]^{p+rl}\Bigl[Bz_2(1+\mu(z))\Bigr]^{q+sl}\\
&& +\dps \sum_{i\in\{0,\ldots,l-1,l+K\}}a'_i\Bigl(Az_1(1+\t(z))\Bigr)^{r(i+1)}     \Bigl(Bz_2(1+\mu(z))\Bigr)^{s(i+1)},\\ 
&&\hfill  A^rB^sz_1^rz_2^s(1+\t(z))^r(1+\mu(z))^s\Bigr)
\end{array}$$
Second members give the equality
$$B\bigl(1+\mu(G(z))\bigr)=A^rB^s(1+\t(z))^r(1+\mu(z))^s \leqno{(II)}$$
Therefore
\begin{equation}\label{II}
B=A^rB^s,\quad{\rm and}\quad 1+\mu(z)=\left(\prod_{j=0}^\infty\Bigl(1+\t(G^j(z))\Bigr)^{r/s^{j+1}}\right)^{-1}.
\end{equation}
First members of the conjugation give
$$\left\{\begin{array}{l}
A\Bigl(z_1^{p+rl}z_2^{q+sl}+\dps\sum_{i\in\{0,\ldots,l-1,l+K\}}a_i\bigl(z_1^rz_2^s\bigr)^{i+1}\Bigr)\bigl(1+\t(G(z))\bigr)\\
\\
\hspace{10mm}= \Bigl(Az_1(1+\t(z))\Bigr)^{p+rl}\Bigl(Bz_2(1+\mu(z))\Bigr)^{q+sl} \\
\\
\hspace{15mm}+\dps \sum_{i\in\{0,\ldots,l-1,l+K\}} a'_i\Bigl(Az_1(1+\t(z))\Bigr)^{r(i+1)}     \Bigl(Bz_2(1+\mu(z))\Bigr)^{s(i+1)}
\end{array}\right.\leqno{(I)}$$
Setting $\d =ps-qr=\pm 1$,  we obtain with (\ref{II}),
$$\left\{\begin{array}{l}
A\Bigl(z_1^{p+rl}z_2^{q+sl}+\dps\sum_{i\in\{0,\ldots,l-1,l+K\}}a_i\bigl(z_1^rz_2^s\bigr)^{i+1}\Bigr)\bigl(1+\t(G(z))\bigr)\\
\\
\hspace{10mm}= A^{p+rl}B^{q+sl}z_1^{p+rl}z_2^{q+sl}(1+\t(z))^{\d /s}\left(\dps\prod_{j=1}^\infty\Bigl(1+\t(G^j(z))\Bigr)^{r/s^{j+1}}\right)^{-(q+sl)}\\
\\ \hspace{15mm}+\dps \sum_{i\in\{0,\ldots,l-1,l+K\}} a'_i B^{i+1}
(z_1^rz_2^s)^{(i+1)} \left(\prod_{j=1}^\infty\Bigl(1+\t(G^j(z))\Bigr)^{r/s^{j+1}}\right)^{-s(i+1)}  
\end{array}\right.\leqno{(I)}$$

\begin{Lem and Def}\label{resonance}  With the previous notations, given $p,q,r,s,l$, the positive integral solutions $(i,j)$ of the system
$$\left\{\begin{array}{lcl}p+rl+i&=&r\g\\q+sl+j&=&s\g\end{array}\right.\leqno{(E)}$$
for which there exists $\g\ge 1$ are all of the form
$$\left\{\begin{array}{lcl}i&=&kr-p\\ j&=&ks-q\end{array}\right., \quad k\ge 1.$$
We have then $\g=k+l$. In particular the least solution is $(r-p,s-q)$. When $(E)$ has a solution we shall say that there is a {\bf resonance}.
\end{Lem and Def}
Proof: We have $\g=l+\frac{p+i}{r}=l+\frac{q+j}{s}$. Since $\frac{p+i}{r}$ and $\frac{q+j}{s}$ are integers, there exists $k,k'\in\bb N$ such that $p+i=kr$ and $q+j=k's$. Moreover $k=\frac{p+i}{r}=\frac{q+j}{s}=k'$ which gives the result. The other assertions are evident.\hfill$\Box$\\

Comparing monomial terms $z_1^{p+rl}z_2^{q+sl}$ in $(I)$ we obtain thanks to lemma \thesection.\ref{resonance}
\begin{equation}
A=A^{p+rl}B^{q+sl}
\end{equation}
By lemma \thesection.\ref{automorphismesdiagonaux},  $A$, $B$ are roots of unity.\\
 Let $Aut(\bb C^2,0)$ be the group of germs of biholomorphisms of $(\bb C^2,0)$ and $Aut(\bb C^2,H,0)$ be the subgroup of $Aut(\bb C^2,0)$ whose germs leave each of the components of the hypersurface $H=\{z_1z_2=0\}$ invariant, i.e. $\f\in Aut(\bb C^2,H,0)$ has the form 
$$\f(z)=\bigl(Az_1(1+\t(z)),Bz_2(1+\mu(z))\bigr).$$
Notice that $\Phi\subset Aut(\bb C^2,H,0)$.\\
Let $Aut_{Id}(\bb C^2,0)$ (resp. $Aut_{Id}(\bb C^2,H,0)$) be the subgroup of $Aut(\bb C^2,0)$ (resp. $Aut(\bb C^2,H,0)$) of germs of biholomorphisms $\f$ tangent to the identity, i.e. $\f\in Aut_{Id}(\bb C^2,H,0)$ if 
$$\f(z)=\bigl(z_1(1+\t(z)),z_2(1+\mu(z))\bigr).$$
\begin{Lem}  Let $\a:Aut_{Id}(\bb C^2,H,0)\to Aut(\bb C^2,H,0)$ the canonical injection and $\b:Aut(\bb C^2,H,0)\to L$ defined by $\b(\f)=\f_{AB}$ the linear part of $\f$. Then, 
$Aut_{Id}(\bb C^2,H,0)$ is a normal subgroup of $Aut(\bb C^2,H,0)$ and we have the exact sequence
$$\{Id\}\to Aut_{Id}(\bb C^2,H,0)\stackrel{\a}{\to}Aut(\bb C^2,H,0)\stackrel{\b}{\to} L\to \{Id\}.$$
\end{Lem}

Replacing $\f$ by $\f\f_{A,B}^{-1}$ we obtain an automorphism  tangent to the identity, therefore we have to determine equivalence classes of the equivalence relation on $\cal G$
$$G\sim G' \quad \Longleftrightarrow \quad \exists\ \f\in Aut_{Id}(\bb C^2,H,0), \quad G'\f=\f G.$$ 
 The equation  $(I)$ becomes

$$\left\{\begin{array}{l}
\Bigl(z_1^{p+rl}z_2^{q+sl}+\dps\sum_{i\in\{0,\ldots,l-1,l+K\}}a_i\bigl(z_1^rz_2^s\bigr)^{i+1}\Bigr)\bigl(1+\t(G(z))\bigr)\\
\\
\hspace{10mm}= z_1^{p+rl}z_2^{q+sl}(1+\t(z))^{\d /s}\left(\dps\prod_{j=1}^\infty\Bigl(1+\t(G^j(z))\Bigr)^{r/s^{j+1}}\right)^{-(q+sl)}\\
\\ \hspace{15mm}+\dps \sum_{i\in\{0,\ldots,l-1,l+K\}} a'_i 
(z_1^r z_2^s)^{i+1} \left(\prod_{j=1}^\infty\Bigl(1+\t(G^j(z))\Bigr)^{r/s^{j+1}}\right)^{-s(i+1)}
\end{array}\right.\leqno{(I)}$$
Notice that we obtain immediately $a_0=a'_0$.\\
The question is to determine the quotient $\cal G/\sim$. We shall see at the end of this section that the equivalence relation is generically trivial.\\

In the following lemma the maximum is due to the fact that we may have $l-d<0$.
\begin{Lem} \label{annulationtheta} Let $\mu=\max\left\{d,l+\dps\left[\frac{l-d}{r+s-1}\right]\right\}$ and $\t(z)=\sum_{i+j\ge 1}t_{ij}z_1^iz_2^j$. If $t_{ij}=0$ for $i+j\le \mu$, then
$$\t=0$$
and $\f$ is linear.
\end{Lem}
Proof: By hypothesis we have
$$\t(G(z))=\left(\sum_{i+j= \mu+1}a_0^it_{ij}\right)(z_1^rz_2^s)^{\mu+1}\quad {\rm mod}\ \frak M^{(r+s)(\mu+1)+1},$$
hence
$$\prod_{j=1}^\infty\Bigl(1+\t(G^j(z))\Bigr)^{r/s^{j+1}}= 1+ \frac{r}{s^2}\left(\sum_{i+j= \mu+1}a_0^it_{ij}\right)(z_1^rz_2^s)^{\mu+1}\quad {\rm mod}\ \frak M^{(r+s)(\mu+1)+1}.$$
We show by induction on $k=i+j\ge \mu+1$ that $t_{ij}=0$.\\
We consider the terms of degree $p+q+(r+s)l+\mu+1$
$$z_1^{p+rl}z_2^{q+sl} \, \frac{\d }{s}\, \sum_{i+j=\mu+1}t_{ij}z_1^iz_2^j$$
and we are looking for other terms of the same degree or bidegrees in $(I)$.

The inequalities
$$\left\{\begin{array}{l}
p+q+(r+s)l+\mu+1<p+q+(r+s)l+(r+s)(\mu+1),\\
\\
p+q+(r+s)l+\mu+1<r+s+(r+s)(\mu+1)
\end{array}\right. $$

show that there is no other term of the same degree when $a_{K+l}=0$. If $a_{K+l}\neq 0$, $l=d+K(r+s-1)$ and it is easy to check that $(r+s)(l+K)\neq p+q+(r+s)l+\mu+1$, hence $t_{ij}=0$ if $i+j=\mu+1$.\\
Suppose that for $k\ge \mu+2$,
$$\t(z)=\sum_{i+j\ge k}t_{ij}z_1^iz_2^j,$$
then the similar inequalities show the result.\hfill$\Box$\\

\begin{Lem} \label{determinationtheta} Let $\mu=\max\left\{d,l+\dps\left[\frac{l-d}{r+s-1}\right]\right\}$.\\
Then, the coefficients $t_{ij}$, for $i+j\le\mu$, with $a_i$ and $a'_i$, $i=0,\ldots,l-1,l+K$, determine uniquely $\t$ hence also $\f$. 
\end{Lem}

Proof:  We show by induction on $k\ge 0$ that the coefficients $t_{ij}$ for $i+j\le\mu$ determine uniquely the coefficients $t_{ij}$ for $i+j\ge \mu+k$. It is sufficient to show  that if the coefficients $t_{ij}$, for $i+j\le \mu+k$ are determined by coefficients $t_{ij}$ for $i+j\le\mu$ then the coefficients $t_{ij}$ for $i+j=\mu+k+1$ are determined by coefficients $t_{ij}$ for $i+j\le \mu+k$. On that purpose we consider homogeneous part of degree $p+q+(r+s)l+\mu+k+1$ in $(I)$ which contains
 the part
$$z_1^{p+rl}z_2^{q+sl} \frac{\d }{s} \left(\sum_{i+j=\mu+k+1}t_{ij}z_1^iz_2^j\right)$$
In order to prove that all other  terms with such degree involve only $t_{ij}$ with $i+j\le\mu+k$, it is sufficient to prove that if $i+j\ge \mu+k+1$ then
$$r+s+(i+j)(r+s)> p+q+(r+s)l+\mu+k+1,$$
and it is sufficient to prove that
$$r+s+(\mu+k+1)(r+s)>p+q+(r+s)l+\mu+k+1.$$
\begin{itemize}
\item If $l\le d$, $\mu=d$, and we have to check that 
$$r+s+(d+k+1)(r+s)>p+q+(r+s)d+d+k+1$$
which is clear;
\item If  $d+K(r+s-1)\le l<d+(K+1)(r+s-1)$, then $\mu=l+K$. We have to check
$$r+s+(l+K+k+1)(r+s)>p+q+(r+s)l+l+K+k+1$$
However this inequality is equivalent to
$$d+(K+k+1)(r+s-1)>l$$
which is satisfied by assumption.
\end{itemize}
\hfill$\Box$\\

\begin{Prop} \label{presque-unicite} Let  $\cal G=\cal G(p,q,r,s,l)$ the family of contracting birational mappings 
$$G(z)=\left(z_1^{p+rl}z_2^{q+sl}+\sum_{i=0}^{l-1}a_i\bigl(z_1^rz_2^s\bigr)^{i+1}+a_{l+K}\bigl(z_1^rz_2^s\bigr)^{l+K+1},z_1^rz_2^s\right),$$
where $K=\max\left\{0,\left[\frac{l-d}{r+s-1}\right]\right\}$, $a_0\in\bb C^\star$, $a_i\in\bb C$, $i=1,\ldots,l-1,l+K$.\\
1) The group $Aut_{Id}(\bb C^2,H,0)\cap\Phi$ is isomorphic to $(\bb C,+)$ and $\t$ is determined by exactly one coefficient of the homogeneous part of degree $l+K$.\\
2) Suppose that  $\frac{l-d}{r+s-1}$ is a non negative integer, i.e. $l=d+K(r+s-1)$, then 
\begin{description}
\item{a)} If there are global vector fields, $Aut_{Id}(\bb C^2,H,0)\cap \Phi $ acts trivially on $\cal G$, in particular $a_{l+K}$ is an effective parameter,
\item{b)} If there are no global vector fields, $Aut_{Id}(\bb C^2,H,0)\cap \Phi$ acts transitively on $\bb C_{a_{l+K}}$, i.e. the complex structure on $S(G)$ does not depend on $a_{l+K}$.
\end{description}
3) Suppose that $\frac{l-d}{r+s-1}$ is not a non negative integer, then $Aut_{Id}(\bb C^2,H,0)\cap\Phi$ acts transitively on $\bb C_{a_{l+K}}$, i.e. the complex structure on $S(G)$ does not depend on $a_{l+K}$.\\ 
\end{Prop}
Proof:  Suppose that $\t\neq 0$ and let  $\g=\min\{i+j\ge 1\mid t_{ij}\neq 0\}$. By lemma \ref{annulationtheta}, $\g\le\mu$. The homogeneous parts of lower degree in $(I)$ which involve $t_{ij}$ with $\g=i+j$  are
\begin{itemize}
\item Case $\g\le l-1$ or $\g=l+K$,
$$a_0z_1^rz_2^s\left(\sum_{i+j=\g}t_{ij}a_0^i\right)(z_1^rz_2^s)^\g + a_\g(z_1^rz_2^s)^{\g+1},\leqno{(A)}$$
$$ z_1^{p+rl}z_2^{q+sl}\frac{\d }{s}\sum_{i+j=\g}t_{ij}z_1^iz_2^j,\leqno{(B)}$$
$$ -\frac{r}{s}a_0z_1^rz_2^s\left(\sum_{i+j=\g}t_{ij}a_0^i\right)(z_1^rz_2^s)^\g + a'_\g(z_1^rz_2^s)^{\g+1}\leqno{(C)}$$
\item Case $\g \ge l$ and $\g\neq l+K$, $(A)$ is replaced by
$$a_0z_1^rz_2^s\left(\sum_{i+j=\g}t_{ij}a_0^i\right)(z_1^rz_2^s)^\g,\leqno{(A')}$$
and $(C)$ by
$$ -\frac{r}{s}a_0z_1^rz_2^s\left(\sum_{i+j=\g}t_{ij}a_0^i\right)(z_1^rz_2^s)^\g\leqno{(C')}$$
\end{itemize}

\begin{itemize}
\item If there is no resonance, the bidegrees of the terms $(A)$ and $(C)$ (resp. $(A')$ and $(C')$) are all distinct of those in $(B)$, therefore we obtain readily 
$$\sum_{i+j=\g}t_{ij}z_1^iz_2^j=0,$$
 hence a contradiction
\item Therefore there is a resonance and there exists a unique coefficient $t_{kr-p,ks-q}\neq 0$ with $k(r+s)-(p+q)=\g$. Then
\begin{itemize}
\item Case $\g\le l-1$ or $\g=l+K$,
$$\begin{array}{l}
a_0z_1^rz_2^s t_{kr-p,ks-q}a_0^{kr-p}(z_1^rz_2^s)^\g+ a_\g(z_1^rz_2^s)^{\g+1}\\
\\
\hspace{5mm}=z_1^{p+rl}z_2^{q+sl}\frac{\d }{s} t_{kr-p,ks-q} z_1^{kr-p}z_2^{ks-q} + a'_\g(z_1^rz_2^s)^{\g+1}-\frac{r}{s}a_0z_1^rz_2^st_{kr-p,ks-q}a_0^{kr-p}(z_1^rz_2^s)^\g
\end{array}$$
\item Case $\g\ge l$, $\g\neq l+K$
$$\begin{array}{l}
a_0z_1^rz_2^s t_{kr-p,ks-q}a_0^{kr-p}(z_1^rz_2^s)^\g\\
\\
\hspace{20mm}=z_1^{p+rl}z_2^{q+sl}\frac{\d }{s} t_{kr-p,ks-q} z_1^{kr-p}z_2^{ks-q} -\frac{r}{s}a_0z_1^rz_2^st_{kr-p,ks-q}a_0^{kr-p}(z_1^rz_2^s)^\g
\end{array}$$
\end{itemize}
\end{itemize}
The equality of degrees implies that $l=d+(k-1)(r+s-1)$, i.e. the surface has twisted vector fields and $\g=l+(k-1)=l+K$ (and the second case never appears). After simplification (recall that $\d=\pm 1$), we obtain 
$$a'_{l+K}=a_{l+K}- t_{kr-p,ks-q}\frac{\d }{s}\Bigl(1-\d (r+s)a_0^{kr-p+1}\Bigr)$$
Let $\frak M=(z_1,z_2)$. Since $\t(z)=0$ mod $\frak M^{l+K}$, $(I)$ gives
$$a'_i=a_i,\quad{\rm for}\quad i=0,\ldots,l-1,$$
therefore, applying corollary 4.\ref{existchamp+},
\begin{itemize}
\item If there are global vector fields, $1-\d (r+s)a_0^{kr-p+1}=0$  and $a'_{l+K}=a_{l+K}$, hence $Aut_{Id}(\bb C^2,H,0)\cap\Phi$ acts trivially;
\item If there are no vector fields, $1-\d (r+s)a_0^{kr-p+1}\neq 0$, and $G_I\cap\Phi$ acts transitively on the line $a_{l+K}\in \bb C$.
\end{itemize}
By lemma 4.\ref{determinationtheta}, $t=t_{(K+1)r-p,(K+1)s-q}\in \bb C$ determines the formal series $\t$. 
It remains to prove that   $\t$ is convergent hence $Aut_{Id}(\bb C^2,H,0)\cap\Phi\simeq\bb C$.
\begin{itemize}
\item If there are global vector fields, there exists a $1$-parameter group of automorphisms, therefore there are such $\t$ and conversely, any $\t$ defines an automorphism of $S(G)$ which is in the identity component of $Aut(S(G))$;
\item If there is no global vector fields, $a_{l+K}$ is a superfluous parameter and all surfaces are isomorphic, therefore there are such isomorphisms.
\end{itemize}

\hfill$\Box$

\subsubsection{The twisting coefficient}

We want to compare birational germs and Favre polynomial germs of the form 
$$F(z_1,z_2)=(\l z_1z_2^\s+P(z_2)+cz_2^{\frac{\s k}{k-1}},z_2^k),\quad P(z_2)=\sum_{i=p+q}^\s b_iz_2^i$$
given in \cite{OT} (see section \ref{sectionOT}). 
The parameter $\l$  determines the twisting coefficient $\kappa$ such that $H^0(S,K_S^{-\mu}\ot L^\kappa)\neq 0$. In birational germs this role is played by the position of the blown-up point $(a_{0},0)$, $a_0\neq 0$, on the root of the branch when there is only one branch. The condition $j<k$ (or $p+q<r+s$ in our notations) implies that the first blowing-up is of the form $(u',v')\mapsto (v',u'v')$ hence we have to consider the germ $\Pi_l\cdots\Pi_{n-1}\bar\s\Pi_0\cdots\Pi_{l-1}$ at the point $(a_{0},0)$. After a change of coordinates $u=z_1+a_{l-1}$, $v=z_2$ we obtain
$$G(z_1,z_2)=\left(\Bigl(z_1z_2^l+\sum_{i=0}^{l-1}a_iz_2^{i+1}+a_{l+K}z_2^{l+K+1}\Bigr)^pz_2^q, \ \Bigl(z_1z_2^l+\sum_{i=0}^{l-1}a_iz_2^{i+1}+a_{l+K}z_2^{l+K+1}\Bigr)^rz_2^s\right)$$
{\bf \boldmath  In case there is no  global vector fields, $a_{l+K}$ is a superfluous parameter, hence we shall suppose that $a_{l+K}=0$}.\\
\begin{Rem} 1) If $index(S)=1$, we have $\l^{-1}=k(S)\k$, i.e. the invariant used here is the inverse of the invariant $\l=\l(S)$ in \cite{DO1}.\\
2) If  $index(S)\neq 1$, $\l(a)$ is defined up to a $(k-1)$-root of unity.
\end{Rem}

\begin{Prop} \label{dependancedekappa} Let $S$ be a surface with GSS associated to the germ
$$G(z_1,z_2)=\left(\Bigl(z_1z_2^l+\sum_{i=0}^{l-1}a_iz_2^{i+1}+a_{l+K}z_2^{l+K+1}\Bigr)^pz_2^q, \ \Bigl(z_1z_2^l+\sum_{i=0}^{l-1}a_iz_2^{i+1}+a_{l+K}z_2^{l+K+1}\Bigr)^rz_2^s\right), \quad a_0\neq 0$$
Let $\mu=index(X)$ be the index of $S$. Then on the corresponding base $B_{J,M}$ of the family $\Phi_{J,M,\s}:\cal S_{J,M,\s}\to B_{J,M}$, the holomorphic function
$$\k=\k_{J,M,\s}:\bb C^\star\to \bb C^\star$$
such that $H^0(S_a,K_{S_a}^{-\mu}\ot L^{\k(a)})\neq 0$ is a monomial holomorphic function of $a_{0}$, where $O_{0}=(a_{0},0)$. More precisely, if $\d=ps-qr$ and $\s=p+q+l-1$,
$$\k_{J,M,\s}(a_{0})=\d^\mu a_{0}^{\mu\left(\frac{r\s}{r+s-1}-p+1\right)},\quad \mu\left(\frac{r\s}{r+s-1}-p+1\right)\in\bb N^\star.$$
In particular $\k$ is surjective.
\end{Prop}
Proof: Let  $\k:=\k_{J,M,\s}$ be the holomorphic function given by \cite{D6} prop.4.24.

Setting 
$$\left\{
\begin{array}{l}
\dps\Bigl(\ \Bigr):=z_1z_2^l+\sum_{i\in\{0,\ldots,l-1,l+K\}}a_iz_2^{i+1} \quad {\rm and}\\
\\
\dps \Bigl[ \ \Bigr]:=\frac{\part}{\part z_2}\Bigl(\ \Bigr)=lz_1z_2^{l-1}+\sum_{i\in\{0,\ldots,l-1,l+K\}}a_i(i+1)z_2^{i}
 \end{array}\right.$$
$$DG(z)=\left(
\begin{array}{cc} p\Bigl(\ \Bigr)^{p-1}z_2^{l+q}&\quad p\Bigl(\ \Bigr)^{p-1}\Bigl[ \ \Bigr]z_2^q+\Bigl(\ \Bigr)^{p}qz_2^{q-1}\\
&\\
r\Bigl(\ \Bigr)^{r-1}z_2^{l+s}&\quad r\Bigl(\ \Bigr)^{r-1}\Bigl[ \ \Bigr]z_2^s+\Bigl(\ \Bigr)^{r}sz_2^{s-1}
\end{array}\right)$$
and
$$\det DG(z)=(ps-qr)\Bigl(\ \Bigr)^{p+r-1}z_2^{l+q+s-1}.$$
Let $\t\in H^0(S,K_S^{-\mu}\ot L^{\k})$, then there exists an invariant germ in a neighbourhood of the origin of the ball still denoted by $\t$ which vanishes only on the curve $\{z_2=0\}$,
$$\t(z)=z_2^\a A(z)\left(\frac{\part}{\part z_1}\wedge\frac{\part}{\part z_1}\right)^{\ot \mu}$$
such that $A(0)\neq 0$. This germ satisfies the condition
$$\t(G(z))=\k\bigl(\det DG(z)\bigr)^\mu \t(z),$$
which is equivalent to
$$\Bigl(\ \Bigr)^{\a r}z_2^{\a s}A(G(z))=\k(ps-qr)^\mu \Bigl(\ \Bigr)^{\mu(p+r-1)}z_2^{\mu(l+q+s-1)+\a}A(z)$$
where $\d:=ps-qr=\pm 1$.
Considering the homogeneous part of lower degree of each member, we obtain
$$\a(r+s-1)=\mu(p+q+l-1+r+s-1)=\mu(\s+r+s-1)\leqno{(1)}$$

$$\k=(ps-qr)^\mu a_{0}^{\a r-\mu(p+r-1)}\leqno{(2)}.$$
By \cite{D6} proposition 4.24, $\k$ vanishes on smaller strata, i.e. when $a_{0}=0$, therefore $\a r-\mu(p+r-1)>0$. We derive the value of $\k$ from $(1)$ and $(2)$.\hfill$\Box$\\
\vspace{2mm}

\subsubsection{Representation of surfaces with one branch and without twisted vector fields by birational germ} 
{\bf\boldmath We suppose that $l-d\not\equiv 0$ mod $r+s-1$, $a_{l+K}=0$}.\\ 
Given $F(z_1,z_2)=\Bigl(\l z_1z_2^\s+\sum_{i=p+q}^\s b_iz_2^i,\  z_2^{r+s}\Bigr)$, there exist germs $G$ which have the same twisting coefficient $\l$ as $F$ by the surjectivity of $\k$ (Prop. 4.\ref{dependancedekappa}). The aim of the sequel of this section is to prove 
\begin{Th} \label{applicationrationnelle} Given $\l$, we choose suitably $a_0\in\bb C^\star$ and $\e$, $\e^{r+s-1}=1$, in such a way that any $G\in \cal G(p,q,r,s,l)$ with parameter $a_0$ is conjugated to a germ $F\in \cal F(\s,k,j)$ with parameter $\l$. Let $\s=p+q+l-1$. Then\\
A) If $r+s-1$ does not divide $l-d$ or $\l\neq 1$ there is a bijective polynomial mapping 
$$\begin{array}{cccc}
f_{a_0,\e}:& \bb C^{l-1}&\longrightarrow&  \bb C^{l-1}\\
&a=(a_1,\ldots,a_{l-1})&\longmapsto&\Bigl(b_{p+q+1}(a),\ldots,b_{p+q+l-1}(a)\Bigr)
\end{array}$$
  such that 
$$G(z_1,z_2)=\left(\Bigl(z_1z_2^l+\sum_{i=0}^{l-1}a_iz_2^{i+1}\Bigr)^pz_2^q, \ \Bigl(z_1z_2^l+\sum_{i=0}^{l-1}a_iz_2^{i+1}\Bigr)^rz_2^s\right)$$
 is conjugated to the polynomial germ
$$F(z_1,z_2)=\Bigl(\l z_1z_2^\s+\sum_{i=p+q}^\s b_iz_2^i,\  z_2^{r+s}\Bigr),$$
where $\l$ as $\k$ depend only on $a_0$ by 4.\ref{dependancedekappa} and $b_{p+q}=1$.\\
B) If $l-d=K(r+s-1)$ and $\l=1$,  there is a bijective polynomial mapping 
$$\begin{array}{cccc}
f_{a_0,\e}:& \bb C^{l-1}\times \bb C&\longrightarrow&  \bb C^{l-1}\times \bb C\\
&a=(a_1,\ldots,a_{l-1},a_{l+K})&\longmapsto&\Bigl(b_{p+q+1}(a),\ldots,b_{p+q+l-1}(a),c(a)\Bigr)
\end{array}$$
such that
$$G(z_1,z_2)=\left(\Bigl(z_1z_2^l+\sum_{i=0}^{l-1}a_iz_2^{i+1}+a_{l+K}z_2^{l+K+1}\Bigr)^pz_2^q, \ \Bigl(z_1z_2^l+\sum_{i=0}^{l-1}a_iz_2^{i+1}+a_{l+K}z_2^{l+K+1}\Bigr)^rz_2^s\right)$$
 is conjugated to the polynomial germ
 $$F(z_1,z_2)=\left(\l z_1z_2^\s+\sum_{k=p+q}^\s b_kz_2^k+ cz_2^{\frac{\s k(S)}{k(S)-1}}, z_2^{r+s}\right),$$
 where $b_{p+q}=1$.
\end{Th}
Proof: Let $\f(z)=\bigl(\f_1(z),Cz_2(1+\mu(z))\bigr)$ be a germ of biholomorphic map which preserves the degeneration set $\{z_2=0\}$.\\
$$\f(G(z))=\left(\f_1\Bigl(G(z)\Bigr),\ C\left\{z_1z_2^l+\sum_{i=0}^{l-1}a_iz_2^{i+1}\right\}^rz_2^s\Bigl(1+\mu(G(z)\Bigr)\right),$$
$$F(\f(z))=\left(\l \f_1(z)C^\s z_2^\s(1+\mu(z))^\s+\sum_{k=p+q}^\s b_kC^kz_2^k(1+\mu(z))^k,\ C^{r+s}z_2^{r+s}(1+\mu(z))^{r+s}\right).$$
Comparing the right members we have
$$\left\{z_1z_2^{l-1}+\sum_{i=0}^{l-1}a_iz_2^{i}\right\}^r\Bigl(1+\mu(G(z)\Bigr)= C^{r+s-1}(1+\mu(z))^{r+s}.\leqno{(II)}$$
Constant parts give the condition
\begin{equation} \label{a0C}
a_0^r=C^{r+s-1}
\end{equation}
therefore $C$ is determined up to a root of unity $\e$ such that $\e^{r+s-1}=1$. In other terms if we choose a local determination of the $(r+s-1)$-root $a_0^{1/(r+s-1)}$,
\begin{equation} \label{C}
C=\e a_0^{r/(r+s-1)},\quad \e^{r+s-1}=1.
\end{equation}
Moreover the equation
$$\left\{1+\frac{1}{a_0}\left(\sum_{i=1}^{l-1}a_iz_2^{i}+z_1z_2^{l-1}\right)\right\}^r\Bigl(1+\mu(G(z)\Bigr)= (1+\mu(z))^{r+s}.\leqno{(II)}$$
 has the solution
$$1+\mu(z) = \prod_{j=0}^\infty\left\{1+\frac{1}{a_0}\left(\sum_{i=1}^{l-1}a_i(G_2^j(z))^i+G_1^j(z)(G_2^j(z))^{l-1}\right)\right\}^{\frac{r}{(r+s)^{j+1}}}$$

Left members give the equality
$$\left\{\begin{array}{l}
\dps\f_1\left(a_0^p\left\{1+\frac{1}{a_0}\left(\sum_{i=1}^{l-1}a_iz_2^{i}+z_1z_2^{l-1}\right)\right\}^pz_2^{p+q}, \ a_0^r\left\{1+\frac{1}{a_0}\left(\sum_{i=1}^{l-1}a_iz_2^{i}+z_1z_2^{l-1}\right)\right\}^rz_2^{r+s}\right)\\
\\
\hspace{20mm}\dps = \l \f_1(z)\Bigl(C z_2(1+\mu(z)\Bigr)^{p+q+l-1}+\sum_{k=p+q}^{p+q+l-1} b_kC^kz_2^k(1+\mu(z))^k
\end{array}\right.\leqno{(I)}$$
We want to express the coefficients $b_k$ with the $a_i$'s, however the coefficients $A_{ij}$ of the series
$$\f_1(z_1,z_2)=\sum_{i,j}A_{ij}z_1^iz_2^j$$
depend also on $a_i$'s. For example, considering homogeneous parts of bidegree $(0,p+q)$, we have,
$$A_{10}a_0^p=b_{p+q}C^{p+q}=C^{p+q}\leqno{(R_{0})}$$
hence with (\ref{a0C}),
\begin{equation}\label{A10} \dps A_{10}=\e a_0^{\frac{p-\d}{r+s-1}}.\end{equation}
If $p>0$, $r+s>p+q+1$ and $l\ge 2$, homogeneous part of bidegree $(0,p+q+1)$ gives
$$A_{10}a_0^{p-1}pa_1=b_{p+q+1}C^{p+q+1}+C^{p+q}\frac{(p+q)r}{r+s}\frac{a_1}{a_0}$$
therefore by $(R_{0})$,
$$b_{p+q+1}=\frac{\d a_1}{Ca_0(r+s)}.\leqno{(R_1)}$$
Comparing terms of bidegree $(1,p+q+l-1)$ we obtain
$$A_{10}pa_0^{p-1}=\l A_{10}C^{p+q+l-1}+\frac{r(p+q)}{r+s}\frac{C^{p+q}}{a_0}$$
therefore with $(\ref{A10})$ and (\ref{C}), and since $k=k(S)=r+s$,
\begin{equation}\label{lambda}
\l =\frac{\d}{\e^{\s}k}\, a_0^{p-1-\frac{r\s}{k-1}}.
\end{equation}
where $\d =ps-qr$ (with formula $\l^{-1}=\k k$, notice that we recover  the result of Prop 4.\ref{dependancedekappa}) .\\

In order to express the coefficients $b_{p+q+j}$, $j\ge 1$, as polynomials of variables $a_1,\ldots,a_{l-1}$, it is also necessary to express  the coefficients $A_{ij}$ involved in the relations as polynomials of the same variables $a_1,\ldots,a_{l-1}$. Therefore we have to determine the set of points $(i,j)\in\bb N\times\bb N$ which occur as indices of the $A_{ij}$'s in the relations. \\
Let $E_0$ be the subset of indices $(i,j)$ which occur in homogeneous part of bidegree $(0,k)$ for $p+q\le k\le p+q+l-1$ in equation $(I)$. We have
$$E_0=\{(i,j)\mid p+q\le i(p+q)+j(r+s)\le p+q+l-1\}$$

Then we define a translation
$$T(i,j)=(i,j+p+q+l-1)$$
and we want to determine which coefficients $A_{\a\b}$ are involved on the homogeneous part of bidegree $T(i,j)$. On that purpose we define a sequence $(E_m)_{m\ge -1}$ of increasing subsets of $\bb N\times\bb N$, starting with $E_{-1}=\emptyset$,
$$E_m=\left\{(i,j)\mid  i(p+q)+j(r+s)\le (p+q+l-1)\left(1+\frac{1}{r+s}+\cdots+\frac{1}{(r+s)^m}\right)\right\}, \quad m\ge 0,$$
and 
$$E_\infty:=\left\{(i,j)\mid  i(p+q)+j(r+s)< (p+q+l-1)\frac{r+s}{r+s-1}\right\}.$$

{\bf For any polynomial expression $Q(z_1,z_2)$ we denote by $\la Q(z_1,z_2)\ra_{a,b}$ the homogeneous part of bidegree $(a,b)$}.\\
\begin{Lem} \label{populationdecoefficients} Suppose $l-d\not\equiv 0$ mod $r+s-1$. Let $(i,j)\in E_m$, $m\ge 0$.\\
1)  If $i\ge 2$ then for any $(\a,\b)$, the homogeneous parts of bidegree $(i,j+p+q+l-1)$
satisfy
$$\left\la A_{\a\b}\ a_0^{p\a+r\b}\left\{1+\frac{1}{a_0}\left(\sum_{i=1}^{l-1}a_iz_2^{i}+z_1z_2^{l-1}\right)\right\}^{p\a+r\b}z_2^{\a(p+q)+\b(r+s)}\right\ra_{i,j+p+q+l-1}= 0.$$
2) If  $i=1$, and homogeneous part of bidegree $(i,j+p+q+l-1)$
satisfies
$$\left\la A_{\a\b}\ a_0^{p\a+r\b}\left\{1+\frac{1}{a_0}\left(\sum_{i=1}^{l-1}a_iz_2^{i}+z_1z_2^{l-1}\right)\right\}^{p\a+r\b}z_2^{\a(p+q)+\b(r+s)}\right\ra_{i,j+p+q+l-1}\neq 0$$
 then $(\a,\b)\in E_1$.\\
 Moreover
 \begin{itemize}
 \item If $m=0$, then $(\a,\b)\in E_0$,
 \item If $\a=1$, then $\b\le j/(r+s)$, in particular if $j\neq 0$, then $\b\neq j$,
 \item If $\a=0$, then $\b<j$ or $\Bigl\{(i,j)=(1,1)\ {\rm and}\  (\a,\b)=(0,1)\Bigr\}$.
 \end{itemize}
3) If $i=0$ and
$$\left\la A_{0\b}\ a_0^{r\b}\left\{1+\frac{1}{a_0}\left(\sum_{i=1}^{l-1}a_iz_2^{i}+z_1z_2^{l-1}\right)\right\}^{r\b}z_2^{\b(r+s)}\right\ra_{0,j+p+q+l-1}\neq 0,$$
the following two conditions cannot be fulfilled  at the same time
\begin{itemize}
\item $i=\a=0$, and  $j=\b$,
 \item $\b(r+s)=j+p+q+l-1$,
 \end{itemize}
  i.e. if the coefficient $A_{0,j}$ appears two times when considering homogeneous part of bidegree $(0,j+p+q+l-1)$, one occurrence is multiplied by a non constant polynomial in $a_1,\ldots,a_j$.\\
4) If $(i,j)\in E_m$ and $i\ge 2$, then $A_{ij}=0$.\\
5) If  $(0,j)\in E_m\setminus E_{m-1}$, $m\ge 0$ and $(\a,\b)$  satisfies 
$$ \a(p+q)+\b(r+s)=j+p+q+l-1,$$
 then 
 \begin{itemize}
 \item $(\a,\b)\in E_{m+1}\setminus E_m$
 \item $\a=0$ or $\a=1$ and $(\a,\b)$ is unique.
 \end{itemize}
 In other words, in homogeneous part of bidegree $(0,j+p+q+l-1)$, there are, modulo $\frak M=(a_1,\ldots,a_{l-1})$, at most two coefficients which occur:  $A_{0,j}$ and perhaps another $A_{\a\b}$ with $\a=0$ or $\a=1$.
\end{Lem}

Proof:  If 
$$\left\la A_{\a\b}\ a_0^{p\a+r\b}\left\{1+\frac{1}{a_0}\left(\sum_{i=1}^{l-1}a_iz_2^{i}+z_1z_2^{l-1}\right)\right\}^{p\a+r\b}z_2^{\a(p+q)+\b(r+s)}\right\ra_{i,j+p+q+l-1}\neq 0$$
then, $p\a+r\b\ge i$, and the least degree in $z_2$ is
$$(l-1)i+\a(p+q)+\b(r+s).$$
Since $\la .\ra_{i,j+p+q+l-1}\neq 0$, 
$$(l-1)i+\a(p+q)+\b(r+s)\le j+p+q+l-1\leqno{(\ast)}$$
however by assumption $(i,j)\in E_m\subset E_\infty$,
$$j< \frac{p+q+l-1}{r+s-1}-i\ \frac{p+q}{r+s}$$
therefore
$$(l-1)i+i\ \frac{p+q}{r+s}+\a(p+q)+\b(r+s)< (p+q+l-1)\left(1+\frac{1}{r+s-1}\right)\leqno{(\ast\ast)}$$
1) If $i\ge 2$, we have, by $(\ast\ast)$, 
$$2(l-1)+2\ \frac{p+q}{r+s}+\a(p+q)+\b(r+s)< (p+q+l-1)\left(1+\frac{1}{r+s-1}\right)$$
Since $p+q<r+s$,
$$(l-1)\left(1-\frac{1}{r+s-1}\right)+(p+q)\left(\frac{2}{r+s}-\frac{1}{r+s-1}+\a+\b-1\right)< 0$$
which is impossible.\\
2) Suppose that $(i,j)\in  E_m$, and $i=1$  then 
$$j< \frac{p+q+l-1}{r+s}\left(1+\frac{1}{r+s-1}\right)-\frac{p+q}{r+s}$$
and by ($\ast$)
$$(l-1)+\frac{p+q}{r+s}+\a(p+q)+\b(r+s)< (p+q+l-1)\left(1+\frac{1}{r+s}+\frac{1}{(r+s)(r+s-1)}\right)$$
which is equivalent to
$$\begin{array}{l}
\dps(l-1)\left(1-\frac{1}{(r+s)(r+s-1)}\right)+(p+q)\left(\frac{1}{r+s}-\frac{1}{(r+s)(r+s-1)}\right)+\a(p+q)+\b(r+s)\\
\\
\dps\hfill< (p+q+l-1)\left(1+\frac{1}{r+s}\right)
\end{array}$$
hence 
$$\a(p+q)+\b(r+s)\le (p+q+l-1)\left(1+\frac{1}{r+s}\right)$$
and $(\a,\b)\in E_1$.\\
If $m=0$, the result derives from the definition of $E_0$ and $(\ast)$.\\
If in ($\ast$), $\a=1$, $\b(r+s)\le j$.\\
If in ($\ast$), $\a=0$, $\b(r+s)\le j+(p+q)$. If moreover $\b\ge j$, then
\begin{itemize}
\item $j=0$ and $\b(r+s)\le p+q$ which is impossible because $\a=0$,
\item $j=1$, and $\b=1$.
\end{itemize}
3) Let $(0,j)\in E_m$, with $m\ge 0$ is minimal. We have
$$j(r+s)\le (p+q+l-1)\left(1+\frac{1}{r+s}+\cdots+\frac{1}{(r+s)^m}\right).$$
 If $j(r+s)=j+p+q+l-1$, then
$$j\le (p+q+l-1)\left(\frac{1}{r+s}+\cdots+\frac{1}{(r+s)^m}\right)$$
hence
$$j(r+s)\le (p+q+l-1)\left(1+\frac{1}{r+s}+\cdots+\frac{1}{(r+s)^{m-1}}\right).$$
and $(0,j)\in E_{m-1}$ which is contradictory.\\
4) Let $(i,j)\in E_m$ with $i\ge 2$ and consider part of bidegree $(i,j+p+q+l-1)$. By 1), left member of $(I)$ gives no contribution, we show now that
$$\left\la b_kC^kz_2^k(1+\mu(z))^k\right\ra_{i,j+p+q+l-1}=0.$$
In fact, the monomials which contain $z_1^i$ contain $z_2$ at the power at least $k+i(l-1)$ with $p+q\le k\le p+q+l-1$ and it is sufficient to show that 
$$j+p+q+l-1<k+i(l-1).$$
Moreover,  $k\ge p+q$ and $i\ge 2$, hence it is sufficient to prove that $j+p+q+l-1<p+q+2(l-1)$, i.e.
$$j<l-1.\leqno{(\maltese)}$$
By assumption, $(i,j)\in E_m$, therefore
$$j\le \frac{1}{r+s}(p+q+l-1)\left(1+\cdots+\frac{1}{(r+s)^m}\right)-i\frac{p+q}{r+s}$$
and  condition $(\maltese)$ is satisfied if
$$\frac{1}{r+s}(p+q+l-1)\left(1+\cdots+\frac{1}{(r+s)^m}\right)<(l-1) +2\frac{p+q}{r+s}$$
which is clearly satisfied since $r+s\ge 2$.
Finally we obtain
$$0=\l A_{ij}C^{p+q+l-1}z_1^iz_2^{j+p+q+l-1}$$
and $A_{ij}=0$.\\
5) If $(0,j)\in E_m\setminus E_{m-1}$, 
$$\begin{array}{l}
\dps(p+q+l-1)\left(1+\frac{1}{r+s}+\cdots+\frac{1}{(r+s)^{m-1}}\right)<j(r+s)\\
\\
\hspace{60mm} \dps\le (p+q+l-1)\left(1+\frac{1}{r+s}+\cdots+\frac{1}{(r+s)^m}\right)
\end{array}$$
By hypothesis,  $(\a,\b)$ satisfies
$$\a(p+q)+\b(r+s)= j+p+q+l-1$$
therefore the following holds
$$\begin{array}{l}
\dps(p+q+l-1)\left(1+\frac{1}{r+s}+\cdots+\frac{1}{(r+s)^{m}}\right)<\a(p+q)+\b(r+s)\\
\\
\hspace{60mm} \dps\le (p+q+l-1)\left(1+\frac{1}{r+s}+\cdots+\frac{1}{(r+s)^{m+1}}\right)
\end{array}$$
i.e. $(\a,\b)\in E_{m+1}\setminus E_m$. By 4), $\a=0$ or $\a=1$.\\
If $k=(p+q)+\b(r+s)=\b'(r+s)$ then $p+q$ is a multiple of $r+s$ which is impossible, therefore we have the unicity of $(\a,\b)$.
\hfill$\Box$\\

\begin{Lem}\label{frationnelle} The linear system with coefficients in $\bb C[a_1,\ldots,a_{l-1}]$ and unknowns  $$b_{p+q+1},\ldots,b_{p+q+l-1}\quad {\rm and} \quad A_{ij}, \quad (i,j)\in E_\infty$$ is a Cramer system of order $l-1+\Card(E_\infty)$. More precisely, modulo $\frak M$, its determinant is
$$\D=C^{p+q+1}\cdots C^{p+q+l-1}(\l C^{p+q+l-1})^{\Card E_\infty}\neq 0 \quad {\rm mod}\ \frak M:=(a_1,\ldots,a_{l-1})$$
and $b_k=\frac{B_k}{\D}$, $k=p+q+1,\ldots,p+q+l-1$, $A_{ij}=\frac{B_{ij}}{\D}$, $(i,j)\in E_\infty$, with $B_k,B_{ij}\in \bb C[a_1,\ldots,a_{l-1}]$. In particular, $f_{a_0,\e}$ is rational.
\end{Lem}
Proof: We order the unknowns in the following way: First unknowns $b_{p+q+1},\ldots,b_{p+q+l-1}$, after coefficients $A_{0j}\neq 0$, with $(0,j)\in E_0$ then $(0,j)\in E_1\setminus E_0$, \ldots $(0,j)\in E_{m+1}\setminus E_m$, exhausting $E_\infty$. Finally coefficients $A_{1j}$, with $j$ in the decreasing order. We have the same number of equations and of unknowns, therefore we have a linear system of order $l-1+\Card(E_\infty)$. Let $\frak M=(a_1,\ldots,a_{l-1})$.  In order to prove that we have a Cramer system it is sufficient to prove that modulo $\frak M$ the determinant $\D$ is nonzero. Therefore we consider the equation $(I)$ modulo $\frak M$, i.e.
$$\left\{\begin{array}{l}
\dps\f_1\left(a_0^p\left\{1+\frac{z_1z_2^{l-1}}{a_0}\right\}^pz_2^{p+q}, \ a_0^r\left\{1+\frac{z_1z_2^{l-1}}{a_0}\right\}^rz_2^{r+s}\right)\\
\\
\hspace{20mm}\dps = \l \f_1(z)\Bigl(C z_2(1+\mu(z)\Bigr)^{p+q+l-1}+\sum_{k=p+q}^{p+q+l-1} b_kC^kz_2^k(1+\mu(z))^k \qquad {\rm mod}\ \frak M
\end{array}\right.\leqno{(I_{\frak M})}$$
where, in the infinite product $1+\mu(z)$, 
$$\begin{array}{lcl}
G(z_1,z_2)&=&\Bigl((z_1z_2^l+a_0z_2)^pz_2^{q}, \ (z_1z_2^l+a_0z_2)^rz_2^{s}\Bigr)\\
&&\\
&=& \dps\left(a_0^p\left(1+\frac{z_1z_2^{l-1}}{a_0}\right)^pz_2^{p+q}, \ a_0^r\left(1+\frac{z_1z_2^{l-1}}{a_0}\right)^rz_2^{r+s}\right) \quad {\rm mod}\  \frak M
\end{array}$$
which provides
$$\begin{array}{lcl}
1+\mu(z)&=&\dps\left\{1+\frac{z_1z_2^{l-1}}{a_0}\right\}^\frac{r}{r+s}\left\{1+\frac{a_0^{p+r(l-1)}\left(1+\frac{z_1z_2^{l-1}}{a_0}\right)^{p+r(l-1)}z_2^{p+q+(r+s)(l-1)}}{a_0}\right\}^{\frac{r}{(r+s)^2}}\cdots\\
&&\\
&=&\dps1+\frac{r}{r+s}\frac{z_1z_2^{l-1}}{a_0}+\frac{r a_0^{p+r(l-1)-1}}{(r+s)^2}z_2^{p+q+(r+s)(l-1)}\\
&&\\
&&\hspace{3mm}+\dps\frac{r[p+r(l-1)]a_0^{p+r(l-1)-2}}{(r+s)^2}z_1z_2^{p+q+(r+s+1)(l-1)}+\cdots\quad {\rm mod}\  \frak M
\end{array}$$
in particular, in the development of $1+\mu(z)$ the least degree in $z_2$ is $p+q+(r+s)(l-1)$.\\

By construction, the diagonal of the matrix is
$$C^{p+q+1},\ldots,C^{p+q+l-1},\l C^{p+q+l-1}, \ldots,\l C^{p+q+l-1}$$
and the square submatrix, of order $l-1$ corresponding to the unknowns
$$b_{p+q+i}, \quad i=1,\ldots,l-1,$$
 is diagonal because $p+q+(r+s)(l-1)>p+q+l-1$ and no term comes from $(1+\mu(z))$.\\
We shall show that after some linear combinations of the lines, we obtain an upper triangular matrix, which yieds $\D\neq 0$.\\
Let $(0,j)\in E_m\setminus E_{m-1}$ ($E_{-1}:=\emptyset$). Since $A_{0,j}\neq 0$,  the homogeneous part of bidegree $(0,j+p+q+l-1)$ is by lemma \ref{populationdecoefficients}, 5)
$$A_{\a\b}a_0^{p\a+r\b}z_2^{\a(p+q)+\b(r+s)}=\l A_{0j}z_2^j(Cz_2)^{p+q+l-1},\quad {\rm mod}\ \frak M$$
with $(\a,\b)\in E_{m+1}\setminus E_m$, if such $(\a,\b)$ exists, or
$$0=\l A_{0j}z_2^j(Cz_2)^{p+q+l-1},\quad {\rm mod}\ \frak M$$
otherwise. A term $b_iz_2^i(1+\mu(z))^i$ has no part of homogeneous bidegree $(0,m)$ because $j+p+q+l-1<2(p+q)+(r+s)(l-1)$. Therefore, with the chosen order on the unknowns, all coefficients of the linear equation are \underline{over} the diagonal of the matrix.\\
 Remain homogeneous parts of bidegree $(1,j+p+q+l-1)$ involving $A_{1,j}$ for $j\ge 1$. We have
$$\begin{array}{lcl}
(1+\mu(z))^i&=&\dps1+\frac{ir}{r+s}\frac{z_1z_2^{l-1}}{a_0}+\frac{ir a_0^{p+r(l-1)-1}}{(r+s)^2}z_2^{p+q+(r+s)(l-1)}\\
&&\\
&&\hspace{3mm}+\dps\frac{ir[p+r(l-1)]a_0^{p+r(l-1)-2}}{(r+s)^2}z_1z_2^{p+q+(r+s+1)(l-1)}+\cdots\quad {\rm mod}\  \frak M
\end{array}$$

It is easy to check that for $i\ge p+q$,  
$$i+p+q+(r+s+1)(l-1)>j+(p+q+l-1),$$
 therefore
the only terms which may be involved in homogeneous part of bidegree $(1,j+p+q+l-1)$ are
$$b_iC^iz_2^i\frac{ir}{r+s}\frac{z_1z_2^{l-1}}{a_0}, \quad {\rm where}\quad p+q\le i\le p+q+l-1$$
 therefore 
 $$i=j+p+q.$$
 We have still to check that $j\le l-1$. Since $j\ge 1$, $l\ge 2$. If $(1,j)\in E_\infty$, then 
 $$p+q+(r+s)j\le (p+q+l-1)\frac{r+s}{r+s-1}$$
 This equivalent to
 $$j\le \frac{p+q}{(r+s-1)(r+s)}+\frac{l-1}{r+s-1}$$
As
$$ \frac{p+q}{(r+s-1)(r+s)}+\frac{l-1}{r+s-1}\le \frac{1}{r+s}+\frac{l-1}{r+s-1}$$
taking if necessary the integral part of the last member, the inequality $j\le l-1$ is still fulfilled.\\
 
Now, there are two possibilities\\
\begin{enumerate}
\item There is no $(\a,\b)$ such that $\a(p+q)+\b(r+s)=j+p+q$.  Therefore 
$$0=\l A_{1j}z_1z_2^j(Cz_2)^{p+q+l-1} + b_{p+q+j}C^{p+q+j}z_2^{p+q+j}\frac{(p+q+j)r}{r+s}\frac{z_1z_2^{l-1}}{a_0}\quad {\rm mod}\  \frak M$$

The $j$-th equation (which gives the $j$-th line $L_j$ of the matrix) is
$$0=b_{p+q+j}C^{p+q+j} \quad {\rm mod}\ \frak M$$
therefore substracting $\frac{(p+q+j)r}{a_0(r+s)}L_j$ we remove the coefficient $b_{p+q+j}C^{p+q+j}\frac{(p+q+j)r}{a_0(r+s)}$ which was under the diagonal.

\item There exists $(\a,\b)$ such that $\a(p+q)+\b(r+s)=p+q+j$. By lemma \ref{populationdecoefficients}, 4),  there is at most two such coefficients $(0,\b)$ and $(1,\b')$. By the choice of the ordering, and lemma \ref{populationdecoefficients}, 2),  $A_{1,\b'}>A_{1j}$ and the coefficient of $A_{1,\b'}$
$$-a_0^{p + r\b'-1}(p+r\b')$$
 is over the diagonal.\\
Then mod $\frak M$, the homogeneous part of bidegree $(1,j+p+q+l-1)$ is
$$\begin{array}{l}
\dps A_{0\b}\,a_0^{r\b}\,r\b\, \frac{z_1z_2^{l-1}}{a_0}z_2^{\b(r+s)}+A_{1\b'}\,a_0^{p+r\b'}\,(p+r\b')\, \frac{z_1z_2^{l-1}}{a_0}z_2^{(p+q)+\b'(r+s)}\\
\\
\hspace{25mm}\dps =\l A_{1j}z_1z_2^j(Cz_2)^{p+q+l-1} + b_{p+q+j}C^{p+q+j}z_2^{p+q+j}\, \frac{(p+q+j)r}{r+s}\,\frac{z_1z_2^{l-1}}{a_0}
\end{array}$$
hence
$$A_{0\b}\,a_0^{r\b-1}\, r\b +A_{1\b'}\,a_0^{p+r\b'-1}\,(p+r\b') =\l A_{1j}C^{p+q+l-1}+ b_{p+q+j}C^{p+q+j}\frac{(p+q+j)r}{r+s}\frac{1}{a_0}$$
where perhaps one of the coefficients $A_{0\b}=0$ or $A_{1\b'}=0$.
The $j$-th equation derived from the homogeneous part of bidegree $(0,p+q+j)$ is
$$A_{0\b}a_0^{r\b}+A_{1\b'}a_0^{p + r\b'}=b_{p+q+j}C^{p+q+j} \quad {\rm mod}\ \frak M$$
 and if $A_{0\b}=0$, it remains to substract $\frac{(p+q+j)r}{a_0(r+s)}L_j$ to obtain a triangular matrix.
If $A_{0\b}\neq 0$, we have two coefficients under the diagonal: $A_{0\b}$ and $b_{p+q+j}$. However (miracle !)
$$r\b=\frac{(p+q+j)r}{r+s}$$
therefore substrating $\frac{(p+q+j)r}{a_0(r+s)}L_j=\frac{r\b}{a_0}L_j$ we remove both coefficients, obtaining the desired upper triangular matrix. 
\end{enumerate}
We conclude that $\D=C^{p+q+1}\cdots C^{p+q+l-1}(\l C^{p+q+l-1})^{\Card E_\infty}\neq 0$.
The second member of the Cramer system is nonzero and involves $A_{10}$ and $b_{p+q}=1$, therefore solutions of the system are rational fractions in variables $a_1,\ldots,a_{l-1}$.
\hfill$\Box$\\

Consider the restriction of the equivalence relation defined by $L$. Since $a_0$ is fixed, lemma 4.\ref{automorphismesdiagonaux} shows that we have the extra condition $A=B$ and
$$a=(a_1,\ldots,a_{l-1})\sim a'=(a'_1,\ldots,a'_{l-1}) \Longleftrightarrow  a'_i=B^{i}a_i,\quad{\rm for}\quad i=1,\ldots,l-1,l+K,$$
where
$$B^{k-1}=B^{r+s-1}=1$$
and 
$$B^\s=B^{p+q+l-1}=B^{p+q+(r+s)l-1}=1.$$
Let $\Pi_L:\bb C^{l-1}\to \bb C^{l-1}/L$  be the canonical mapping (when there are twisted vector fields, $\s=(K+1)(r+s-1)$ and $L=\bb Z_{k-1}$). Similarly, consider the restriction to $\bb C^{l-1}$ of the equivalence relation of Favre germs given by lemma 3.\ref{condconj}. We have $\e^{k-1}=1$ and if we fix $\l$ (recall that by Prop.  \ref{dependancedekappa}, $\l$ depends on $a_0$) then $\e^{p+q+l-1}=\e^\s=1$ and 
$$b=(b_{p+q+1},\ldots,b_{p+q+l-1})\sim b'=(b'_{p+q+1},\ldots,b'_{p+q+l-1}) \Longleftrightarrow b'_{p+q+i}=\e^i b_{p+q+i},\hfill 1\le i\le l-1.$$
We see that the equivalence relations $a\sim a'$ and $b\sim b'$ on $\bb C^{l-1}$ are equal.

\subsubsection{Explicit construction of the isomorphic polynomial mapping (no global twisted vector fields)}
In this section we show that $f_{a_0,\e}$ is polynomial. We still suppose that there is .
\begin{Prop}\label{applicationpolynomiale} We choose $a_0\in\bb C^\star$ and $\e$ such that $\e^{r+s-1}=1$. Let  $\s=p+q+l-1$ and suppose that $r+s-1$ does not divide $l-d$. Then there is a  bijective triangular polynomial mapping 
$$\begin{array}{cccc}
f_{a_0,\e}:& \bb C^{l-1}&\longrightarrow&  \bb C^{l-1}\\
&a=(a_1,\ldots,a_{l-1})&\longmapsto&\Bigl(b_{p+q+1}(a),\ldots,b_{p+q+l-1}(a)\Bigr)
\end{array}$$
  such that for $j=1,\ldots,l-1$,
  $$b_{p+q+j}(a)=\frac{\d \, a_j}{C^j(r+s)a_0}+R_j(a_1,\ldots,a_{j-1}),$$
$$G(z_1,z_2)=\left(\Bigl(z_1z_2^l+\sum_{i=0}^{l-1}a_iz_2^{i+1}\Bigr)^pz_2^q, \ \Bigl(z_1z_2^l+\sum_{i=0}^{l-1}a_iz_2^{i+1}\Bigr)^rz_2^s\right)$$
 is conjugated to the polynomial germ
$$F(z_1,z_2)=\Bigl(\l z_1z_2^\s+\sum_{i=p+q}^\s b_iz_2^i, z_2^{r+s}\Bigr),$$
where $\l$ depends only on $a_0$ by lemma 4.\ref{dependancedekappa}.
\end{Prop}
Proof:  Denote by $f_{a_0,\e}$ be the rational mapping of lemma 4.\ref{frationnelle}.\\
For $1\le j\le l-1$,  the homogeneous part of bidegree $(0,p+q+j)$ is:
$$\begin{array}{l}\dps b_{p+q+j}C^{p+q+j}z_2^{p+q+j}+\sum_{j'=1}^{j-1}b_{p+q+j'}C^{p+q+j'}z_2^{p+q+j'}P_{jj'}(a_1,\ldots,a_{j-j'})z_2^{j-j'}\\
\\
\dps-\sum_{(\a,\b)\neq (1,0)\atop {\a(p+q)\atop{+\b(r+s)\atop\le p+q+j}}}A_{\a\b}a_0^{\a p+\b r}z_2^{\a(p+q)+\b(r+s)}\left\la \left\{1+\frac{1}{a_0}\sum_{i=1}^{l-1}a_iz_2^i\right\}^{\a p+\b r}\right\ra_{\bigl(0,p+q+j-\a(p+q)-\b(r+s)\bigr)}\\
\\
\dps =A_{10}z_2^{p+q}\left\la a_0^p\left\{1+\frac{1}{a_0}\sum_{i=1}^{l-1}a_iz_2^i\right\}^{p}\right\ra_{(0,j)}-C^{p+q}z_2^{p+q}\Bigl\la(1+\mu(z))^{p+q}\Bigr\ra_{(0,j)}
\end{array}$$
After cancellation of $z_2^{p+q+j}$ and recalling that $A_{10}a_0^p=C^{p+q}$, we obtain the $j$-th equation
$$\begin{array}{l}\dps b_{p+q+j}C^{p+q+j}+\sum_{j'=1}^{j-1}b_{p+q+j'}C^{p+q+j'}P_{jj'}(a_1,\ldots,a_{j-j'})\\
\\
\dps-\sum_{(\a,\b)\neq (1,0)\atop {\a(p+q)\atop{+\b(r+s)\atop\le p+q+j}}}A_{\a\b}a_0^{\a p+\b r}\frac{\left\la \left\{1+\dps\frac{1}{a_0}\sum_{i=1}^{l-1}a_iz_2^i\right\}^{\a p+\b r}\right\ra_{\bigl(0,p+q+j-\a(p+q)-\b(r+s)\bigr)}}{z_2^{p+q+j-\a(p+q)-\b(r+s)}}\\
\\
\dps =C^{p+q}\left(\frac{\left\la \left\{1+\frac{1}{a_0}\sum_{i=1}^{l-1}a_iz_2^i\right\}^{p}\right\ra_{(0,j)}-\Bigl\la(1+\mu(z))^{p+q}\Bigr\ra_{(0,j)}}{z_2^j}\right)
\end{array}$$
We show  that for $j=1,\ldots,l-1$,
$$b_{p+q+j}=b_{p+q+j}(a_1,\ldots,a_j)= \frac{\d a_j}{Ca_0(r+s)}+R_j(a_1,\ldots,a_{j-1})\quad {\rm with}\ R_j\in \bb C[a_1,\ldots,a_{j-1}].$$
In fact
\begin{itemize}
\item For $j'=1,\ldots,j-1$, $P_{jj'}\in \bb C[a_1,\ldots,a_{j-1}]$,
\item Since $p+q+j-\a(p+q)-\b(r+s)<j$,
$$\frac{\left\la \left\{1+\frac{1}{a_0}\sum_{i=1}^{l-1}a_iz_2^i\right\}^{\a p+\b r}\right\ra_{\bigl(0,p+q+j-\a(p+q)-\b(r+s)\bigr)}}{z_2^{p+q+j-\a(p+q)-\b(r+s)}}\in \bb C[a_1,\ldots,a_{j-1}],$$
\item Clearly $mod\ \frak M_{j-1}:=\bigl(a_1,\ldots,a_{j-1}\bigr)$,
$$\frac{\left\la \left\{1+\frac{1}{a_0}\sum_{i=1}^{l-1}a_iz_2^i\right\}^{p}\right\ra_{(0,j)}}{z_2^j}=\frac{\left\la \left\{1+\frac{1}{a_0}\sum_{i=1}^{j}a_iz_2^i\right\}^{p}\right\ra_{(0,j)}}{z_2^j}=\frac{pa_j}{a_0}$$
\item From the definition of $\mu$, $mod\ \frak M_{j-1}$,
$$\frac{ \Bigl\la(1+\mu(z))^{p+q}\Bigr\ra_{(0,j)}}{z_2^j}=\frac{\left\la \left\{1+\frac{1}{a_0}\sum_{i=1}^{j}a_{i}z_2^{i}\right\}^{\frac{r(p+q)}{r+s}}\cdots\right\ra}{z_2^j}=\frac{r(p+q)a_j}{(r+s)a_0}$$
\end{itemize}
Therefore modulo $\frak M_{j-1}$,
$$b_{p+q+j}C^{j}=\frac{\d \, a_j}{(r+s)a_0}$$
and there exists a polynomial $R_j(a_1,\ldots,a_{j-1})\in\bb C[a_1,\ldots,a_{j-1}]$ without constant term such that
$$b_{p+q+j}=b_{p+q+j}(a_1,\ldots,a_j)=\frac{\d \, a_j}{C^j(r+s)a_0}+R_j(a_1,\ldots,a_{j-1}).\leqno{(T)}$$

\hfill$\Box$\\

\begin{Cor} We choose $a_0\in\bb C^\star$, $\e$ such that $\e^{r+s-1}=1$. Let  $\s=p+q+l-1$ and suppose that $r+s-1$ does not divide $l-d$. Then 
\DIAGV{70}
{\bb C^{l-1}}\n{}\n{\Earv{f_{a_0,\e}}{80}}\n{}\n{\bb C^{l-1}}\nn
{}\n{\Searv{\Pi_L}{50}}\n{}\n{\Swarv{\Pi_L}{50}}\nn
{}\n{}\n{\bb C^{l-1}/\bb Z_{k-1}\quad}
\diag
where  $$f_{a_0,\e}:\bb C^{l-1}\to \bb C^{l-1},\quad (a_1,\ldots,a_{l-1})\mapsto (b_{p+q+1},\ldots,b_{p+q+l-1})$$
is a commutative diagram and $f_{a_0,\e}$ is an isomorphic polynomial mapping.
 \end{Cor}

\subsubsection{Explicit construction of the isomorphic polynomial mapping (there exists global twisted vector fields)}
{\bf\boldmath We suppose that $l-d=K(r+s-1)$ i.e. there are non trivial global twisted vector fields}. We have
$$l+K=d+Kk(S), \quad \s=p+q+l-1=(k(S)-1)(K+1),\quad \frac{\s k(S)}{k(S)-1}=k(S)(K+1).$$
We denote by
$$\left(\sum\right):=\left(\sum_{i=1}^{l-1} a_iz_2^{i}+a_{l+K}\, z_2^{l+K}+z_1z_2^{l-1}\right)$$
the equation $(I)$ is now
$$\left\{\begin{array}{l}
\dps\f_1\left(a_0^p\left\{1+\frac{1}{a_0}\left(\sum\right)\right\}^pz_2^{p+q}, a_0^r\left\{1+\frac{1}{a_0}\left(\sum\right)\right\}^rz_2^{r+s}\right)\\
\\
\hspace{5mm}\dps = \l \f_1(z)\Bigl(C z_2(1+\mu(z)\Bigr)^{\s}+\sum_{i=p+q}^{\s} b_i\Bigl(Cz_2(1+\mu(z))\Bigr)^i
+c\Bigl(Cz_2(1+\mu(z))\Bigr)^{\frac{\s k}{k-1}}\end{array}\right.\leqno{(I)}$$
If $m\ge \frac{\s }{k-1}$, then  $(0,m)\not\in E_\infty$, in fact
$$m(r+s)\ge \frac{\s k}{k-1}=(p+q+l-1)\frac{r+s}{r+s-1},$$
therefore the coefficients $a_{l+K}$ and $c$ doesn't occur in the previous calculations and we obtain by similar arguments a polynomial mapping $f_{a_0,\e}$.

\begin{Lem} \label{casl=d+K(r+s-1)} Suppose $l=d+K(r+s-1)$. Let $\frak M=\frak M_{l-1}=(a_1,\ldots,a_{l-1})$ and $(i,j)$ such that
$$\begin{array}{l}\left\la A_{ij}a_0^{pi+rj}\left\{1+\frac{1}{a_0}\left(\sum_{m=1}^{l-1}a_mz_2^m+a_{l+K}z_2^{l+K}+z_1z_2^{l-1}\right)\right\}^{pi+rj}z_2^{i(p+q)+j(r+s)}\right\ra_{(0,\frac{\s k}{k-1})}\neq 0,\\ \hfill{\rm mod}\quad \frak M 
\end{array}$$
then, $(i,j)=(1,0)$ or $(i,j)=(0,\frac{\s}{k-1})$. More precisely homogeneous part of bidegree $(0,\frac{\s k}{k-1})$ is 
$$c\, C^{\frac{\s k}{k-1}}=A_{10}\, p a_0^{p-1}a_{l+K}+A_{0,\frac{\s}{k-1}}C^\s(1-\l) \quad {\rm mod}\ \frak M.$$
In particular if there are global vector fields, i.e. $\l=1$,
$$c\, C^{\frac{\s k}{k-1}}=A_{10}\, p a_0^{p-1}a_{l+K}\quad {\rm mod}\ \frak M.$$
\end{Lem}
Proof: 1) If $(i,j)\in E_\infty$, then $i(p+q)+j(r+s)<\frac{\s k}{k-1}$ and $i\le 1$ by lemma 4.\ref{populationdecoefficients} 4).
\begin{itemize}
\item Case $i=1$: Since 
$$l+K+(p+q)+j(r+s)=\frac{\s k}{k-1}+jk\ge \frac{\s k}{k-1}$$
we have equality if $j=0$ hence $(i,j)=(1,0)$.
\item  Case $i=0$: then $1\le j< \frac{\s}{k-1}$ and mod $\frak M$,
$$\left\la A_{0j}a_0^{rj}\left\{1+\frac{a_{l+K}z_2^{l+K}}{a_0}\right\}^{rj}z_2^{j(r+s)}\right\ra_{(0,\frac{\s k}{k-1})}\neq 0,\\ \hfill{\rm mod}\quad \frak M $$
 In the left member the possible powers of $z_2$ are of the form  $\a(l+K)+jk$ with $\a\ge 0$ and $j\ge 1$ such that
$$\a(l+K)+jk=\frac{\s k}{k-1}.$$
Since $\a(l+K)+jk=(d+Kk)\a+jk$, we derive that  $\a\ge 1$ is impossible, therefore $\a=0$ and 
$j=\frac{\s}{k-1}$.
\end{itemize}
2) From 1) we deduce that there exists a polynomial $P$ in variables $a_1,\ldots,a_{l-1}$ such that the coefficients of $z_2^{\frac{\s k}{k-1}}$ in $(I)$ give the equality
$$A_{10}\, p a_0^{p-1}a_{l+K}+A_{0,\frac{\s}{k-1}}a_0^{\frac{\s r}{k-1}}=\l A_{0\frac{\s}{k-1}} C^\s + c\, C^{\frac{\s k}{k-1}} +P(a_1,\ldots,a_{l-1}).$$
By equation (\ref{a0C}), 
$$a_0^{\frac{\s r}{k-1}}-\l C^\s=C^\s(1-\l)$$
which gives the result.
\hfill$\Box$\\
\begin{Prop}  If $l-d=K(r+s-1)$ and $\l=1$, there is a bijective triangular polynomial mapping 
$$\begin{array}{cccc}
g_{a_0,\e}:& \bb C^{l-1}\times \bb C&\longrightarrow&  \bb C^{l-1}\times \bb C\\
&a=(a_1,\ldots,a_{l-1},a_{l+K})&\longmapsto&\Bigl(b_{p+q+1}(a),\ldots,b_{p+q+l-1}(a),c(a)\Bigr)
\end{array}$$
such that,
  $$b_{p+q+j}(a)=\frac{\d \, a_j}{C^j(r+s)a_0}+R_j(a_1,\ldots,a_{j-1}), \quad  j=1,\ldots,l-1,$$ 
  $$c(a)=C^{-\frac{\s k}{k-1}}A_{10}\, p a_0^{p-1}a_{l+K}+R(a_1,\ldots,a_{l-1}),$$
$$G(z_1,z_2)=\left(\Bigl(z_1z_2^l+\sum_{i=0}^{l-1}a_iz_2^{i+1}+a_{l+K}z_2^{l+K+1}\Bigr)^pz_2^q, \ \Bigl(z_1z_2^l+\sum_{i=0}^{l-1}a_iz_2^{i+1}+a_{l+K}z_2^{l+K+1}\Bigr)^rz_2^s\right)$$
 is conjugated to the polynomial germ
 $$F(z_1,z_2)=\left(\l z_1z_2^\s+\sum_{k=p+q}^\s b_kz_2^k+ cz_2^{\frac{\s k(S)}{k(S)-1}}, z_2^{r+s}\right).$$
\end{Prop}
Proof: We have a bijective polynomial map 
$$f_{a_0,\e}:\bb C^{l-1}\to \bb C^{l-1}, \quad a\mapsto b=f_{a_0,\e}(a).$$
From lemma 5.\ref{casl=d+K(r+s-1)}, when $a=(a_1,\ldots,a_{l-1})$ is fixed and $a_{l+K}\in \bb C$, the mapping $c:\bb C\to\bb C$, $a_{l+K}\mapsto c=c(a_{l+K})=C^{-\frac{\s k}{k-1}}A_{10}\, p a_0^{p-1}a_{l+K}$ is linear hence bijective.\\
\hfill$\Box$

\begin{Cor} Any surface with GSS with one branch admits a special birational structure.
\end{Cor}

\begin{Cor} The intersection $A:=Aut(\bb C^2,H,0)\cap \Phi$ is the trivial group or a group isomorphic to $(\bb C,+)$. Moreover
\begin{itemize}
\item if $k-1$ does not divide $\frak s=p+q+l-1$, the canonical mapping 
$$g:\cal G/A=\cal G(p,q,r,s,l)/A\to U_{k,\frak s,m_1}/\bb Z_{k-1}$$
  to the Oeljeklaus-Toma coarse moduli space of marked surfaces $(S,C_0)$ with one branch 
  $$U_{k,\frak s,m_1}/\bb Z_{k-1}=\bb C^\star\times\bb C^{l-1}/\bb Z_{k-1}$$
   is isomorphic and  there is a polynomial lifting 
   $$(\l,b):\bb C^\star\times \bb C^{l-1}\to \bb C^\star\times \bb C^{l-1}$$
    which is a covering such that
\DIAGV{50}
{\bb C^\star\times \bb C^{l-1}}\n{}\n{}\n{\Ear{(\l,b)}}\n{}\n{}\n{\bb C^\star\times \bb C^{l-1}}\nn
{\sar}\n{}\n{}\n{}\n{}\n{}\n{\sar}\nn
{\cal G /A}\n{}\n{}\n{\Ear{g}}\n{}\n{}\n{U_{k,\frak s,m_1}/\bb Z_{k-1}}
\diag
is commutative,
\item  if $k-1$ divides $\frak s=p+q+l-1$, we have similar results for 
$$U_{k,\frak s,m_1}^{\l\neq 0,c=0}/\bb Z_{k-1}\quad {\rm and}\quad U_{k,\frak s,m_1}^{\l=1}/\bb Z_{k-1}.$$
\end{itemize}
\end{Cor}
\begin{Cor} Let $\cal S_{J,\s}\to B_J$ be a large family with $\s=Id$. Let $T_{J,\s}$ the hypersurface where cocycles $[\t^i]$ and $[\mu^i]$ are not independent. Then for each stratum $B_{J,M}$, the trace  $T_{J,\s}\cap B_{J,M}$ on $B_{J,M}$  is equal to the inverse image of the ramification set by the lift of the canonical mapping i.e.
\begin{itemize}
\item If $k-1$ does not divide $\frak s$,
$$T_{J,\s}\cap B_{J,M}=(\l,b)^{-1}(T_{k,\frak s,m_1}),$$
\item If $k-1$ divides $\frak s$
$$T_{J,\s}\cap B_{J,M}=(\l,b)^{-1}(T_{k,\frak s,m_1}^{\l\neq 1,c=0}).$$
\end{itemize}
In particular in $B_J$ there is no curve  over which the surfaces are isomorphic.
\end{Cor}

\subsection{Special birational structures on Kato surfaces}
Let $S$ be a  Kato surface with $b_2(S)=n$, let $D_0,\ldots,D_{n-1}$ be its rational curves, $p:\tilde S\to S$ its universal cover, $C_0$ a lift of $D_0$, $C_i$, $i\in\bb Z$, the rational curves in $\tilde S$ in the canonical order and $\frak U=(U_i)_{0\le i\le n-1}$ an Enoki covering of $S$ such that $U_0$ contains $C_0$ with a  deleted disc (see the construction of these surfaces). 
If $S$ is associated to a germ $F=\Pi\s$ where $\s$ is birational, then $S$ is endowed with a birational structures as well as $\tilde S$ and $p$  is a $\bigl(Bir(\bb P^2(\bb C)),\bb P^2(\bb C)\bigr)$-morphism. 
\begin{Def} A birational structure on a surface with GSS $S$ will be called {\bf special} if there is  a contracting germ $F=\Pi\s$ with $F$ or equivalently $\s$ birational and $S=S(F)$.
\end{Def}
This definition is independent of the numbering. 
Each open set $U_i$ is covered by two charts $U_i'$ and $U_i''$ with local coordinates $\f'_i=(u_i',v_i'):U_i'\to\bb C^2$ and $ \f_i''=(u_i'',v_i'') :U_i''\to \bb C^2$ respectively. We denote by $\f_i=(u_i,v_i)$  the  local coordinates whose domain contains the blown-up point $O_i\in C_i$. Here $\f_i=Id$. Then, with the identification 
$$i:\bb C^2\simeq \{[z_0:z_1:z_2]\in\bb P^2(\bb C)\mid z_2=1\}\subset \bb P^2(\bb C),$$ 
$$b_{i,i+1}(u_{i+1},v_{i+1})=\f_i\circ\f_{i+1}^{-1}(u_{i+1},v_{i+1})=\Pi_{i+1}(u_{i+1},v_{i+1}), \ i=0,\ldots,n-2, \quad \b_{n-1,0}(u_{0},v_{0})=\s\circ\Pi_0(u_{0},v_{0}).$$
If $S$ contains a GSS  but is not minimal the order on the curves is no more total.

\begin{Lem}\label{mermap} Let $S=S(\Pi,\s)$ be a  surface containing a GSS (not necessarily minimal) such that $n=b_2(S)$. If $\s$ is birational,  there exist for any $j\in\bb Z$
\begin{itemize}
\item a meromorphic developing map $\widehat{Dev}_j:\hat S_{C_j}\to \bb P^2$, locally biholomorphic outside the rational curves, such that  $\widehat{Dev}_j(C_j)$ is the rational curve $\{z_1=0\}\subset \bb P^2(\bb C)$ and $O_j:=\widehat{Dev}_j(\hat O_{C_j})=[a_j:b_j:1]$,
\item a birational mappings $G_{j}: \bb P^2(\bb C) \to \bb P^2(\bb C)$, holomorphic in a neighbourhood of $\hat O_{C_j}$.
\end{itemize}
such that the following diagrams are commutative:
\DIAGV{70}
{\tilde S}\n{}\n{\Earv{\tilde g}{90}}\n{}\n{\tilde S}\n{}\n{}\n{\tilde S}\n{\Ear{\tilde g}}\n{\tilde S}\nn
{\Sar{p_{C_j}}}\n{}\n{}\n{}\n{\saR{p_{C_{j+n}}}}\n{}\n{}\n{\Sar{p_{C_j}}}\n{}\n{\saR{p_{C_j}}}\nn
{\hat S_{C_j}}\n{}\n{\Earv{\s_{j}^{j+n}}{90}}\n{}\n{\hat S_{ C_{j+n}}}\n{}\n{}\n{\hat S_{C_j}}\n{\Ear{F_{C_j}}}\n{\hat S_{C_j}}\nn
{}\n{\seaRv{\widehat{Dev}_j}{40}}\n{}\n{\swaRv{\widehat{Dev}_{j+n}}{40}}\n{}\n{}\n{}\n{\Sar{\widehat{Dev}_j}}\n{}\n{\saR{\widehat{Dev}_j}}\nn
{}\n{}\n{\bb P^2(\bb C)}\n{}\n{}\n{}\n{}\n{\bb P^2(\bb C)}\n{\Earv{G_{j}}{30}}\n{\bb P^2(\bb C)}
\diag
Moreover if $S$ is minimal (i.e. if $S$ is a Kato surface) then $O_j$ is a fixed point of $G_j$.
\end{Lem}
Proof: To simplify the notations we may suppose that $j=0$. Recall that by construction $W_0\subset \hat S_{C_0}$ and then we apply lemma 2.\ref{Extbir}. Since all birational transition functions are isomorphic outside the curves, the developing map $\widehat{Dev}_0$ is a local biholomorphism outside the rational curves. 
We denote a blowup $\Pi_{i+1}:U_{i+1}\to B_i\subset W_i$ by 
$$\Pi_{i+1}:(u'_{i+1},v'_{i+1})\mapsto (u'_{i+1}v'_{i+1}+a_i,v'_{i+1})=(u_i,v_i),$$
$$\Pi_{i+1}:(u''_{i+1},v''_{i+1})\mapsto  (v''_{i+1}+a_i,u''_{i+1}v''_{i+1})=(u_i,v_i)$$
 with inverse 
 $$\Pi_{i+1}^{-1}:(u_i,v_i,)\mapsto (u'_{i+1},v'_{i+1})=\left(\frac{u_i-a_i}{v_i}, v_i\right),$$
$$\Pi_{i+1}^{-1}:(u_i,v_i,)\mapsto (u''_{i+1},v''_{i+1})=\left(\frac{v_i}{u_i-a_i},u_i-a_i\right)$$
On $U_0$, $\widehat{Dev}_0$ is defined in the following way: If $\hat O_{C_0}$ is in the chart $(u'_0,v'_0)$,
$$\widehat{Dev}_0(u'_0,v'_0)=[u'_0:v'_0:1], \quad \widehat{Dev}_0({u_0''},{v_0''})=[\frac{1}{u''_0}:u''_0v''_0:1]=[1:{u''_0}^2v''_0:u''_0].$$
 The construction is similar if $\hat O_{C_0}$ is in the chart $(u''_0,v''_0)$.\\
On   $\bigcup_{k<0}U_k$, we have
 $$\begin{array}{lccl}
{\widehat{Dev}_0}:&U_k &\to &\bb P^2(\bb C)\\ &(u_k,v_k)&\mapsto& \widehat{Dev}_0(u_k,v_k)=i\circ b_{0,-1}\circ\cdots\circ b_{k+1,k}(u_k,v_k)\end{array}$$
i.e. for $k>-n-1$,
$$\widehat{Dev}_0(u_k,v_k)=i\circ \Pi_0^{-1}\circ\s_0^{-1}\circ\Pi_{-1}^{-1}\circ\cdots\circ\Pi_{k+1}^{-1}(u_k,v_k)
$$ and 
$\s_0\Pi_0:U_0\to U_{-1}$ is the composition of $\Pi_0:U_{0}\to B$ and of $\s_0:B\to U_{-1}$ which  is birational induced by $\s$. The image $\s_0\Pi_0(U_0)$ is a ball in $W_{-1}$; $\s_0\Pi_0$ being birational $\Pi_0^{-1}\s_0^{-1}$ extends to $U_{-1}$. \\
The points $(a_k,0)\in W_k$ are indeterminacy points of $\Pi^{-1}_{k+1}$ however do not belong to $U_k$. Therefore $\Pi_{k+1}^{-1}(U_k)$ has an empty intersection with $C_{k+1}$ and $\widehat{Dev}_0$ is holomorphic.\\
 The upper parts of the diagrams are commutative by \cite{D1}, p30; to see the commutativity of the lower parts it is sufficient to check it on the chart $(u_0,v_0)$.
 \hfill $\Box$
 The following theorem shows that we recover a GSS in $S$ thanks to a small sphere centered at $\widehat{Dev}_j(\hat O_j)$.
\begin{Th} Let $S=S(\Pi,\s)$ be a  Kato surface such that $n=b_2(S)$. If $\s$ is birational,  there exist for any $j\in\bb Z$ a meromorphic developing map $\widetilde{Dev}_j:\tilde S\to \bb P^2$, locally biholomorphic outside the rational curves, such that  $\widetilde{Dev}_j(C_j)$ is the rational curve $\{z_1=0\}\subset \bb P^2(\bb C)$ and $O_j:=\widetilde{Dev}_j(p_{C_j}^{-1}(\hat O_{C_j}))=[a_j:b_j:1]$, i.e. $\widetilde{Dev}_j$ blows down an infinite number of curves. Moreover for any small ball $B_j$ centered at $O_j$, $p(\widetilde{Dev}_j^{-1}(\part B_j)$ is a GSS in $S$.
\end{Th}
Proof: $\widetilde{Dev}_j=\widehat{Dev}_j \circ p_{C_j}$ has the expected properties.\hfill$\Box$

\end{document}

%% file: diagram.tex
% DIAGRAM MACROS (version 2.2)

% If a document has been previously realized using version 1.0 or 1.1
% of the DIAGRAM macros and if you have to recompile it using now
% version 2.2, it is compulsory to modify in this document each
% occurence of \Diag, \Diagv, \diaG, \diagsV and \move
% according to the new rules.

% Typeset the document DIAGRAM.READ.ME to get user's information
% about the following macros.

% MACROS FOR DRAWING IN-TEXT PICTURES

% \tlowername{P}{f} puts the name f under the picture P
\newcommand{\tlowername}[2]%
{$\stackrel{\makebox[1pt]{#1}}%
{\begin{picture}(0,0)%
\put(0,0){\makebox(0,6)[t]{\makebox[1pt]{$#2$}}}%
\end{picture}}$}%

% \tcase{P} draws the picture P with length 20pt
\newcommand{\tcase}[1]{\makebox[23pt]%
{\raisebox{2.5pt}{#1{20}}}}%

% \Tcase{P}{f} draws the picture P with upper name f
% and length 20pt.
\newcommand{\Tcase}[2]{\makebox[23pt]%
{\raisebox{2.5pt}{$\stackrel{#2}{#1{20}}$}}}%

% \tbicase{P} draws the bi-picture P with length 20pt
\newcommand{\tbicase}[1]{\makebox[23pt]%
{\raisebox{1pt}{#1{20}}}}%

% \Tbicase{P}{f}{g} draws the bi-picture P with names f, g
% and length 20pt.
\newcommand{\Tbicase}[3]{\makebox[23pt]{\raisebox{-7pt}%
{$\stackrel{#2}{\mbox{\tlowername{#1{20}}{\scriptstyle{#3}}}}$}}}%

% IN-TEXT ARROWS

% \AR{n} draws an arrow of length n units
\newcommand{\AR}[1]%
{\begin{picture}(#1,0)%
\put(0,0){\vector(1,0){#1}}%
\end{picture}}%

% \DOTAR{n} draws a dotted  arrow of length n units
\newcommand{\DOTAR}[1]%
{\NUMBEROFDOTS=#1%
\divide\NUMBEROFDOTS by 3%
\begin{picture}(#1,0)%
\multiput(0,0)(3,0){\NUMBEROFDOTS}{\circle*{1}}%
\put(#1,0){\vector(1,0){0}}%
\end{picture}}%

% \MONO{n} draws a monomorphism of length n units
\newcommand{\MONO}[1]%
{\begin{picture}(#1,0)%
\put(0,0){\vector(1,0){#1}}%
\put(2,-2){\line(0,1){4}}%
\end{picture}}%

% \EPI{n} draws an epimorphism of length n units
\newcommand{\EPI}[1]%
{\begin{picture}(#1,0)(-#1,0)%
\put(-#1,0){\vector(1,0){#1}}%
\put(-6,-2){\line(0,1){4}}%
\end{picture}}%

% \BIMO{n} draws a bimorphism of length n units
\newcommand{\BIMO}[1]%
{\begin{picture}(#1,0)(-#1,0)%
\put(-#1,0){\vector(1,0){#1}}%
\put(-6,-2){\line(0,1){4}}%
\put(-#1,-2){\hspace{2pt}\line(0,1){4}}%
\end{picture}}%

% \BIAR{n} draws a pair of arrows of length n units
\newcommand{\BIAR}[1]%
{\begin{picture}(#1,4)%
\put(0,0){\vector(1,0){#1}}%
\put(0,4){\vector(1,0){#1}}%
\end{picture}}%

% \EQL{n} draws an equality of length n units
\newcommand{\EQL}[1]%
{\begin{picture}(#1,0)%
\put(0,1){\line(1,0){#1}}%
\put(0,-1){\line(1,0){#1}}%
\end{picture}}%

% \ADJAR{n} draws a pair of adjoint arrows of length n units
\newcommand{\ADJAR}[1]%
{\begin{picture}(#1,4)%
\put(0,0){\vector(1,0){#1}}%
\put(#1,4){\vector(-1,0){#1}}%
\end{picture}}%

% All the following cammands produce arrows with length 20pt
% between 1.5pt spaces.

% arrow
\newcommand{\ar}{\tcase{\AR}}%

% arrow with upper name [1]
\newcommand{\Ar}[1]{\Tcase{\AR}{#1}}%

% dotted arrow
\newcommand{\dotar}{\tcase{\DOTAR}}%

% dotted arrow with upper name [1]
\newcommand{\Dotar}[1]{\Tcase{\DOTAR}{#1}}%

% monomorphism
\newcommand{\mono}{\tcase{\MONO}}%

% monomorphism with upper name [1]
\newcommand{\Mono}[1]{\Tcase{\MONO}{#1}}%

% epimorphism
\newcommand{\epi}{\tcase{\EPI}}%

% epimorphism with upper name [1]
\newcommand{\Epi}[1]{\Tcase{\EPI}{#1}}%

% bimorphism
\newcommand{\bimo}{\tcase{\BIMO}}%

% bimorphism with upper name [1]
\newcommand{\Bimo}[1]{\Tcase{\BIMO}{#1}}%

% isomorphism
\newcommand{\iso}{\Tcase{\AR}{\cong}}%

% isomorphism with upper name [1]
\newcommand{\Iso}[1]{\Tcase{\AR}{\cong{#1}}}%

% pair of arrows
\newcommand{\biar}{\tbicase{\BIAR}}%

% pair of arrows with names [1],[2]
\newcommand{\Biar}[2]{\Tbicase{\BIAR}{#1}{#2}}%

% equality
\newcommand{\eql}{\tcase{\EQL}}%

% pair of adjoint arrows
\newcommand{\adjar}{\tbicase{\ADJAR}}%

% pair of adjoint arrows with names [1],[2]
\newcommand{\Adjar}[2]{\Tbicase{\ADJAR}{#1}{#2}}%

% IN-TEXT BACK ARROWS

% \BKAR{n} draws a pointing back arrow of length n units
\newcommand{\BKAR}[1]%
{\begin{picture}(#1,0)%
\put(#1,0){\vector(-1,0){#1}}%
\end{picture}}%

% \BKDOTAR{n} draws a backward dotted  arrow
% of length n units
\newcommand{\BKDOTAR}[1]%
{\NUMBEROFDOTS=#1%
\divide\NUMBEROFDOTS by 3%
\begin{picture}(#1,0)%
\multiput(#1,0)(-3,0){\NUMBEROFDOTS}{\circle*{1}}%
\put(0,0){\vector(-1,0){0}}%
\end{picture}}%

% \BKMONO{n} draws a pointing back monomorphism of length n units
\newcommand{\BKMONO}[1]%
{\begin{picture}(#1,0)(-#1,0)%
\put(0,0){\vector(-1,0){#1}}%
\put(-2,-2){\line(0,1){4}}%
\end{picture}}%

% \BKEPI{n} draws a pointing back epimorphism of length n units
\newcommand{\BKEPI}[1]%
{\begin{picture}(#1,0)%
\put(#1,0){\vector(-1,0){#1}}%
\put(6,-2){\line(0,1){4}}%
\end{picture}}%

% \BKBIMO{n} draws a pointing back bimorphism of length n units
\newcommand{\BKBIMO}[1]%
{\begin{picture}(#1,0)%
\put(#1,0){\vector(-1,0){#1}}%
\put(6,-2){\line(0,1){4}}%
\put(#1,-2){\hspace{-2pt}\line(0,1){4}}%
\end{picture}}%

% \BKBIAR{n} draws a pair of pointing back arrows of length n units
\newcommand{\BKBIAR}[1]%
{\begin{picture}(#1,4)%
\put(#1,0){\vector(-1,0){#1}}%
\put(#1,4){\vector(-1,0){#1}}%
\end{picture}}%

% \BKADJAR{n} draws a pair of adjoint arrows of length n units
\newcommand{\BKADJAR}[1]%
{\begin{picture}(#1,4)%
\put(0,4){\vector(1,0){#1}}%
\put(#1,0){\vector(-1,0){#1}}%
\end{picture}}%

% All the following commands produce back arrows with length 20pt
% between 1.5pt spaces.

% back arrow
\newcommand{\bkar}{\tcase{\BKAR}}%

% back arrow with upper name [1]
\newcommand{\Bkar}[1]{\Tcase{\BKAR}{#1}}%

% backward dotted arrow
\newcommand{\bkdotar}{\tcase{\BKDOTAR}}%

% backward dotted arrow with upper name [1]
\newcommand{\Bkdotar}[1]{\Tcase{\BKDOTAR}{#1}}%

% back monomorphism
\newcommand{\bkmono}{\tcase{\BKMONO}}%

% back monomorphism with upper name [1]
\newcommand{\Bkmono}[1]{\Tcase{\BKMONO}{#1}}%

% back epimorphism
\newcommand{\bkepi}{\tcase{\BKEPI}}%

% back epimorphism with upper name [1]
\newcommand{\Bkepi}[1]{\Tcase{\BKEPI}{#1}}%

% back bimorphism
\newcommand{\bkbimo}{\tcase{\BKBIMO}}%

% back bimorphism with upper name [1]
\newcommand{\Bkbimo}[1]{\Tcase{\BKBIMO}{\hspace{9pt}#1}}%

% back isomorphism
\newcommand{\bkiso}{\Tcase{\BKAR}{\cong}}%

% back isomorphism with upper name [1]
\newcommand{\Bkiso}[1]{\Tcase{\BKAR}{\cong{#1}}}%

% pair of back arrows
\newcommand{\bkbiar}{\tbicase{\BKBIAR}}%

% pair of back arrows with names [1],[2]
\newcommand{\Bkbiar}[2]{\Tbicase{\BKBIAR}{#1}{#2}}%

% back equality
\newcommand{\bkeql}{\tcase{\EQL}}%

% back pair of adjoint arrows
\newcommand{\bkadjar}{\tbicase{\BKADJAR}}%

% back pair of adjoint arrows with names [1],[2]
\newcommand{\Bkadjar}[2]{\Tbicase{\BKADJAR}{#1}{#2}}%

% MACROS FOR DRAWING HORIZONTAL PICTURES

% \lowername{P}{f} puts the name f under the picture P
\newcommand{\lowername}[2]%
{$\stackrel{\makebox[1pt]{#1}}%
{\begin{picture}(0,0)%
\truex{600}%
\put(0,0){\makebox(0,\value{x})[t]{\makebox[1pt]{$#2$}}}%
\end{picture}}$}%

% \hcase{P}{n} draws the picture P with length n units
\newcommand{\hcase}[2]%
{\makebox[0pt]%
{\raisebox{-1pt}[0pt][0pt]{#1{#2}}}}%

% \Hcase{P}{f}{n} draws the picture P with upper name f
% and length n units.
\newcommand{\Hcase}[3]%
{\makebox[0pt]
{\raisebox{-1pt}[0pt][0pt]%
{$\stackrel{\makebox[0pt]{$\textstyle{#2}$}}{#1{#3}}$}}}%

% \hcasE{P}{f}{n} draws the picture P with lower name f
% and length n units.
\newcommand{\hcasE}[3]%
{\makebox[0pt]%
{\raisebox{-9pt}[0pt][0pt]%
{\lowername{#1{#3}}{#2}}}}%

% \hbicase{P}{n} draws the bi-picture P with length n units.
\newcommand{\hbicase}[2]%
{\makebox[0pt]%
{\raisebox{-2.5pt}[0pt][0pt]{#1{#2}}}}%

% \Hbicase{P}{f}{g}{n} draws the bi-picture P with names f, g
% and length n units.
\newcommand{\Hbicase}[4]%
{\makebox[0pt]
{\raisebox{-10.5pt}[0pt][0pt]%
{$\stackrel{\makebox[0pt]{$\textstyle{#2}$}}%
{\mbox{\lowername{#1{#4}}{#3}}}$}}}%

% EAST ARROWS

% \EAR{n} draws an east arrow of length n units
\newcommand{\EAR}[1]%
{\begin{picture}(#1,0)%
\put(0,0){\vector(1,0){#1}}%
\end{picture}}%

% \EDOTAR draws a dotted arrow of length n units
\newcommand{\EDOTAR}[1]%
{\truex{100}\truey{300}%
\NUMBEROFDOTS=#1%
\divide\NUMBEROFDOTS by \value{y}%
\begin{picture}(#1,0)%
\multiput(0,0)(\value{y},0){\NUMBEROFDOTS}%
{\circle*{\value{x}}}%
\put(#1,0){\vector(1,0){0}}%
\end{picture}}%

% \EMONO{n} draws an east monomorphism of length n units
\newcommand{\EMONO}[1]%
{\begin{picture}(#1,0)%
\put(0,0){\vector(1,0){#1}}%
\truex{300}\truey{600}%
\put(\value{x},-\value{x}){\line(0,1){\value{y}}}%
\end{picture}}%

% \EEPI{n} draws an east epimorphism of length n units
\newcommand{\EEPI}[1]%
{\begin{picture}(#1,0)(-#1,0)%
\put(-#1,0){\vector(1,0){#1}}%
\truex{300}\truey{600}\truez{800}%
\put(-\value{z},-\value{x}){\line(0,1){\value{y}}}%
\end{picture}}%

% \EBIMO{n} draws an east bimorphism of length n units
\newcommand{\EBIMO}[1]%
{\begin{picture}(#1,0)(-#1,0)%
\put(-#1,0){\vector(1,0){#1}}%
\truex{300}\truey{600}\truez{800}%
\put(-\value{z},-\value{x}){\line(0,1){\value{y}}}%
\put(-#1,-\value{x}){\hspace{3pt}\line(0,1){\value{y}}}%
\end{picture}}%

% \EBIAR{n} draws an east pair of arrows of length n units
\newcommand{\EBIAR}[1]%
{\truex{400}%
\begin{picture}(#1,\value{x})%
\put(0,0){\vector(1,0){#1}}%
\put(0,\value{x}){\vector(1,0){#1}}%
\end{picture}}%

% \EEQL{n} draws an east equality of length n units
\newcommand{\EEQL}[1]%
{\begin{picture}(#1,0)%
\truex{200}%
\put(0,\value{x}){\line(1,0){#1}}%
\put(0,0){\line(1,0){#1}}%
\end{picture}}%

% \EADJAR{n} draws an east pair of adjoint arrows of length n
% units
\newcommand{\EADJAR}[1]%
{\truex{400}%
\begin{picture}(#1,\value{x})%
\put(0,0){\vector(1,0){#1}}%
\put(#1,\value{x}){\vector(-1,0){#1}}%
\end{picture}}%

% All the following commands produce east arrows
% raised at the adequate position and of formal dimensions (0,0)

% east arrow of length [1]x100 units
\newcommand{\earv}[1]{\hcase{\EAR}{#100}}%

% east arrow of adjusted length
\newcommand{\ear}%
{\hspace{\SOURCE\unitlength}%
\hcase{\EAR}{\ARROWLENGTH}}%

% east arrow with upper name [1] and length [2]x100 units
\newcommand{\Earv}[2]{\Hcase{\EAR}{#1}{#200}}%

% east arrow with upper name [1] and adjusted length
\newcommand{\Ear}[1]%
{\hspace{\SOURCE\unitlength}%
\Hcase{\EAR}{#1}{\ARROWLENGTH}}%

% east arrow with lower name [1] and length [2]x100 units
\newcommand{\eaRv}[2]{\hcasE{\EAR}{#1}{#200}}%

% east arrow with lower name [1] and adjusted length
\newcommand{\eaR}[1]%
{\hspace{\SOURCE\unitlength}%
\hcasE{\EAR}{#1}{\ARROWLENGTH}}%

% east dotted arrow of length [1]x100 units
\newcommand{\edotarv}[1]{\hcase{\EDOTAR}{#100}}%

% east dotted arrow of adjusted length
\newcommand{\edotar}%
{\hspace{\SOURCE\unitlength}%
\hcase{\EDOTAR}{\ARROWLENGTH}}%

% east dotted arrow with upper name [1] and length [2]x100 units
\newcommand{\Edotarv}[2]{\Hcase{\EDOTAR}{#1}{#200}}%

% east dotted arrow with upper name [1] andadjusted length
\newcommand{\Edotar}[1]%
{\hspace{\SOURCE\unitlength}%
\Hcase{\EDOTAR}{#1}{\ARROWLENGTH}}%

% east dotted arrow with lower name [1] and length [2]x100 units
\newcommand{\edotaRv}[2]{\hcasE{\EDOTAR}{#1}{#200}}%

% east dotted arrow with lower name [1] and adjusted length
\newcommand{\edotaR}[1]%
{\hspace{\SOURCE\unitlength}%
\hcasE{\EDOTAR}{#1}{\ARROWLENGTH}}%

% east monomorphism with length [1]x100 units
\newcommand{\emonov}[1]{\hcase{\EMONO}{#100}}%

% east monomorphism with adjusted length
\newcommand{\emono}%
{\hspace{\SOURCE\unitlength}%
\hcase{\EMONO}{\ARROWLENGTH}}%

% east monomorphism with upper name [1] and length [2]x100 units
\newcommand{\Emonov}[2]{\Hcase{\EMONO}{#1}{#200}}%

% east monomorphism with upper name [1] and adjusted length
\newcommand{\Emono}[1]%
{\hspace{\SOURCE\unitlength}%
\Hcase{\EMONO}{#1}{\ARROWLENGTH}}%

% east monomorphism with lower name [1] and length [2]x100 units
\newcommand{\emonOv}[2]{\hcasE{\EMONO}{#1}{#200}}%

% east monomorphism with lower name [1] and adjusted length
\newcommand{\emonO}[1]%
{\hspace{\SOURCE\unitlength}%
\hcasE{\EMONO}{#1}{\ARROWLENGTH}}%

% east epimorphism with length [1]x100 units
\newcommand{\eepiv}[1]{\hcase{\EEPI}{#100}}%

% east epimorphism with adjusted length
\newcommand{\eepi}%
{\hspace{\SOURCE\unitlength}%
\hcase{\EEPI}{\ARROWLENGTH}}%

% east epimorphism with upper name [1] and length [2]x100 units
\newcommand{\Eepiv}[2]{\Hcase{\EEPI}{#1}{#200}}%

% east epimorphism with upper name [1] and adjusted length
\newcommand{\Eepi}[1]%
{\hspace{\SOURCE\unitlength}%
\Hcase{\EEPI}{#1}{\ARROWLENGTH}}%

% east epimorphism with lower name [1] and length [2]x100 units
\newcommand{\eepIv}[2]{\hcasE{\EEPI}{#1}{#200}}%

% east epimorphism with lower name [1] and adjusted length
\newcommand{\eepI}[1]%
{\hspace{\SOURCE\unitlength}%
\hcasE{\EEPI}{#1}{\ARROWLENGTH}}%

% east bimorphism with length [1]x100 units
\newcommand{\ebimov}[1]{\hcase{\EBIMO}{#100}}%

% east bimorphism with adjusted length
\newcommand{\ebimo}%
{\hspace{\SOURCE\unitlength}%
\hcase{\EBIMO}{\ARROWLENGTH}}%

% east bimorphism with upper name [1] and length [2]x100 units
\newcommand{\Ebimov}[2]{\Hcase{\EBIMO}{#1}{#200}}%

% east bimorphism with upper name [1] and adjusted length
\newcommand{\Ebimo}[1]%
{\hspace{\SOURCE\unitlength}%
\Hcase{\EBIMO}{#1}{\ARROWLENGTH}}%

% east bimorphism with lower name [1] and length [2]x100 units
\newcommand{\ebimOv}[2]{\hcasE{\EBIMO}{#1}{#200}}%

% east bimorphism with lower name [1] and adjusted length
\newcommand{\ebimO}[1]%
{\hspace{\SOURCE\unitlength}%
\hcasE{\EBIMO}{#1}{\ARROWLENGTH}}%

% east isomorphism with length [1]x100 units
\newcommand{\eisov}[1]{\Hcase{\EAR}{\cong}{#100}}%

% east isomorphism with adjusted length
\newcommand{\eiso}%
{\hspace{\SOURCE\unitlength}%
\Hcase{\EAR}{\cong}{\ARROWLENGTH}}%

% east isomorphism with upper name [1] and length [2]x100 units
\newcommand{\Eisov}[2]{\Hcase{\EAR}{\cong#1}{#200}}%

% east isomorphism with upper name [1] and adjusted length
\newcommand{\Eiso}[1]%
{\hspace{\SOURCE\unitlength}%
\Hcase{\EAR}{\cong#1}{\ARROWLENGTH}}%

% east isomorphism with lower name [1] and length [2]x100 units
\newcommand{\eisOv}[2]{\hcasE{\EAR}{\cong#1}{#200}}%

% east isomorphism with lower name [1] and adjusted length
\newcommand{\eisO}[1]%
{\hspace{\SOURCE\unitlength}%
\hcasE{\EAR}{\cong#1}{\ARROWLENGTH}}%

% pair of east arrows of length [1]x100 units
\newcommand{\ebiarv}[1]{\hbicase{\EBIAR}{#100}}%

% pair of east arrows of adjusted length
\newcommand{\ebiar}%
{\hspace{\SOURCE\unitlength}%
\hbicase{\EBIAR}{\ARROWLENGTH}}%

% pair of east arrows with names [1] [2] and length [3]x100 units
\newcommand{\Ebiarv}[3]{\Hbicase{\EBIAR}{#1}{#2}{#300}}%

% pair of east arrows with names [1] [2] and adjusted length
\newcommand{\Ebiar}[2]%
{\hspace{\SOURCE\unitlength}%
\Hbicase{\EBIAR}{#1}{#2}{\ARROWLENGTH}}%

% east equality of length [1]x100 units
\newcommand{\eeqlv}[1]{\hcase{\EEQL}{#100}}%

% east equality of adjusted length
\newcommand{\eeql}%
{\hspace{\SOURCE\unitlength}%
\hbicase{\EEQL}{\ARROWLENGTH}}%

% east pair of adjoint arrows of length [1]x100 units
\newcommand{\eadjarv}[1]{\hbicase{\EADJAR}{#100}}%

% east pair of adjoint arrows of adjusted length
\newcommand{\eadjar}%
{\hspace{\SOURCE\unitlength}%
\hbicase{\EADJAR}{\ARROWLENGTH}}%

% east pair of adjoint arrows with names [1][2] and length [3]x100units
\newcommand{\Eadjarv}[3]{\Hbicase{\EADJAR}{#1}{#2}{#300}}%

% east pair of adjoint arrows with names [1]J[2] and adjusted length
\newcommand{\Eadjar}[2]%
{\hspace{\SOURCE\unitlength}%
\Hbicase{\EADJAR}{#1}{#2}{\ARROWLENGTH}}%

% WEST ARROWS

% \WAR{n} draws a pointing back arrow of length n units
\newcommand{\WAR}[1]%
{\begin{picture}(#1,0)%
\put(#1,0){\vector(-1,0){#1}}%
\end{picture}}%

% \WDOTAR draws a dotted arrow of length n units
\newcommand{\WDOTAR}[1]%
{\truex{100}\truey{300}%
\NUMBEROFDOTS=#1%
\divide\NUMBEROFDOTS by \value{y}%
\begin{picture}(#1,0)%
\multiput(#1,0)(-\value{y},0){\NUMBEROFDOTS}%
{\circle*{\value{x}}}%
\put(0,0){\vector(-1,0){0}}%
\end{picture}}%

% \WMONO{n} draws a pointing back monomorphism of length n units
\newcommand{\WMONO}[1]%
{\begin{picture}(#1,0)(-#1,0)%
\put(0,0){\vector(-1,0){#1}}%
\truex{300}\truey{600}%
\put(-\value{x},-\value{x}){\line(0,1){\value{y}}}%
\end{picture}}%

% \WEPI{n} draws a pointing back epimorphism of length n units
\newcommand{\WEPI}[1]%
{\begin{picture}(#1,0)%
\put(#1,0){\vector(-1,0){#1}}%
\truex{300}\truey{600}\truez{800}%
\put(\value{z},-\value{x}){\line(0,1){\value{y}}}%
\end{picture}}%

% \WBIMO{n} draws a pointing back bimorphism of length n units
\newcommand{\WBIMO}[1]%
{\begin{picture}(#1,0)%
\put(#1,0){\vector(-1,0){#1}}%
\truex{300}\truey{600}\truez{800}%
\put(\value{z},-\value{x}){\line(0,1){\value{y}}}%
\put(#1,-\value{x}){\hspace{-3pt}\line(0,1){\value{y}}}%
\end{picture}}%

% \WBIAR{n} draws a pair of pointing back arrows of length n units
\newcommand{\WBIAR}[1]%
{\truex{400}%
\begin{picture}(#1,\value{x})%
\put(#1,0){\vector(-1,0){#1}}%
\put(#1,\value{x}){\vector(-1,0){#1}}%
\end{picture}}%

% \WADJAR{n} draws a pair of adjoint arrows of length n units
\newcommand{\WADJAR}[1]%
{\truex{400}%
\begin{picture}(#1,\value{x})%
\put(0,\value{x}){\vector(1,0){#1}}%
\put(#1,0){\vector(-1,0){#1}}%
\end{picture}}%

% All the following commands produce west arrows
% raised at the adequate position and of formal dimensions (0,0)

%  west arrow of length [1]x100 units
\newcommand{\warv}[1]{\hcase{\WAR}{#100}}%

%  west arrow of adjusted length
\newcommand{\war}%
{\hspace{\SOURCE\unitlength}%
\hcase{\WAR}{\ARROWLENGTH}}%

% west arrow with upper name [1] and length [2]x100 units
\newcommand{\Warv}[2]{\Hcase{\WAR}{#1}{#200}}%

% west arrow with upper name [1] and adjusted length
\newcommand{\War}[1]%
{\hspace{\SOURCE\unitlength}%
\Hcase{\WAR}{#1}{\ARROWLENGTH}}%

% west arrow with lower name [1] and length [2]x100 units
\newcommand{\waRv}[2]{\hcasE{\WAR}{#1}{#200}}%

% west arrow with lower name [1] and adjusted length
\newcommand{\waR}[1]%
{\hspace{\SOURCE\unitlength}%
\hcasE{\WAR}{#1}{\ARROWLENGTH}}%

% west dotted arrow of length [1]x100 units
\newcommand{\wdotarv}[1]{\hcase{\WDOTAR}{#100}}%

% west dotted arrow of adjusted length
\newcommand{\wdotar}%
{\hspace{\SOURCE\unitlength}%
\hcase{\WDOTAR}{\ARROWLENGTH}}%

% west dotted arrow with upper name [1] and length [2]x100 units
\newcommand{\Wdotarv}[2]{\Hcase{\WDOTAR}{#1}{#200}}%

% west dotted arrow with upper name [1] andadjusted length
\newcommand{\Wdotar}[1]%
{\hspace{\SOURCE\unitlength}%
\Hcase{\WDOTAR}{#1}{\ARROWLENGTH}}%

% west dotted arrow with lower name [1] and length [2]x100 units
\newcommand{\wdotaRv}[2]{\hcasE{\WDOTAR}{#1}{#200}}%

% west dotted arrow with lower name [1] and adjusted length
\newcommand{\wdotaR}[1]%
{\hspace{\SOURCE\unitlength}%
\hcasE{\WDOTAR}{#1}{\ARROWLENGTH}}%

% west monomorphism with length [1]x100 units
\newcommand{\wmonov}[1]{\hcase{\WMONO}{#100}}%

% west monomorphism with adjusted length
\newcommand{\wmono}%
{\hspace{\SOURCE\unitlength}%
\hcase{\WMONO}{\ARROWLENGTH}}%

% west monomorphism with upper name [1] and length [2]x100 units
\newcommand{\Wmonov}[2]{\Hcase{\WMONO}{#1}{#200}}%

% west monomorphism with upper name [1] and adjusted length
\newcommand{\Wmono}[1]%
{\hspace{\SOURCE\unitlength}%
\Hcase{\WMONO}{#1}{\ARROWLENGTH}}%

% west monomorphism with lower name [1] and length [2]x100 units
\newcommand{\wmonOv}[2]{\hcasE{\WMONO}{#1}{#200}}%

% west monomorphism with lower name [1] and adjusted length
\newcommand{\wmonO}[1]%
{\hspace{\SOURCE\unitlength}%
\hcasE{\WMONO}{#1}{\ARROWLENGTH}}%

% west epimorphism with length [1]x100 units
\newcommand{\wepiv}[1]{\hcase{\WEPI}{#100}}%

% west epimorphism with adjusted length
\newcommand{\wepi}%
{\hspace{\SOURCE\unitlength}%
\hcase{\WEPI}{\ARROWLENGTH}}%

% west epimorphism with upper name [1] and length [2]x100 units
\newcommand{\Wepiv}[2]{\Hcase{\WEPI}{#1}{#200}}%

% west epimorphism with upper name [1] and adjusted length
\newcommand{\Wepi}[1]%
{\hspace{\SOURCE\unitlength}%
\Hcase{\WEPI}{#1}{\ARROWLENGTH}}%

% west epimorphism with lower name [1] and length [2]x100 units
\newcommand{\wepIv}[2]{\hcasE{\WEPI}{#1}{#200}}%

% west epimorphism with lower name [1] and adjusted length
\newcommand{\wepI}[1]%
{\hspace{\SOURCE\unitlength}%
\hcasE{\WEPI}{#1}{\ARROWLENGTH}}%

% west bimorphism with length [1]x100 units
\newcommand{\wbimov}[1]{\hcase{\WBIMO}{#100}}%

% west bimorphism with adjusted length
\newcommand{\wbimo}%
{\hspace{\SOURCE\unitlength}%
\hcase{\WBIMO}{\ARROWLENGTH}}%

% west bimorphism with upper name [1] and length [2]x100 units
\newcommand{\Wbimov}[2]{\Hcase{\WBIMO}{#1}{#200}}%

% west bimorphism with upper name [1] and adjusted length
\newcommand{\Wbimo}[1]%
{\hspace{\SOURCE\unitlength}%
\Hcase{\WBIMO}{#1}{\ARROWLENGTH}}%

% west bimorphism with lower name [1] and length [2]x100 units
\newcommand{\wbimOv}[2]{\hcasE{\WBIMO}{#1}{#200}}%

% west bimorphism with lower name [1] and adjusted length
\newcommand{\wbimO}[1]%
{\hspace{\SOURCE\unitlength}%
\hcasE{\WBIMO}{#1}{\ARROWLENGTH}}%

% west isomorphism with length [1]x100 units
\newcommand{\wisov}[1]{\Hcase{\WAR}{\cong}{#100}}%

% west isomorphism with adjusted length
\newcommand{\wiso}%
{\hspace{\SOURCE\unitlength}%
\Hcase{\WAR}{\cong}{\ARROWLENGTH}}%

% west isomorphism with upper name [1] and length [2]x100 units
\newcommand{\Wisov}[2]{\Hcase{\WAR}{\cong#1}{#200}}%

% west isomorphism with upper name [1] and adjusted length
\newcommand{\Wiso}[1]%
{\hspace{\SOURCE\unitlength}%
\Hcase{\WAR}{#1}{\ARROWLENGTH}}%

% west isomorphism with lower name [1] and length [2]x100 units
\newcommand{\wisOv}[2]{\hcasE{\WAR}{\cong#1}{#200}}%

% west isomorphism with lower name [1] and adjusted length
\newcommand{\wisO}[1]%
{\hspace{\SOURCE\unitlength}%
\hcasE{\WAR}{#1}{\ARROWLENGTH}}%

% pair of west arrows of length [1]x100 units
\newcommand{\wbiarv}[1]{\hbicase{\WBIAR}{#100}}%

% pair of west arrows of adjusted length
\newcommand{\wbiar}%
{\hspace{\SOURCE\unitlength}%
\hbicase{\WBIAR}{\ARROWLENGTH}}%

% pair of west arrows with names [1] [2] and length [3]x100 units
\newcommand{\Wbiarv}[3]{\Hbicase{\WBIAR}{#1}{#2}{#300}}%

% pair of west arrows with names [1] [2] and adjusted length
\newcommand{\Wbiar}[2]%
{\hspace{\SOURCE\unitlength}%
\Hbicase{\WBIAR}{#1}{#2}{\ARROWLENGTH}}%

% west equality of length [1]x100 units
\newcommand{\weqlv}[1]{\hbicase{\EEQL}{#100}}%

% west equality of adjusted length
\newcommand{\weql}%
{\hspace{\SOURCE\unitlength}%
\hbicase{\EEQL}{\ARROWLENGTH}}%

% west pair of adjoint arrows of length [1]x100 units
\newcommand{\wadjarv}[1]{\hbicase{\WADJAR}{#100}}%

% west pair of adjoint arrows of adjusted length
\newcommand{\wadjar}%
{\hspace{\SOURCE\unitlength}%
\hbicase{\WADJAR}{\ARROWLENGTH}}%

% west pair of adjoint arrows with names [1][2] and length [3]x100units
\newcommand{\Wadjarv}[3]{\Hbicase{\WADJAR}{#1}{#2}{#300}}%

% west pair of adjoint arrows with names [1]J[2] and adjusted length
\newcommand{\Wadjar}[2]%
{\hspace{\SOURCE\unitlength}%
\Hbicase{\WADJAR}{#1}{#2}{\ARROWLENGTH}}%

% MACROS FOR DRAWING VERTICAL PICTURES

% \vcase{P}{n} draws the vertical picture P with length n units.
\newcommand{\vcase}[2]{#1{#2}}%

% \Vcase{P}{f}{n} draws the vertical picture P
% with left name f and length n units.
\newcommand{\Vcase}[3]{\makebox[0pt]%
{\makebox[0pt][r]{\raisebox{0pt}[0pt][0pt]{${#2}\hspace{2pt}$}}}#1{#3}}%

% \vcasE{P}{f}{n} draws the vertical picture P
% with right name f and length n units.
\newcommand{\vcasE}[3]{\makebox[0pt]%
{#1{#3}\makebox[0pt][l]{\raisebox{0pt}[0pt][0pt]{\hspace{2pt}$#2$}}}}%

% \vbicase{P}{n} draws the vertical bi-picture P with length n units.
\newcommand{\vbicase}[2]{\makebox[0pt]{{#1{#2}}}}%

% \Vbicase{P}{f}{g}{n} draws the vertical bi-picture P
% with names f, g and length n units.
\newcommand{\Vbicase}[4]{\makebox[0pt]%
{\makebox[0pt][r]{\raisebox{0pt}[0pt][0pt]{$#2$\hspace{4pt}}}#1{#4}%
\makebox[0pt][l]{\raisebox{0pt}[0pt][0pt]{\hspace{5pt}$#3$}}}}%

% SOUTH ARROWS

% \SAR{n} draws a south arrow of length n units
% and centers it in a box of width  0pt and height 0pt
\newcommand{\SAR}[1]%
{\begin{picture}(0,0)%
\put(0,0){\makebox(0,0)%
{\begin{picture}(0,#1)%
\put(0,#1){\vector(0,-1){#1}}%
\end{picture}}}\end{picture}}%

% \SDOTAR{n} draws a south dotted arrow of length n units
% and centers it in a box of width 0pt and height 0pt
\newcommand{\SDOTAR}[1]%
{\truex{100}\truey{300}%
\NUMBEROFDOTS=#1%
\divide\NUMBEROFDOTS by \value{y}%
\begin{picture}(0,0)%
\put(0,0){\makebox(0,0)%
{\begin{picture}(0,#1)%
\multiput(0,#1)(0,-\value{y}){\NUMBEROFDOTS}%
{\circle*{\value{x}}}%
\put(0,0){\vector(0,-1){0}}%
\end{picture}}}\end{picture}}%

% \SMONO{n} draws a south monomorphism of length n units
% and centers it in a box of width 0pt and height 0pt
\newcommand{\SMONO}[1]%
{\begin{picture}(0,0)%
\put(0,0){\makebox(0,0)%
{\begin{picture}(0,#1)%
\put(0,#1){\vector(0,-1){#1}}%
\truex{300}\truey{600}%
\put(0,#1){\begin{picture}(0,0)%
\put(-\value{x},-\value{x}){\line(1,0){\value{y}}}\end{picture}}%
\end{picture}}}\end{picture}}%

% \SEPI{n} draws a south epimorphism of length n units
% and centers it in a box of width 0pt and height 0pt
\newcommand{\SEPI}[1]%
{\begin{picture}(0,0)%
\put(0,0){\makebox(0,0)%
{\begin{picture}(0,#1)%
\put(0,#1){\vector(0,-1){#1}}%
\truex{300}\truey{600}\truez{800}%
\put(-\value{x},\value{z}){\line(1,0){\value{y}}}%
\end{picture}}}\end{picture}}%

% \SBIMO{n} draws a south bimorphism of length n units
% and centers it in a box of width 0pt and height 0pt
\newcommand{\SBIMO}[1]%
{\begin{picture}(0,0)%
\put(0,0){\makebox(0,0)%
{\begin{picture}(0,#1)%
\put(0,#1){\vector(0,-1){#1}}%
\truex{300}\truey{600}\truez{800}%
\put(0,#1){\begin{picture}(0,0)%
\put(-\value{x},-\value{x}){\line(1,0){\value{y}}}\end{picture}}%
\put(-\value{x},\value{z}){\line(1,0){\value{y}}}%
\end{picture}}}\end{picture}}%

% \SBIAR{n} draws a pair of south arrows of length n units
% and centers it in a box of width 0pt and height 0pt
\newcommand{\SBIAR}[1]%
{\begin{picture}(0,0)%
\truex{200}%
\put(0,0){\makebox(0,0)%
{\begin{picture}(0,#1)\put(-\value{x},#1){\vector(0,-1){#1}}%
\put(\value{x},#1){\vector(0,-1){#1}}%
\end{picture}}}\end{picture}}%

% \SEQL{n} draws a vertical equality of length n units
% and centers it in a box of width 0pt and height 0pt
\newcommand{\SEQL}[1]%
{\begin{picture}(0,0)%
\truex{100}%
\put(0,0){\makebox(0,0)%
{\begin{picture}(0,#1)\put(-\value{x},#1){\line(0,-1){#1}}%
\put(\value{x},#1){\line(0,-1){#1}}%
\end{picture}}}\end{picture}}%

% \SADJAR{n} draws a pair of vertical adjoint arrows of length n units
% and centers it in a box of width 0pt and height 0pt
\newcommand{\SADJAR}[1]{\begin{picture}(0,0)%
\truex{200}%
\put(0,0){\makebox(0,0)%
{\begin{picture}(0,#1)\put(-\value{x},#1){\vector(0,-1){#1}}%
\put(\value{x},0){\vector(0,1){#1}}%
\end{picture}}}\end{picture}}%

% All the following commands produce south arrows
% centered in a box of width O pt and height 0pt

% south arrow of length [1]x100 units
\newcommand{\sarv}[1]{\vcase{\SAR}{#100}}%

% south arrow of length 5000 units
\newcommand{\sar}{\sarv{50}}%

% south arrow with left name [1] and length [2]x100 units
\newcommand{\Sarv}[2]{\Vcase{\SAR}{#1}{#200}}%

% south arrow with left name [1] and length 5000 units
\newcommand{\Sar}[1]{\Sarv{#1}{50}}%

% south arrow with right name [1] and length [2]x100 units
\newcommand{\saRv}[2]{\vcasE{\SAR}{#1}{#200}}%

% south arrow with right name [1] and length 5000 units
\newcommand{\saR}[1]{\saRv{#1}{50}}%

% south dotted arrow of length [1]x100 units
\newcommand{\sdotarv}[1]{\vcase{\SDOTAR}{#100}}%

% south dotted arrow of length 5000 units
\newcommand{\sdotar}{\sdotarv{50}}%

% south dotted arrow with left name [1] and length [2]x100 units
\newcommand{\Sdotarv}[2]{\Vcase{\SDOTAR}{#1}{#200}}%

% south dotted arrow with left name [1] and length 5000 units
\newcommand{\Sdotar}[1]{\Sdotarv{#1}{50}}%

% south dotted arrow with right name [1] and length [2]x100 units
\newcommand{\sdotaRv}[2]{\vcasE{\SDOTAR}{#1}{#200}}%

% south dotted arrow with right name [1] and length 5000 units
\newcommand{\sdotaR}[1]{\sdotaRv{#1}{50}}%

% south monomorphism of length [1]x100 units
\newcommand{\smonov}[1]{\vcase{\SMONO}{#100}}%

% south monomorphism of length 5000 units
\newcommand{\smono}{\smonov{50}}%

% south monomorphism with left name [1] and length [2]x100 units
\newcommand{\Smonov}[2]{\Vcase{\SMONO}{#1}{#200}}%

% south monomorphism with left name [1] and length 5000 units
\newcommand{\Smono}[1]{\Smonov{#1}{50}}%

% south monomorphism with right name [1] and length [2]x100 units
\newcommand{\smonOv}[2]{\vcasE{\SMONO}{#1}{#200}}%

% south monomorphism with right name [1] and length 5000 units
\newcommand{\smonO}[1]{\smonOv{#1}{50}}%

% south epimorphism of length [1]x100 units
\newcommand{\sepiv}[1]{\vcase{\SEPI}{#100}}%

% south epimorphism of length 5000 units
\newcommand{\sepi}{\sepiv{50}}%

% south epimorphism with left name [1] and length [2]x100 units
\newcommand{\Sepiv}[2]{\Vcase{\SEPI}{#1}{#200}}%

% south epimorphism with left name [1] and length 5000 units
\newcommand{\Sepi}[1]{\Sepiv{#1}{50}}%

% south epimorphism with right name [1] and length [2]x100 units
\newcommand{\sepIv}[2]{\vcasE{\SEPI}{#1}{#200}}%

% south epimorphism with right name [1] and length 5000 units
\newcommand{\sepI}[1]{\sepIv{#1}{50}}%

% south bimorphism of length [1]x100 units
\newcommand{\sbimov}[1]{\vcase{\SBIMO}{#100}}%

% south bimorphism of length 5000 units
\newcommand{\sbimo}{\sbimov{50}}%

% south bimorphism with left name [1] and length [2]x100 units
\newcommand{\Sbimov}[2]{\Vcase{\SBIMO}{#1}{#200}}%

% south bimorphism with left name [1] and length 5000 units
\newcommand{\Sbimo}[1]{\Sbimov{#1}{50}}%

% south bimorphism with right name [1] and length [2]x100 units
\newcommand{\sbimOv}[2]{\vcasE{\SBIMO}{#1}{#200}}%

% south bimorphism with right name [1] and length 5000 units
\newcommand{\sbimO}[1]{\sbimOv{#1}{50}}%

% south isomorphism of length [1]x100 units
\newcommand{\sisov}[1]{\vcasE{\SAR}{\cong}{#100}}%

% south isomorphism of length 5000 units
\newcommand{\siso}{\sisov{50}}%

% south isomorphism with left name [1] and length [2]x100 units
\newcommand{\Sisov}[2]%
{\Vbicase{\SAR}{#1\hspace{-2pt}}{\hspace{-2pt}\cong}{#200}}%

% south isomorphism with left name [1] and length 5000 units
\newcommand{\Siso}[1]{\Sisov{#1}{50}}%

% pair of south arrows of length [1]x100 units
\newcommand{\sbiarv}[1]{\vbicase{\SBIAR}{#100}}%

% pair of south arrows of length 5000 units
\newcommand{\sbiar}{\sbiarv{50}}%

% pair of south arrows with names [1],[2] and length [3]x100 units
\newcommand{\Sbiarv}[3]{\Vbicase{\SBIAR}{#1}{#2}{#300}}%

% pair of south arrows with names [1],[2] and length 5000 units
\newcommand{\Sbiar}[2]{\Sbiarv{#1}{#2}{50}}%

% south equality of length [1]x100 units
\newcommand{\seqlv}[1]{\vbicase{\SEQL}{#100}}%

% south equality of length 50 pt
\newcommand{\seql}{\seqlv{50}}%

% south pair of adjoint arrows of length [1]x100 units
\newcommand{\sadjarv}[1]{\vbicase{\SADJAR}{#100}}%

% south pair of adjoint arrows of length 5000 units
\newcommand{\sadjar}{\sadjarv{50}}%

% south pair of adjoint arrows with names [1],[2] and length [3]x100 units
\newcommand{\Sadjarv}[3]{\Vbicase{\SADJAR}{#1}{#2}{#300}}%

% south pair of adjoint arrows with names [1],[2] and length 5000 units
\newcommand{\Sadjar}[2]{\Sadjarv{#1}{#2}{50}}%

% NORTH ARROWS

% \NAR{n} draws a north arrow of length n pt
% and centers it in a box of width 0pt and height 0pt
\newcommand{\NAR}[1]%
{\begin{picture}(0,0)%
\put(0,0){\makebox(0,0)%
{\begin{picture}(0,#1)\put(0,0){\vector(0,1){#1}}%
\end{picture}}}\end{picture}}%

% \NDOTAR{n} draws a north dotted arrow of length n units
% and centers it in a box of width 0pt and height 0pt
\newcommand{\NDOTAR}[1]%
{\truex{100}\truey{300}%
\NUMBEROFDOTS=#1%
\divide\NUMBEROFDOTS by \value{y}%
\begin{picture}(0,0)%
\put(0,0){\makebox(0,0)%
{\begin{picture}(0,#1)%
\multiput(0,0)(0,\value{y}){\NUMBEROFDOTS}%
{\circle*{\value{x}}}%
\put(0,#1){\vector(0,1){0}}%
\end{picture}}}\end{picture}}%

% \NMONO{n} draws a north monomorphism of length n pt
% and centers it in a box of width 0pt and height 0pt
\newcommand{\NMONO}[1]%
{\begin{picture}(0,0)%
\put(0,0){\makebox(0,0)%
{\begin{picture}(0,#1)%
\put(0,0){\vector(0,1){#1}}%
\truex{300}\truey{600}%
\put(-\value{x},\value{x}){\line(1,0){\value{y}}}%
\end{picture}}}%
\end{picture}}%

% \NEPI{n} draws a north epimorphism of length n pt
% and centers it in a box of width 0pt and height 0pt
\newcommand{\NEPI}[1]%
{\begin{picture}(0,0)%
\put(0,0){\makebox(0,0)%
{\begin{picture}(0,#1)%
\put(0,0){\vector(0,1){#1}}%
\truex{300}\truey{600}\truez{800}%
\put(0,#1){\begin{picture}(0,0)%
\put(-\value{x},-\value{z}){\line(1,0){\value{y}}}\end{picture}}%
\end{picture}}}\end{picture}}%

% \NBIMO{n} draws a north bimorphism of length n pt
% and centers it in a box of width 0pt and height 0pt
\newcommand{\NBIMO}[1]%
{\begin{picture}(0,0)%
\put(0,0){\makebox(0,0)%
{\begin{picture}(0,#1)%
\put(0,0){\vector(0,1){#1}}%
\truex{300}\truey{600}\truez{800}%
\put(-\value{x},\value{x}){\line(1,0){\value{y}}}%
\put(0,#1){\begin{picture}(0,0)%
\put(-\value{x},-\value{z}){\line(1,0){\value{y}}}\end{picture}}%
\end{picture}}}\end{picture}}%

% \NBIAR{n} draws a pair of north arrows of length n pt
% and centers it in a box of width 0pt and height 0pt
\newcommand{\NBIAR}[1]%
{\begin{picture}(0,0)%
\truex{200}%
\put(0,0){\makebox(0,0)%
{\begin{picture}(0,#1)\put(-\value{x},0){\vector(0,1){#1}}%
\put(\value{x},0){\vector(0,1){#1}}%
\end{picture}}}\end{picture}}%

% \NADJAR{n} draws a pair of vertical adjoint arrows of length n pt
% and centers it in a box of width 0pt and height 0pt
\newcommand{\NADJAR}[1]{\begin{picture}(0,0)%
\truex{200}%
\put(0,0){\makebox(0,0)%
{\begin{picture}(0,#1)\put(\value{x},#1){\vector(0,-1){#1}}%
\put(-\value{x},0){\vector(0,1){#1}}%
\end{picture}}}\end{picture}}%

% All the following commands produce north arrows
% centered in a box of width O pt and height 0 pt

% north arrow of length [1]x100 units
\newcommand{\narv}[1]{\vcase{\NAR}{#100}}%

% north arrow of length 5000 units
\newcommand{\nar}{\narv{50}}%

% north arrow with left name [1] and length [2]x100 units
\newcommand{\Narv}[2]{\Vcase{\NAR}{#1}{#200}}%

% north arrow with left name [1] and length 5000 units
\newcommand{\Nar}[1]{\Narv{#1}{50}}%

% north arrow with right name [1] and length [2]x100 units
\newcommand{\naRv}[2]{\vcasE{\NAR}{#1}{#200}}%

% north arrow with right name [1] and length 5000 units
\newcommand{\naR}[1]{\naRv{#1}{50}}%

% north dotted arrow of length [1]x100 units
\newcommand{\ndotarv}[1]{\vcase{\NDOTAR}{#100}}%

% north dotted arrow of length 5000 units
\newcommand{\ndotar}{\ndotarv{50}}%

% north dotted arrow with left name [1] and length [2]x100 units
\newcommand{\Ndotarv}[2]{\Vcase{\NDOTAR}{#1}{#200}}%

% north dotted arrow with left name [1] and length 5000 units
\newcommand{\Ndotar}[1]{\Ndotarv{#1}{50}}%

% north dotted arrow with right name [1] and length [2]x100 units
\newcommand{\ndotaRv}[2]{\vcasE{\NDOTAR}{#1}{#200}}%

% north dotted arrow with right name [1] and length 5000 units
\newcommand{\ndotaR}[1]{\ndotaRv{#1}{50}}%

% north monomorphism of length [1]x100 units
\newcommand{\nmonov}[1]{\vcase{\NMONO}{#100}}%

% north monomorphism of length 5000 units
\newcommand{\nmono}{\nmonov{50}}%

% north monomorphism with left name [1] and length [2]x100 units
\newcommand{\Nmonov}[2]{\Vcase{\NMONO}{#1}{#200}}%

% north monomorphism with left name [1] and length 5000 units
\newcommand{\Nmono}[1]{\Nmonov{#1}{50}}%

% north monomorphism with right name [1] and length [2]x100 units
\newcommand{\nmonOv}[2]{\vcasE{\NMONO}{#1}{#200}}%

% north monomorphism with right name [1] and length 5000 units
\newcommand{\nmonO}[1]{\nmonOv{#1}{50}}%

% north epimorphism of length [1]x100 units
\newcommand{\nepiv}[1]{\vcase{\NEPI}{#100}}%

% north epimorphism of length 5000 units
\newcommand{\nepi}{\nepiv{50}}%

% north epimorphism with left name [1] and length [2]x100 units
\newcommand{\Nepiv}[2]{\Vcase{\NEPI}{#1}{#200}}%

% north epimorphism with left name [1] and length 5000 units
\newcommand{\Nepi}[1]{\Nepiv{#1}{50}}%

% north epimorphism with right name [1] and length [2]x100 units
\newcommand{\nepIv}[2]{\vcasE{\NEPI}{#1}{#200}}%

% north epimorphism with right name [1] and length 5000 units
\newcommand{\nepI}[1]{\nepIv{#1}{50}}%

% north bimorphism of length [1]x100 units
\newcommand{\nbimov}[1]{\vcase{\NBIMO}{#100}}%

% north bimorphism of length 5000 units
\newcommand{\nbimo}{\nbimov{50}}%

% north bimorphism with left name [1] and length [2]x100 units
\newcommand{\Nbimov}[2]{\Vcase{\NBIMO}{#1}{#200}}%

% north bimorphism with left name [1] and length 5000 units
\newcommand{\Nbimo}[1]{\Nbimov{#1}{50}}%

% north bimorphism with right name [1] and length [2]x100 units
\newcommand{\nbimOv}[2]{\vcasE{\NBIMO}{#1}{#200}}%

% north bimorphism with right name [1] and length 5000 units
\newcommand{\nbimO}[1]{\nbimOv{#1}{50}}%

% north isomorphism of length [1]x100 units
\newcommand{\nisov}[1]{\vcasE{\NAR}{\cong}{#100}}%

% north isomorphism of length 5000 units
\newcommand{\niso}{\nisov{50}}%

% north isomorphism with left name [1] and length [2]x100 units
\newcommand{\Nisov}[2]%
{\Vbicase{\NAR}{#1\hspace{-2pt}}{\hspace{-2pt}\cong}{#200}}%

% north isomorphism with left name [1] and length 5000 units
\newcommand{\Niso}[1]{\Nisov{#1}{50}}%

% pair of north arrows of length [1]x100 units
\newcommand{\nbiarv}[1]{\vbicase{\NBIAR}{#100}}%

% pair of north arrows of length 5000 units
\newcommand{\nbiar}{\nbiarv{50}}%

% pair of north arrows with names [1],[2] and length [3]x100 units
\newcommand{\Nbiarv}[3]{\Vbicase{\NBIAR}{#1}{#2}{#300}}%

% pair of north arrows with names [1],[2] and length 5000 units
\newcommand{\Nbiar}[2]{\Nbiarv{#1}{#2}{50}}%

% north equality of length [1]x100 units
\newcommand{\neqlv}[1]{\vbicase{\SEQL}{#100}}%

% north equality of length 50 pt
\newcommand{\neql}{\neqlv{50}}%

% north pair of adjoint arrows of length [1]x100 units
\newcommand{\nadjarv}[1]{\vbicase{\NADJAR}{#100}}%

% north pair of adjoint arrows of length 5000 units
\newcommand{\nadjar}{\nadjarv{50}}%

% north pair of adjoint arrows with names [1],[2] and length [3]x100 units
\newcommand{\Nadjarv}[3]{\Vbicase{\NADJAR}{#1}{#2}{#300}}%

% north pair of adjoint arrows with names [1],[2] and length 5000 units
\newcommand{\Nadjar}[2]{\Nadjarv{#1}{#2}{50}}%

% MACROS FOR  FIRST DIAGONAL PICTURES

% \fdcase{P}{f}{g} draws the picture P with names f, g
\newcommand{\fdcase}[3]{\begin{picture}(0,0)%
\put(0,-150){#1}%
\truex{200}\truey{600}\truez{600}%
\put(-\value{x},-\value{x}){\makebox(0,\value{z})[r]{${#2}$}}%
\put(\value{x},-\value{y}){\makebox(0,\value{z})[l]{${#3}$}}%
\end{picture}}%

% \fdbicase{P}{f}{g} draws the bipicture P with names f, g
\newcommand{\fdbicase}[3]{\begin{picture}(0,0)%
\put(0,-150){#1}%
\truex{800}\truey{50}%
\put(-\value{x},\value{y}){${#2}$}%
\truex{200}\truey{950}%
\put(\value{x},-\value{y}){${#3}$}%
\end{picture}}%

% NORTH-EAST ARROWS

% \NEAR draws a north-east arrow
\newcommand{\NEAR}{\begin{picture}(0,0)%
\put(-2900,-2900){\vector(1,1){5800}}%
\end{picture}}%

% \NEDOTAR draws a north-east dotted arrow
\newcommand{\NEDOTAR}%
{\truex{100}\truey{212}%
\NUMBEROFDOTS=5800%
\divide\NUMBEROFDOTS by \value{y}%
\begin{picture}(0,0)%
\multiput(-2900,-2900)(\value{y},\value{y}){\NUMBEROFDOTS}%
{\circle*{\value{x}}}%
\put(2900,2900){\vector(1,1){0}}%
\end{picture}}%

% \NEMONO draws a north-east monomorphism
\newcommand{\NEMONO}{\begin{picture}(0,0)%
\put(-2900,-2900){\vector(1,1){5800}}%
\put(-2900,-2900){\begin{picture}(0,0)%
\truex{141}%
\put(\value{x},\value{x}){\makebox(0,0){$\times$}}%
\end{picture}}\end{picture}}%

% \NEEPI draws a north-east epimorphism
\newcommand{\NEEPI}{\begin{picture}(0,0)%
\put(-2900,-2900){\vector(1,1){5800}}%
\put(2900,2900){\begin{picture}(0,0)%
\truex{545}%
\put(-\value{x},-\value{x}){\makebox(0,0){$\times$}}%
\end{picture}}\end{picture}}%

% \NEBIMO draws a north-east bimorphism
\newcommand{\NEBIMO}{\begin{picture}(0,0)%
\put(-2900,-2900){\vector(1,1){5800}}%
\put(2900,2900){\begin{picture}(0,0)%
\truex{545}%
\put(-\value{x},-\value{x}){\makebox(0,0){$\times$}}%
\end{picture}}
\put(-2900,-2900){\begin{picture}(0,0)%
\truex{141}%
\put(\value{x},\value{x}){\makebox(0,0){$\times$}}%
\end{picture}}\end{picture}}%

% \NEBIAR draws a pair of north-east arrows
\newcommand{\NEBIAR}{\begin{picture}(0,0)%
\put(-2900,-2900){\begin{picture}(0,0)%
\truex{141}%
\put(-\value{x},\value{x}){\vector(1,1){5800}}%
\put(\value{x},-\value{x}){\vector(1,1){5800}}%
\end{picture}}\end{picture}}%

% \NEEQL draws a north-east equality
\newcommand{\NEEQL}{\begin{picture}(0,0)%
\put(-2900,-2900){\begin{picture}(0,0)%
\truex{70}%
\put(-\value{x},\value{x}){\line(1,1){5800}}%
\put(\value{x},-\value{x}){\line(1,1){5800}}%
\end{picture}}\end{picture}}%

% \NEADJAR draws a north-east pair of adjoint arrows
\newcommand{\NEADJAR}{\begin{picture}(0,0)%
\put(-2900,-2900){\begin{picture}(0,0)%
\truex{141}%
\put(\value{x},-\value{x}){\vector(1,1){5800}}%
\end{picture}}%
\put(2900,2900){\begin{picture}(0,0)%
\truex{141}%
\put(-\value{x},\value{x}){\vector(-1,-1){5800}}%
\end{picture}}\end{picture}}%

% \NEARV{n} draws a north-east arrow of length nx100 units
\newcommand{\NEARV}[1]{\begin{picture}(0,0)%
\put(0,0){\makebox(0,0){\begin{picture}(#1,#1)%
\put(0,0){\vector(1,1){#1}}\end{picture}}}%
\end{picture}}%

% All the following commands draw north-east arrows of horizontal
% extent 5800 units and of formal dimensions (0,0).

% north-east arrow
\newcommand{\near}{\fdcase{\NEAR}{}{}}%

% north-east arrow with upper name [1]
\newcommand{\Near}[1]{\fdcase{\NEAR}{#1}{}}%

% north-east arrow with lower name [1]
\newcommand{\neaR}[1]{\fdcase{\NEAR}{}{#1}}%

% north-east dotted arrow
\newcommand{\nedotar}{\fdcase{\NEDOTAR}{}{}}%

% north-east dotted arrow with upper name [1]
\newcommand{\Nedotar}[1]{\fdcase{\NEDOTAR}{#1}{}}%

% north-east dotted arrow with lower name [1]
\newcommand{\nedotaR}[1]{\fdcase{\NEDOTAR}{}{#1}}%

% north-east monomorphism
\newcommand{\nemono}{\fdcase{\NEMONO}{}{}}%

% north-est monomorphism with upper name [1]
\newcommand{\Nemono}[1]{\fdcase{\NEMONO}{#1}{}}%

% north-east monomorphism with lower name [1]
\newcommand{\nemonO}[1]{\fdcase{\NEMONO}{}{#1}}%

% north-east epimorphism
\newcommand{\neepi}{\fdcase{\NEEPI}{}{}}%

% north-east epimorphism with upper name [1]
\newcommand{\Neepi}[1]{\fdcase{\NEEPI}{#1}{}}%

% north-east epimorphism with lower name [1]
\newcommand{\neepI}[1]{\fdcase{\NEEPI}{}{#1}}%

% north-east bimorphism
\newcommand{\nebimo}{\fdcase{\NEBIMO}{}{}}%

% north-east bimorphism with upper name [1]
\newcommand{\Nebimo}[1]{\fdcase{\NEBIMO}{#1}{}}%

% north-east bimorphism with lower name [1]
\newcommand{\nebimO}[1]{\fdcase{\NEBIMO}{}{#1}}%

% north-east isomorphism
\newcommand{\neiso}{\fdcase{\NEAR}{\hspace{-2pt}\cong}{}}%

% north-east isomorphism with name [1]
\newcommand{\Neiso}[1]{\fdcase{\NEAR}{\hspace{-2pt}\cong}{#1}}%

% pair of north-east arrows
\newcommand{\nebiar}{\fdbicase{\NEBIAR}{}{}}%

% pair of north-east arrows with names [1], [2]
\newcommand{\Nebiar}[2]{\fdbicase{\NEBIAR}{#1}{#2}}%

% north-east equality
\newcommand{\neeql}{\fdbicase{\NEEQL}{}{}}%

% north-east pair of adjoint arrows
\newcommand{\neadjar}{\fdbicase{\NEADJAR}{}{}}%

% north-east pair of adjoint arrows with names [1], [2]
\newcommand{\Neadjar}[2]{\fdbicase{\NEADJAR}{#1}{#2}}%

% the following commands produce a north-east arrow of variable length
% and of formal dimensions (0,0).

% north-east arrow of horizontal length [1]x100 units
\newcommand{\nearv}[1]{\fdcase{\NEARV{#100}}{}{}}%

% north-east arrow with upper name [1] and horizontal length [2]x100 units
\newcommand{\Nearv}[2]{\fdcase{\NEARV{#200}}{#1}{}}%

% north-east arrow with lower name [1] and horizontal length [2]x100 units
\newcommand{\neaRv}[2]{\fdcase{\NEARV{#200}}{}{#1}}%

% SOUTH-WEST ARROWS

% \SWAR draws a south-west arrow
\newcommand{\SWAR}{\begin{picture}(0,0)%
\put(2900,2900){\vector(-1,-1){5800}}%
\end{picture}}%

% \SWDOTAR draws a south-west dotted arrow
\newcommand{\SWDOTAR}%
{\truex{100}\truey{212}%
\NUMBEROFDOTS=5800%
\divide\NUMBEROFDOTS by \value{y}%
\begin{picture}(0,0)%
\multiput(2900,2900)(-\value{y},-\value{y}){\NUMBEROFDOTS}%
{\circle*{\value{x}}}%
\put(-2900,-2900){\vector(-1,-1){0}}%
\end{picture}}%

% \SWMONO draws a south-west monomorphism
\newcommand{\SWMONO}{\begin{picture}(0,0)%
\put(2900,2900){\vector(-1,-1){5800}}%
\put(2900,2900){\begin{picture}(0,0)%
\truex{141}%
\put(-\value{x},-\value{x}){\makebox(0,0){$\times$}}%
\end{picture}}\end{picture}}%

% \SWEPI draws a south-west epimorphism
\newcommand{\SWEPI}{\begin{picture}(0,0)%
\put(2900,2900){\vector(-1,-1){5800}}%
\put(-2900,-2900){\begin{picture}(0,0)%
\truex{525}%
\put(\value{x},\value{x}){\makebox(0,0){$\times$}}%
\end{picture}}\end{picture}}%

% \SWBIMO draws a south-west bimorphism
\newcommand{\SWBIMO}{\begin{picture}(0,0)%
\put(2900,2900){\vector(-1,-1){5800}}%
\put(2900,2900){\begin{picture}(0,0)%
\truex{141}%
\put(-\value{x},-\value{x}){\makebox(0,0){$\times$}}%
\end{picture}}%
\put(-2900,-2900){\begin{picture}(0,0)%
\truex{525}%
\put(\value{x},\value{x}){\makebox(0,0){$\times$}}%
\end{picture}}\end{picture}}%

% \SWBIAR draws a pair of south-west arrows
\newcommand{\SWBIAR}{\begin{picture}(0,0)%
\put(2900,2900){\begin{picture}(0,0)%
\truex{141}%
\put(\value{x},-\value{x}){\vector(-1,-1){5800}}%
\put(-\value{x},\value{x}){\vector(-1,-1){5800}}%
\end{picture}}\end{picture}}%

% \SWADJAR draws a south-west pair of adjoint arrows
\newcommand{\SWADJAR}{\begin{picture}(0,0)%
\put(-2900,-2900){\begin{picture}(0,0)%
\truex{141}%
\put(-\value{x},\value{x}){\vector(1,1){5800}}%
\end{picture}}%
\put(2900,2900){\begin{picture}(0,0)%
\truex{141}%
\put(\value{x},-\value{x}){\vector(-1,-1){5800}}%
\end{picture}}\end{picture}}%

% \SWARV{n} draws a south-west arrow of length nx100 units
\newcommand{\SWARV}[1]{\begin{picture}(0,0)%
\put(0,0){\makebox(0,0){\begin{picture}(#1,#1)%
\put(#1,#1){\vector(-1,-1){#1}}\end{picture}}}%
\end{picture}}%

% All the following commands draw south-west arrows of horizontal
% extent 5800 units and formal dimensions (0,0).

% south-west arrow
\newcommand{\swar}{\fdcase{\SWAR}{}{}}%

% south-west arrow with upper name [1]
\newcommand{\Swar}[1]{\fdcase{\SWAR}{#1}{}}%

% south-west arrow with lower name [1]
\newcommand{\swaR}[1]{\fdcase{\SWAR}{}{#1}}%

% south-west dotted arrow
\newcommand{\swdotar}{\fdcase{\SWDOTAR}{}{}}%

% south-west dotted arrow with upper name [1]
\newcommand{\Swdotar}[1]{\fdcase{\SWDOTAR}{#1}{}}%

% south-west dotted arrow with lower name [1]
\newcommand{\swdotaR}[1]{\fdcase{\SWDOTAR}{}{#1}}%

% south-west monomorphism
\newcommand{\swmono}{\fdcase{\SWMONO}{}{}}%

% north-est monomorphism with upper name [1]
\newcommand{\Swmono}[1]{\fdcase{\SWMONO}{#1}{}}%

% south-west monomorphism with lower name [1]
\newcommand{\swmonO}[1]{\fdcase{\SWMONO}{}{#1}}%

% south-west epimorphism
\newcommand{\swepi}{\fdcase{\SWEPI}{}{}}%

% south-west epimorphism with upper name [1]
\newcommand{\Swepi}[1]{\fdcase{\SWEPI}{#1}{}}%

% south-west epimorphism with lower name [1]
\newcommand{\swepI}[1]{\fdcase{\SWEPI}{}{#1}}%

% south-west bimorphism
\newcommand{\swbimo}{\fdcase{\SWBIMO}{}{}}%

% south-west bimorphism with upper name [1]
\newcommand{\Swbimo}[1]{\fdcase{\SWBIMO}{#1}{}}%

% south-west bimorphism with lower name [1]
\newcommand{\swbimO}[1]{\fdcase{\SWBIMO}{}{#1}}%

% south-west isomorphism
\newcommand{\swiso}{\fdcase{\SWAR}{\hspace{-2pt}\cong}{}}%

% south-west isomorphism with name [1]
\newcommand{\Swiso}[1]{\fdcase{\SWAR}{\hspace{-2pt}\cong}{#1}}%

% pair of south-west arrows
\newcommand{\swbiar}{\fdbicase{\SWBIAR}{}{}}%

% pair of south-west arrows with names [1], [2]
\newcommand{\Swbiar}[2]{\fdbicase{\SWBIAR}{#1}{#2}}%

% south-west equality
\newcommand{\sweql}{\fdbicase{\NEEQL}{}{}}%

% south-west pair of adjoint arrows
\newcommand{\swadjar}{\fdbicase{\SWADJAR}{}{}}%

% south-west pair of adjoint arrows with names [1], [2]
\newcommand{\Swadjar}[2]{\fdbicase{\SWADJAR}{#1}{#2}}%

% the following commands produce a south-west arrow of variable length
% and formal dimensions (0,0).

% south-west arrow of horizontal length [1]x100 units
\newcommand{\swarv}[1]{\fdcase{\SWARV{#100}}{}{}}%

% south-west arrow with upper name [1] and horizontal length [2]x100 units
\newcommand{\Swarv}[2]{\fdcase{\SWARV{#200}}{#1}{}}%

% south-west arrow with lower name [1] and horizontal length [2]x100 units
\newcommand{\swaRv}[2]{\fdcase{\SWARV{#200}}{}{#1}}%

% MACROS FOR  SECOND DIAGONAL PICTURES

% \sdcase{P}{f}{g} draws the picture P with names f, g
\newcommand{\sdcase}[3]{\begin{picture}(0,0)%
\put(0,-150){#1}%
\truex{100}\truez{600}%
\put(\value{x},\value{x}){\makebox(0,\value{z})[l]{${#2}$}}%
\truex{300}\truey{800}%
\put(-\value{x},-\value{y}){\makebox(0,\value{z})[r]{${#3}$}}%
\end{picture}}%

% \sdbicase{P}{f}{g} draws the bipicture P with names f, g
\newcommand{\sdbicase}[3]{\begin{picture}(0,0)%
\put(0,-150){#1}%
\truex{250}\truey{600}\truez{850}%
\put(\value{x},\value{x}){\makebox(0,\value{y})[l]{${#2}$}}%
\put(-\value{x},-\value{z}){\makebox(0,\value{y})[r]{${#3}$}}%
\end{picture}}%

% SOUT-EAST ARROWS

% \SEAR draws a south-east arrow
\newcommand{\SEAR}{\begin{picture}(0,0)%
\put(-2900,2900){\vector(1,-1){5800}}%
\end{picture}}%

% \SEDOTAR draws a south-east dotted arrow
\newcommand{\SEDOTAR}%
{\truex{100}\truey{212}%
\NUMBEROFDOTS=5800%
\divide\NUMBEROFDOTS by \value{y}%
\begin{picture}(0,0)%
\multiput(-2900,2900)(\value{y},-\value{y}){\NUMBEROFDOTS}%
{\circle*{\value{x}}}%
\put(2900,-2900){\vector(1,-1){0}}%
\end{picture}}%

% \SEMONO draws a south-east monomorphism
\newcommand{\SEMONO}{\begin{picture}(0,0)%
\put(-2900,2900){\vector(1,-1){5800}}%
\put(-2900,2900){\begin{picture}(0,0)%
\truex{141}%
\put(\value{x},-\value{x}){\makebox(0,0){$\times$}}%
\end{picture}}\end{picture}}%

% \SEEPI draws a south-east epimorphism
\newcommand{\SEEPI}{\begin{picture}(0,0)%
\put(-2900,2900){\vector(1,-1){5800}}%
\put(2900,-2900){\begin{picture}(0,0)%
\truex{525}%
\put(-\value{x},\value{x}){\makebox(0,0){$\times$}}%
\end{picture}}\end{picture}}%

% \SEBIMO draws a south-east bimorphism
\newcommand{\SEBIMO}{\begin{picture}(0,0)%
\put(-2900,2900){\vector(1,-1){5800}}%
\put(-2900,2900){\begin{picture}(0,0)%
\truex{141}%
\put(\value{x},-\value{x}){\makebox(0,0){$\times$}}%
\end{picture}}%
\put(2900,-2900){\begin{picture}(0,0)%
\truex{525}%
\put(-\value{x},\value{x}){\makebox(0,0){$\times$}}%
\end{picture}}\end{picture}}%

% \SEBIAR draws a pair of south-east arrows
\newcommand{\SEBIAR}{\begin{picture}(0,0)%
\put(-2900,2900){\begin{picture}(0,0)%
\truex{141}
\put(-\value{x},-\value{x}){\vector(1,-1){5800}}%
\put(\value{x},\value{x}){\vector(1,-1){5800}}%
\end{picture}}\end{picture}}%

% \SEEQL draws a south-east equality
\newcommand{\SEEQL}{\begin{picture}(0,0)%
\put(-2900,2900){\begin{picture}(0,0)%
\truex{70}%
\put(-\value{x},-\value{x}){\line(1,-1){5800}}%
\put(\value{x},\value{x}){\line(1,-1){5800}}%
\end{picture}}\end{picture}}%

% \SEADJAR draws a south-east pair of adjoint arrows
\newcommand{\SEADJAR}{\begin{picture}(0,0)%
\put(-2900,2900){\begin{picture}(0,0)%
\truex{141}%
\put(-\value{x},-\value{x}){\vector(1,-1){5800}}%
\end{picture}}%
\put(2900,-2900){\begin{picture}(0,0)%
\truex{141}%
\put(\value{x},\value{x}){\vector(-1,1){5800}}%
\end{picture}}\end{picture}}%

% \SEARV{n} draws a south-east arrow of length nx100 units
\newcommand{\SEARV}[1]{\begin{picture}(0,0)%
\put(0,0){\makebox(0,0){\begin{picture}(#1,#1)%
\put(0,#1){\vector(1,-1){#1}}\end{picture}}}%
\end{picture}}%

% All the following commands draw south-east arrows of horizontal
% extent 5800 units and of formal dimensions (0,0).

% south-east arrow
\newcommand{\sear}{\sdcase{\SEAR}{}{}}%

% south-east arrow with upper name [1]
\newcommand{\Sear}[1]{\sdcase{\SEAR}{#1}{}}%

% south-east arrow with lower name [1]
\newcommand{\seaR}[1]{\sdcase{\SEAR}{}{#1}}%

% south-east dotted arrow
\newcommand{\sedotar}{\sdcase{\SEDOTAR}{}{}}%

% south-east dotted arrow with upper name [1]
\newcommand{\Sedotar}[1]{\sdcase{\SEDOTAR}{#1}{}}%

% south-east dotted arrow with lower name [1]
\newcommand{\sedotaR}[1]{\sdcase{\SEDOTAR}{}{#1}}%

% south-east monomorphism
\newcommand{\semono}{\sdcase{\SEMONO}{}{}}%

% south-east monomorphism with upper name [1]
\newcommand{\Semono}[1]{\sdcase{\SEMONO}{#1}{}}%

% south-east monomorphism with lower name [1]
\newcommand{\semonO}[1]{\sdcase{\SEMONO}{}{#1}}%

% south-east epimorphism
\newcommand{\seepi}{\sdcase{\SEEPI}{}{}}%

% south-east epimorphism with upper name [1]
\newcommand{\Seepi}[1]{\sdcase{\SEEPI}{#1}{}}%

% south-east epimorphism with lower name [1]
\newcommand{\seepI}[1]{\sdcase{\SEEPI}{}{#1}}%

% south-east bimorphism
\newcommand{\sebimo}{\sdcase{\SEBIMO}{}{}}%

% south-east bimorphism with upper name [1]
\newcommand{\Sebimo}[1]{\sdcase{\SEBIMO}{#1}{}}%

% south-east bimorphism with lower name [1]
\newcommand{\sebimO}[1]{\sdcase{\SEBIMO}{}{#1}}%

% south-east isomorphism
\newcommand{\seiso}{\sdcase{\SEAR}{\hspace{-2pt}\cong}{}}%

% south-east isomorphism with name [1]
\newcommand{\Seiso}[1]{\sdcase{\SEAR}{\hspace{-2pt}\cong}{#1}}%

% pair of south-east arrows
\newcommand{\sebiar}{\sdbicase{\SEBIAR}{}{}}%

% pair of south-east arrows with names [1], [2]
\newcommand{\Sebiar}[2]{\sdbicase{\SEBIAR}{#1}{#2}}%

% south-east equality
\newcommand{\seeql}{\sdbicase{\SEEQL}{}{}}%

% south-east pair of adjoint arrows
\newcommand{\seadjar}{\sdbicase{\SEADJAR}{}{}}%

% south-east pair of adjoint arrows with names [1], [2]
\newcommand{\Seadjar}[2]{\sdbicase{\SEADJAR}{#1}{#2}}%

% the following commands produce a south-east arrow of variable length
% and of formal dimensions (0,0).

% south-east arrow of horizontal length [1]x100 units
\newcommand{\searv}[1]{\sdcase{\SEARV{#100}}{}{}}%

% south-east arrow with upper name [1] and horizontal length [2]x100 units
\newcommand{\Searv}[2]{\sdcase{\SEARV{#200}}{#1}{}}%

% south-east arrow with lower name [1] and horizontal length [2]x100 units
\newcommand{\seaRv}[2]{\sdcase{\SEARV{#200}}{}{#1}}%

% NORTH-WEST ARROWS

% \NWAR draws a north-west arrow
\newcommand{\NWAR}{\begin{picture}(0,0)%
\put(2900,-2900){\vector(-1,1){5800}}%
\end{picture}}%

% \NWDOTAR draws a north-west dotted arrow
\newcommand{\NWDOTAR}%
{\truex{100}\truey{212}%
\NUMBEROFDOTS=5800%
\divide\NUMBEROFDOTS by \value{y}%
\begin{picture}(0,0)%
\multiput(2900,-2900)(-\value{y},\value{y}){\NUMBEROFDOTS}%
{\circle*{\value{x}}}%
\put(-2900,2900){\vector(-1,1){0}}%
\end{picture}}%

% \NWMONO draws a north-west monomorphism
\newcommand{\NWMONO}{\begin{picture}(0,0)%
\put(2900,-2900){\vector(-1,1){5800}}%
\put(2900,-2900){\begin{picture}(0,0)%
\truex{141}%
\put(-\value{x},\value{x}){\makebox(0,0){$\times$}}%
\end{picture}}\end{picture}}%

% \NWEPI draws a north-west epimorphism
\newcommand{\NWEPI}{\begin{picture}(0,0)%
\put(2900,-2900){\vector(-1,1){5800}}%
\put(-2900,2900){\begin{picture}(0,0)%
\truex{525}%
\put(\value{x},-\value{x}){\makebox(0,0){$\times$}}%
\end{picture}}\end{picture}}%

% \NWBIMO draws a north-west bimorphism
\newcommand{\NWBIMO}{\begin{picture}(0,0)%
\put(2900,-2900){\vector(-1,1){5800}}%
\put(2900,-2900){\begin{picture}(0,0)%
\truex{141}%
\put(-\value{x},\value{x}){\makebox(0,0){$\times$}}%
\end{picture}}%
\put(-2900,2900){\begin{picture}(0,0)%
\truex{525}%
\put(\value{x},-\value{x}){\makebox(0,0){$\times$}}%
\end{picture}}\end{picture}}%

% \NWBIAR draws a pair of north-west arrows
\newcommand{\NWBIAR}{\begin{picture}(0,0)%
\put(2900,-2900){\begin{picture}(0,0)%
\truex{141}%
\put(\value{x},\value{x}){\vector(-1,1){5800}}%
\end{picture}}%
\put(2900,-2900){\begin{picture}(0,0)%
\truex{141}
\put(-\value{x},-\value{x}){\vector(-1,1){5800}}%
\end{picture}}\end{picture}}%

% \NWADJAR draws a north-west pair of adjoint arrows
\newcommand{\NWADJAR}{\begin{picture}(0,0)%
\put(-2900,2900){\begin{picture}(0,0)%
\truex{141}%
\put(\value{x},\value{x}){\vector(1,-1){5800}}%
\end{picture}}%
\put(2900,-2900){\begin{picture}(0,0)%
\truex{141}%
\put(-\value{x},-\value{x}){\vector(-1,1){5800}}%
\end{picture}}\end{picture}}%

% \NWARV{n} draws a north-west arrow of length nx100 units
\newcommand{\NWARV}[1]{\begin{picture}(0,0)%
\put(0,0){\makebox(0,0){\begin{picture}(#1,#1)%
\put(#1,0){\vector(-1,1){#1}}\end{picture}}}%
\end{picture}}%

% All the following commands draw north-west arrows of horizontal
% extent 5800 units and of formal dimensions (0,0).

% north-west arrow
\newcommand{\nwar}{\sdcase{\NWAR}{}{}}%

% north-west arrow with upper name [1]
\newcommand{\Nwar}[1]{\sdcase{\NWAR}{#1}{}}%

% north-west arrow with lower name [1]
\newcommand{\nwaR}[1]{\sdcase{\NWAR}{}{#1}}%

% north-west dotted arrow
\newcommand{\nwdotar}{\sdcase{\NWDOTAR}{}{}}%

% north-west dotted arrow with upper name [1]
\newcommand{\Nwdotar}[1]{\sdcase{\NWDOTAR}{#1}{}}%

% north-west dotted arrow with lower name [1]
\newcommand{\nwdotaR}[1]{\sdcase{\NWDOTAR}{}{#1}}%

% north-west monomorphism
\newcommand{\nwmono}{\sdcase{\NWMONO}{}{}}%

% north-west monomorphism with upper name [1]
\newcommand{\Nwmono}[1]{\sdcase{\NWMONO}{#1}{}}%

% north-west monomorphism with lower name [1]
\newcommand{\nwmonO}[1]{\sdcase{\NWMONO}{}{#1}}%

% north-west epimorphism
\newcommand{\nwepi}{\sdcase{\NWEPI}{}{}}%

% north-west epimorphism with upper name [1]
\newcommand{\Nwepi}[1]{\sdcase{\NWEPI}{#1}{}}%

% north-west epimorphism with lower name [1]
\newcommand{\nwepI}[1]{\sdcase{\NWEPI}{}{#1}}%

% north-west bimorphism
\newcommand{\nwbimo}{\sdcase{\NWBIMO}{}{}}%

% north-west bimorphism with upper name [1]
\newcommand{\Nwbimo}[1]{\sdcase{\NWBIMO}{#1}{}}%

% north-west bimorphism with lower name [1]
\newcommand{\nwbimO}[1]{\sdcase{\NWBIMO}{}{#1}}%

% north-west isomorphism
\newcommand{\nwiso}{\sdcase{\NWAR}{\hspace{-2pt}\cong}{}}%

% north-west isomorphism with name [1]
\newcommand{\Nwiso}[1]{\sdcase{\NWAR}{\hspace{-2pt}\cong}{#1}}%

% pair of north-west arrows
\newcommand{\nwbiar}{\sdbicase{\NWBIAR}{}{}}%

% pair of north-west arrows with names [1], [2]
\newcommand{\Nwbiar}[2]{\sdbicase{\NWBIAR}{#1}{#2}}%

% north-west equality
\newcommand{\nweql}{\sdbicase{\SEEQL}{}{}}%

% north-west pair of adjoint arrows
\newcommand{\nwadjar}{\sdbicase{\NWADJAR}{}{}}%

% north-west pair of adjoint arrows with names [1], [2]
\newcommand{\Nwadjar}[2]{\sdbicase{\NWADJAR}{#1}{#2}}%

% the following commands produce a north-west arrow of variable length
% and of formal dimensions (0,0).

% north-west arrow of horizontal length [1]x100 units
\newcommand{\nwarv}[1]{\sdcase{\NWARV{#100}}{}{}}%

% north-west arrow with upper name [1] and horizontal length [2]x100 units
\newcommand{\Nwarv}[2]{\sdcase{\NWARV{#200}}{#1}{}}%

% north-west arrow with lower name [1] and horizontal length [2]x100 units
\newcommand{\nwaRv}[2]{\sdcase{\NWARV{#200}}{}{#1}}%

% HORIZONTAL RECTANGLE DIAGONAL ARROWS

% All the following commands produce arrows of horizontal extent 13200
% units, oriented in the directions of the diagonals of an horizontal
% rectangle of sides 2, 1. The results are positioned in a box of width
% 0pt and height 0pt.

% \ENEAR{f}{g} draws a east-north-east arrow with names f, g
\newcommand{\ENEAR}[2]%
{\makebox[0pt]{\begin{picture}(0,0)%
\put(0,-150){\makebox(0,0){\begin{picture}(0,0)%
\put(-6600,-3300){\vector(2,1){13200}}%
\truex{200}\truey{800}\truez{600}%
\put(-\value{x},\value{x}){\makebox(0,\value{z})[r]{${#1}$}}%
\put(\value{x},-\value{y}){\makebox(0,\value{z})[l]{${#2}$}}%
\end{picture}}}\end{picture}}}%

% east-north-east arrow
\newcommand{\enear}{\ENEAR{}{}}%

% east-north-east arrow with upper name [1]
\newcommand{\Enear}[1]{\ENEAR{#1}{}}%

% east-north-east arrow with lower name [1]
\newcommand{\eneaR}[1]{\ENEAR{}{#1}}%

% \ESEAR{f}{g} draws an east-south-east arrow with names f, g
\newcommand{\ESEAR}[2]%
{\makebox[0pt]{\begin{picture}(0,0)%
\put(0,-150){\makebox(0,0){\begin{picture}(0,0)%
\put(-6600,3300){\vector(2,-1){13200}}%
\truex{200}\truey{800}\truez{600}%
\put(\value{x},\value{x}){\makebox(0,\value{z})[l]{${#1}$}}%
\put(-\value{x},-\value{y}){\makebox(0,\value{z})[r]{${#2}$}}%
\end{picture}}}\end{picture}}}%

% east-south-east arrow
\newcommand{\esear}{\ESEAR{}{}}%

% east-south-east arrow with upper name [1]
\newcommand{\Esear}[1]{\ESEAR{#1}{}}%

% east-south-east arrow with lower name [1]
\newcommand{\eseaR}[1]{\ESEAR{}{#1}}%

% \WNWAR{f}{g} draws a west-north-west arrow with names f, g
\newcommand{\WNWAR}[2]%
{\makebox[0pt]{\begin{picture}(0,0)%
\put(0,-150){\makebox(0,0){\begin{picture}(0,0)%
\put(6600,-3300){\vector(-2,1){13200}}%
\truex{200}\truey{800}\truez{600}%
\put(\value{x},\value{x}){\makebox(0,\value{z})[l]{${#1}$}}%
\put(-\value{x},-\value{y}){\makebox(0,\value{z})[r]{${#2}$}}%
\end{picture}}}\end{picture}}}%

% west-north-west arrow
\newcommand{\wnwar}{\WNWAR{}{}}%

% west-north-west arrow with upper name [1]
\newcommand{\Wnwar}[1]{\WNWAR{#1}{}}%

% west-north-west arrow with lower name [1]
\newcommand{\wnwaR}[1]{\WNWAR{}{#1}}%

% \WSWAR{f}{g} draws a west-south-west arrow with names f, g
\newcommand{\WSWAR}[2]%
{\makebox[0pt]{\begin{picture}(0,0)%
\put(0,-150){\makebox(0,0){\begin{picture}(0,0)%
\put(6600,3300){\vector(-2,-1){13200}}%
\truex{200}\truey{800}\truez{600}%
\put(-\value{x},\value{x}){\makebox(0,\value{z})[r]{${#1}$}}%
\put(\value{x},-\value{y}){\makebox(0,\value{z})[l]{${#2}$}}%
\end{picture}}}\end{picture}}}%

% west-south-west arrow
\newcommand{\wswar}{\WSWAR{}{}}%

% west-south-west arrow with upper name [1]
\newcommand{\Wswar}[1]{\WSWAR{#1}{}}%

% west-south-west arrow with lower name [1]
\newcommand{\wswaR}[1]{\WSWAR{}{#1}}%

% VERTICAL RECTANGLE DIAGONAL ARROWS

% All the following commands produce arrows of horizontal extent 6600
% units, oriented in the directions of the diagonals of a vertical
% rectangle of sides 1, 2. The results are positioned in a box of width
% 0pt and height 0pt.

% \NNEAR{f}{g} draws a north-north-east arrow with names f, g
\newcommand{\NNEAR}[2]%
{\raisebox{-1pt}[0pt][0pt]{\begin{picture}(0,0)%
\put(0,0){\makebox(0,0){\begin{picture}(0,0)%
\put(-3300,-6600){\vector(1,2){6600}}%
\truex{100}\truez{600}%
\put(-\value{x},\value{x}){\makebox(0,\value{z})[r]{${#1}$}}%
\put(\value{x},-\value{z}){\makebox(0,\value{z})[l]{${#2}$}}%
\end{picture}}}\end{picture}}}%

% north-north-east arrow
\newcommand{\nnear}{\NNEAR{}{}}%

% north-north-east arrow with upper name [1]
\newcommand{\Nnear}[1]{\NNEAR{#1}{}}%

% north-north-east arrow with lower name [1]
\newcommand{\nneaR}[1]{\NNEAR{}{#1}}%

% \SSWAR{f}{g} draws a south-south-west arrow with names f, g
\newcommand{\SSWAR}[2]%
{\raisebox{-1pt}[0pt][0pt]{\begin{picture}(0,0)%
\put(0,0){\makebox(0,0){\begin{picture}(0,0)%
\put(3300,6600){\vector(-1,-2){6600}}%
\truex{100}\truez{600}%
\put(-\value{x},\value{x}){\makebox(0,\value{z})[r]{${#1}$}}%
\put(\value{x},-\value{z}){\makebox(0,\value{z})[l]{${#2}$}}%
\end{picture}}}\end{picture}}}%

% south-south-west arrow
\newcommand{\sswar}{\SSWAR{}{}}%

% south-south-west arrow with upper name [1]
\newcommand{\Sswar}[1]{\SSWAR{#1}{}}%

% south-south-west arrow with lower name [1]
\newcommand{\sswaR}[1]{\SSWAR{}{#1}}%

% \SSEAR{f}{g} draws a south-south-east arrow with names f, g
\newcommand{\SSEAR}[2]%
{\raisebox{-1pt}[0pt][0pt]{\begin{picture}(0,0)%
\put(0,0){\makebox(0,0){\begin{picture}(0,0)%
\put(-3300,6600){\vector(1,-2){6600}}%
\truex{200}\truez{600}%
\put(\value{x},\value{x}){\makebox(0,\value{z})[l]{${#1}$}}%
\put(-\value{x},-\value{z}){\makebox(0,\value{z})[r]{${#2}$}}%
\end{picture}}}\end{picture}}}%

% south-south-east arrow
\newcommand{\ssear}{\SSEAR{}{}}%

% south-south-east arrow with upper name [1]
\newcommand{\Ssear}[1]{\SSEAR{#1}{}}%

% south-south-east arrow with lower name [1]
\newcommand{\sseaR}[1]{\SSEAR{}{#1}}%

% \NNWAR{f}{g} draws a north-north-west arrow with names f, g
\newcommand{\NNWAR}[2]%
{\raisebox{-1pt}[0pt][0pt]{\begin{picture}(0,0)%
\put(0,0){\makebox(0,0){\begin{picture}(0,0)%
\put(3300,-6600){\vector(-1,2){6600}}%
\truex{200}\truez{600}%
\put(\value{x},\value{x}){\makebox(0,\value{z})[l]{${#1}$}}%
\put(-\value{x},-\value{z}){\makebox(0,\value{z})[r]{${#2}$}}%
\end{picture}}}\end{picture}}}%

% north-north-west arrow
\newcommand{\nnwar}{\NNWAR{}{}}%

% north-north-west arrow with upper name [1]
\newcommand{\Nnwar}[1]{\NNWAR{#1}{}}%

% north-north-west arrow with lower name [1]
\newcommand{\nnwaR}[1]{\NNWAR{}{#1}}%

% HORIZONTAL CURVED ARROWS

% The following commands produce horizontal curved arrows
% positioned in a box of width 0pt and height 0pt.

% North\East curved arrow with name [1] and length [2]x100 units
\newcommand{\Necurve}[2]%
{\begin{picture}(0,0)%
\truex{1300}\truey{2000}\truez{200}%
\put(0,\value{x}){\oval(#200,\value{y})[t]}%
\put(0,\value{x}){\makebox(0,0){\begin{picture}(#200,0)%
\put(#200,0){\vector(0,-1){\value{z}}}%
\put(0,0){\line(0,-1){\value{z}}}\end{picture}}}%
\truex{2500}%
\put(0,\value{x}){\makebox(0,0)[b]{${#1}$}}%
\end{picture}}%

% North\East curved arrow of length [1]x100 units
\newcommand{\necurve}[1]{\Necurve{}{#1}}%

% North\West curved arrow with name [1] and length [2]x100 units
\newcommand{\Nwcurve}[2]%
{\begin{picture}(0,0)%
\truex{1300}\truey{2000}\truez{200}%
\put(0,\value{x}){\oval(#200,\value{y})[t]}%
\put(0,\value{x}){\makebox(0,0){\begin{picture}(#200,0)%
\put(#200,0){\line(0,-1){\value{z}}}%
\put(0,0){\vector(0,-1){\value{z}}}\end{picture}}}%
\truex{2500}%
\put(0,\value{x}){\makebox(0,0)[b]{${#1}$}}%
\end{picture}}%

% North\West curved arrow of length [1]x100 units
\newcommand{\nwcurve}[1]{\Nwcurve{}{#1}}%

% South\East curved arrow with name [1] and length [2]x100 units
\newcommand{\Securve}[2]%
{\begin{picture}(0,0)%
\truex{1300}\truey{2000}\truez{200}%
\put(0,-\value{x}){\oval(#200,\value{y})[b]}%
\put(0,-\value{x}){\makebox(0,0){\begin{picture}(#200,0)%
\put(#200,0){\vector(0,1){\value{z}}}%
\put(0,0){\line(0,1){\value{z}}}\end{picture}}}%
\truex{2500}%
\put(0,-\value{x}){\makebox(0,0)[t]{${#1}$}}%
\end{picture}}%

% South\East curved arrow of length [1]x100 units
\newcommand{\securve}[1]{\Securve{}{#1}}%

% South\West curved arrow with name [1] and length [2]x100 units
\newcommand{\Swcurve}[2]%
{\begin{picture}(0,0)%
\truex{1300}\truey{2000}\truez{200}%
\put(0,-\value{x}){\oval(#200,\value{y})[b]}%
\put(0,-\value{x}){\makebox(0,0){\begin{picture}(#200,0)%
\put(#200,0){\line(0,1){\value{z}}}%
\put(0,0){\vector(0,1){\value{z}}}\end{picture}}}%
\truex{2500}%
\put(0,-\value{x}){\makebox(0,0)[t]{${#1}$}}%
\end{picture}}%

% South\West curved arrow of length [1]x100 units
\newcommand{\swcurve}[1]{\Swcurve{}{#1}}%

% VERTICAL CURVED ARROWS

% The following commands produce vertical curved arrows
% positioned in a box of width 0pt and height 0pt

% East\South curved arrow with name [1] and length [2]x100 units
\newcommand{\Escurve}[2]%
{\begin{picture}(0,0)%
\truex{1400}\truey{2000}\truez{200}%
\put(\value{x},0){\oval(\value{y},#200)[r]}%
\put(\value{x},0){\makebox(0,0){\begin{picture}(0,#200)%
\put(0,0){\vector(-1,0){\value{z}}}%
\put(0,#200){\line(-1,0){\value{z}}}\end{picture}}}%
\truex{2500}%
\put(\value{x},0){\makebox(0,0)[l]{${#1}$}}%
\end{picture}}%

% East\South curved arrow of length [1]x100 units
\newcommand{\escurve}[1]{\Escurve{}{#1}}%

% East\North curved arrow with name [1] and length [2]x100 units
\newcommand{\Encurve}[2]%
{\begin{picture}(0,0)%
\truex{1400}\truey{2000}\truez{200}%
\put(\value{x},0){\oval(\value{y},#200)[r]}%
\put(\value{x},0){\makebox(0,0){\begin{picture}(0,#200)%
\put(0,0){\line(-1,0){\value{z}}}%
\put(0,#200){\vector(-1,0){\value{z}}}\end{picture}}}%
\truex{2500}%
\put(\value{x},0){\makebox(0,0)[l]{${#1}$}}%
\end{picture}}%

% East\North curved arrow of length [1]x100 units
\newcommand{\encurve}[1]{\Encurve{}{#1}}%

% West\South curved arrow with name [1] and length [2]x100 units
\newcommand{\Wscurve}[2]%
{\begin{picture}(0,0)%
\truex{1300}\truey{2000}\truez{200}%
\put(-\value{x},0){\oval(\value{y},#200)[l]}%
\put(-\value{x},0){\makebox(0,0){\begin{picture}(0,#200)%
\put(0,0){\vector(1,0){\value{z}}}%
\put(0,#200){\line(1,0){\value{z}}}\end{picture}}}%
\truex{2400}%
\put(-\value{x},0){\makebox(0,0)[r]{${#1}$}}%
\end{picture}}%

% West\South curved arrow of length [1]x100 units
\newcommand{\wscurve}[1]{\Wscurve{}{#1}}%

% West\North curved arrow with name [1] and length [2]x100 units
\newcommand{\Wncurve}[2]%
{\begin{picture}(0,0)%
\truex{1300}\truey{2000}\truez{200}%
\put(-\value{x},0){\oval(\value{y},#200)[l]}%
\put(-\value{x},0){\makebox(0,0){\begin{picture}(0,#200)%
\put(0,0){\line(1,0){\value{z}}}%
\put(0,#200){\vector(1,0){\value{z}}}\end{picture}}}%
\truex{2400}%
\put(-\value{x},0){\makebox(0,0)[r]{${#1}$}}%
\end{picture}}%

% West\North curved arrow of length [1]x100 units
\newcommand{\wncurve}[1]{\Wncurve{}{#1}}%

% FREE SLOPE ARROWS

% The following arrows have formal dimensions (0,0)

% arrow with name [1] at position [2,3],
% origin [4,5], slope [6,7] and horizontal extent [8]x100 units,
\newcommand{\Freear}[8]{\begin{picture}(0,0)%
\put(#400,#500){\vector(#6,#7){#800}}%
\truex{#200}\truey{#300}%
\put(\value{x},\value{y}){$#1$}\end{picture}}%

% arrow with origin [1,2], slope [3,4] and horizontal length [5]x100 units
\newcommand{\freear}[5]{\Freear{}{0}{0}{#1}{#2}{#3}{#4}{#5}}%

% INITIALIZATION

\newcount\SCALE%

\newcount\NUMBER%

\newcount\LINE%

\newcount\COLUMN%

\newcount\WIDTH%

\newcount\SOURCE%

\newcount\ARROW%

\newcount\TARGET%

\newcount\ARROWLENGTH%

\newcount\NUMBEROFDOTS%

\newcounter{x}%

\newcounter{y}%

\newcounter{z}%

\newcounter{horizontal}%

\newcounter{vertical}%

\newskip\itemlength%

\newskip\firstitem%

\newskip\seconditem%

\newcommand{\printarrow}{}%

% MACROS FOR DESIGNING DIAGRAMS

% \truex{n} divides nx100 by the scaling factor and puts the result
% in counter x
\newcommand{\truex}[1]{%
\NUMBER=#1%
\multiply\NUMBER by 100%
\divide\NUMBER by \SCALE%
\setcounter{x}{\NUMBER}}%

% \truey{n} divides nx100 by the scaling factor and puts the result
% in counter y
\newcommand{\truey}[1]{%
\NUMBER=#1%
\multiply\NUMBER by 100%
\divide\NUMBER by \SCALE%
\setcounter{y}{\NUMBER}}%

% \truez{n} divides nx100 by the scaling factor and puts the result
% in counter z
\newcommand{\truez}[1]{%
\NUMBER=#1%
\multiply\NUMBER by 100%
\divide\NUMBER by \SCALE%
\setcounter{z}{\NUMBER}}%

% \changecounters computes the values of the various parameters required
% to design an arrow with adjusted length.
\newcommand{\changecounters}[1]{%
\SOURCE=\ARROW%
\ARROW=\TARGET%
\settowidth{\itemlength}{#1}%
\ifdim \itemlength > 2800\unitlength%
\addtolength{\itemlength}{-2800\unitlength}%
\TARGET=\itemlength%
\divide\TARGET by 1310%
\multiply\TARGET by 100%
\divide\TARGET by \SCALE%
\else%
\TARGET=0%
\fi%
\ARROWLENGTH=5000%
\advance\ARROWLENGTH by -\SOURCE%
\advance\ARROWLENGTH by -\TARGET%
\advance\SOURCE by -\TARGET}%

% \initialize initializes the counters required to produce the diagram
% and defines the formal dimensions of the diagram to be (0,0).
\newcommand{\initialize}[1]{%
\LINE=0%
\COLUMN=0%
\WIDTH=0%
\ARROW=0%
\TARGET=0%
\changecounters{#1}%
\renewcommand{\printarrow}{#1}%
\begin{center}%
\vspace{10pt}%
\begin{picture}(0,0)}%

% \DIAG starts the construction of an unscaled diagram
\newcommand{\DIAG}[1]{%
\SCALE=100%
\setlength{\unitlength}{655sp}%
\initialize{\mbox{$#1$}}}%

% \DIAGV{n} starts the construction of a diagram scaled by a factor
% n percent, computes the scaled unit length and the unscaling factor.
\newcommand{\DIAGV}[2]{%
\SCALE=#1%
\setlength{\unitlength}{655sp}%
\multiply\unitlength by \SCALE%
\divide\unitlength by 100%
\initialize{\mbox{$#2$}}}%

% \n introduces the next item of the diagram, computes its parameters
% and prints the previous item.
\newcommand{\n}[1]{%
\changecounters{\mbox{$#1$}}%
\put(\COLUMN,\LINE){\makebox(0,0){\printarrow}}%
\thinlines%
\renewcommand{\printarrow}{\mbox{$#1$}}%
\advance\COLUMN by 4000}%

% \nn prints the  last item of a line and starts a new line.
\newcommand{\nn}[1]{%
\put(\COLUMN,\LINE){\makebox(0,0){\printarrow}}%
\thinlines%
\ifnum \WIDTH < \COLUMN%
\WIDTH=\COLUMN%
\else%
\fi%
\advance\LINE by -4000%
\COLUMN=0%
\ARROW=0%
\TARGET=0%
\changecounters{\mbox{$#1$}}%
\renewcommand{\printarrow}{\mbox{$#1$}}}%

% \conclude prints the last item of the diagram and sets the
% dimensions of the diagram;
\newcommand{\conclude}{%
\put(\COLUMN,\LINE){\makebox(0,0){\printarrow}}%
\thinlines%
\ifnum \WIDTH < \COLUMN%
\WIDTH=\COLUMN%
\else%
\fi%
\setcounter{horizontal}{\WIDTH}%
\setcounter{vertical}{-\LINE}%
\end{picture}}%

% \diag prints the last item of the diagram and takes care of the spacing
\newcommand{\diag}{%
\conclude%
\raisebox{0pt}[0pt][\value{vertical}\unitlength]{}%
\hspace*{\value{horizontal}\unitlength}%
\vspace{10pt}%
\end{center}%
\setlength{\unitlength}{1pt}}%

% \diagv{t}{l}{b} prints the last item of the diagram and
% adds a t points extra space at the top of the diagram
% adds a l points extra space at the left of the diagram
% adds a b points extra space at the bottom of the diagram
\newcommand{\diagv}[3]{%
\conclude%
\NUMBER=#1%
\rule{0pt}{\NUMBER pt}%
\hspace*{-#2pt}%
\raisebox{0pt}[0pt][\value{vertical}\unitlength]{}%
\hspace*{\value{horizontal}\unitlength}%6
\NUMBER=#3%
\advance\NUMBER by 10%
\vspace*{\NUMBER pt}%
\end{center}%
\setlength{\unitlength}{1pt}}%

% \N normalizes the height of a vertex
\newcommand{\N}[1]%
{\raisebox{0pt}[7pt][0pt]{$#1$}}%

% \movename{f}{n}{m} moves the name f of the arrow n points right
% and m points up.
\newcommand{\movename}[3]{%
\hspace{#2pt}%
\raisebox{#3pt}[5pt][2pt]{\raisebox{#3pt}{$#1$}}%
\hspace{-#2pt}}%

% \movearrow{\...}{n}{m} moves arrow \... n points right
% and m points up.
\newcommand{\movearrow}[3]{%
\makebox[0pt]{%
\hspace{#2pt}\hspace{#2pt}%
\raisebox{#3pt}[0pt][0pt]{\raisebox{#3pt}{$#1$}}}}%

% \movevertex{A}{n}{m} moves vertex A n points right and m points up
\newcommand{\movevertex}[3]{%
\mbox{\hspace{#2pt}%
\raisebox{#3pt}{\raisebox{#3pt}{$#1$}}%
\hspace{-#2pt}}}%

% \crosslength{P}{Q} computes the formal dimensions of the
% superposition of pictures P and Q
\newcommand{\crosslength}[2]{%
\settowidth{\firstitem}{#1}%
\settowidth{\seconditem}{#2}%
\ifdim\firstitem < \seconditem%
\itemlength=\seconditem%
\else%
\itemlength=\firstitem%
\fi%
\divide\itemlength by 2%
\hspace{\itemlength}}%

% \cross{P}{Q} superposes pictures P and Q
\newcommand{\cross}[2]{%
\crosslength{\mbox{$#1$}}{\mbox{$#2$}}%
\begin{picture}(0,0)%
\put(0,0){\makebox(0,0){$#1$}}%
\thinlines%
\put(0,0){\makebox(0,0){$#2$}}%
\thinlines%
\end{picture}%
\crosslength{\mbox{$#1$}}{\mbox{$#2$}}}%

% \B prints the next arrow in bold-face type
\newcommand{\B}{\thicklines}

% SPECIAL SYMBOLS

% adjoint symbol
\newcommand{\adj}{\begin{picture}(9,6)%
\put(1,3){\line(1,0){6}}\put(7,0){\line(0,1){6}}%
\end{picture}}%

% commutative diagram symbol
\newcommand{\com}{\begin{picture}(12,8)%
\put(6,4){\oval(8,8)[b]}\put(6,4){\oval(8,8)[r]}%
\put(6,8){\vector(-1,0){2}}\end{picture}}%

% natural transformation with names [1],[2],[3]
\newcommand{\Nat}[3]{\raisebox{-2pt}%
{\begin{picture}(34,15)%
\put(2,10){\vector(1,0){30}}%
\put(2,0){\vector(1,0){30}}%
\put(13,2){$\Downarrow$}%
\put(20,3){$\scriptstyle{#2}$}%
\put(4,11){$\scriptstyle{#1}$}%
\put(4,1){$\scriptstyle{#3}$}%
\end{picture}}}%

% natural transformation
\newcommand{\nat}{\raisebox{-2pt}%
{\begin{picture}(34,10)%
\put(2,10){\vector(1,0){30}}%
\put(2,0){\vector(1,0){30}}%
\put(13,2){$\Downarrow$}%
\end{picture}}}%

% pair of natural transformations with names [1],[2],[3],[4],[5]
\newcommand{\Binat}[5]{\raisebox{-7.5pt}%
{\begin{picture}(34,25)%
\put(2,20){\vector(1,0){30}}%
\put(2,10){\vector(1,0){30}}%
\put(2,0){\vector(1,0){30}}%
\put(13,12){$\Downarrow$}%
\put(13,2){$\Downarrow$}%
\put(20,13){$\scriptstyle{#2}$}%
\put(20,3){$\scriptstyle{#4}$}%
\put(4,21){$\scriptstyle{#1}$}%
\put(4,11){$\scriptstyle{#3}$}%
\put(4,1){$\scriptstyle{#5}$}%
\end{picture}}}%

% pair of natural transformations
\newcommand{\binat}{\raisebox{-7.5pt}%
{\begin{picture}(34,20)%
\put(2,20){\vector(1,0){30}}%
\put(2,10){\vector(1,0){30}}%
\put(2,0){\vector(1,0){30}}%
\put(13,12){$\Downarrow$}%
\put(13,2){$\Downarrow$}%
\end{picture}}}%